\documentclass[a4paper,german,final,oneside,notitlepage,12pt]{article}

\usepackage[T1]{fontenc}
\usepackage{amssymb}
\usepackage{amsthm}
\usepackage{amsmath}
\usepackage[all]{xy}
\usepackage{graphicx}
\usepackage{mathrsfs}
\usepackage{afterpage}
\usepackage{float}
\usepackage[bottom]{footmisc}
\usepackage[margin=2cm,font=normalsize,labelfont=bf]{caption}
\usepackage[normalem]{ulem}
\usepackage[section]{placeins}
\usepackage{verbatim}
\usepackage{amstext}

\makeatletter

\gdef\@ntitle{\@title}
\def\adress#1{\gdef\@adress{#1}}
\def\@adress{}
\def\preprint#1{\gdef\@preprint{#1}}
\def\@preprint{}
\def\keywords#1{\gdef\@keywords{#1}}
\def\@keywords{}

\def\refname{References}

\def\href#1#2{#2}

\def\kohyp{
  \usepackage{hyperref}
  \hypersetup{
    linktocpage = true,
    pdftitle = {\@title},
    pdfauthor = {\@author},
    pdfkeywords = {\@keywords},    
    bookmarksopen = true,
    bookmarksopenlevel = 1
  }}  
\def\showkeywords{\footnotesize\\\\ \textbf{Keywords}: \@keywords.}

\ifnum \c@secnumdepth >0
  \newtheorem{theorem}{Theorem}[section]
\else
  \newtheorem{theorem}{Theorem}
\fi
\newtheorem{definition}[theorem]{Definition}
\newtheorem{rem}[theorem]{Remark}
\newtheorem{corollary}[theorem]{Corollary}
\newtheorem{exa}[theorem]{Example}
\newtheorem{proposition}[theorem]{Proposition}
\newtheorem{lemma}[theorem]{Lemma}

\newenvironment{remark}{\begin{rem}\normalfont}{\end{rem}}
\newenvironment{example}{\begin{exa}\normalfont}{\end{exa}}

\def\N {\mathbb{N}}

\def\R {\mathbb{R}}

\def\id{\mathrm{id}}

\renewcommand{\varepsilon}{\epsilon}

\newcommand{\alxy}[1]{\begin{aligned}\xymatrix{#1}\end{aligned}}
\newcommand{\alxydim}[2]{\begin{aligned}\xymatrix#1{#2}\end{aligned}}

\renewenvironment{proof}{Proof.\ }{\hfill{$\square$}\\}

\newcommand\erf[1]{(\ref{#1})}

\ifnum \c@secnumdepth >0
  \numberwithin{theorem}{section}
  \numberwithin{equation}{section}
  \renewcommand\theequation{\thesection.\arabic{equation}}
\fi

\renewcommand{\emph}[1]{\def\reserved@a{it}\ifx\f@shape\reserved@a\uline{#1}\else\textit{#1}\fi}

\newcommand{\trivlin}{\hbox{$1\hskip -1.2pt\vrule depth 0pt height 1.6ex width 0.7pt \vrule depth 0pt height 0.3pt width 0.12em$}}

\def\tocsection#1{
  \section*{#1}
  \addcontentsline{toc}{section}{#1}}

\makeatother

\renewcommand{\to}{\!\xymatrix@C=0.5cm{\ar[r] &}}
\renewcommand{\mapsto}{\xymatrix@C=0.5cm{\ar@{|->}[r] &}}
\renewcommand{\Rightarrow}{\xymatrix@C=0.5cm{\ar@{=>}[r] &}}
\newcommand{\incl}{\xymatrix@C=0.5cm{\ar@{^(->}[r] &}}

\makeatletter

\def\kobibarxiv#1{\newline\href{http://arxiv.org/abs/#1}{\texttt{arxiv:#1}}}
\def\kobibarxivsection#1{\texttt{[#1]}}
\def\kobiblink#1{\newline\href{#1}{#1}}
\def\nolinks{\def\kobiblink##1{}}
\def\noarxiv{\def\kobibarxiv##1{}\def\kobibarxivsection##1{}}

\nolinks
\noarxiv

\newcommand{\tql}{\textquotedblleft{}}
\newcommand{\tqr}{\textquotedblright}

\def\mytitle{}
\def\zmptitle{
  \begin{tabular}{cc}
    \begin{minipage}[c]{0.4\textwidth}
      \begin{flushleft}
        \includegraphics[width=110pt]{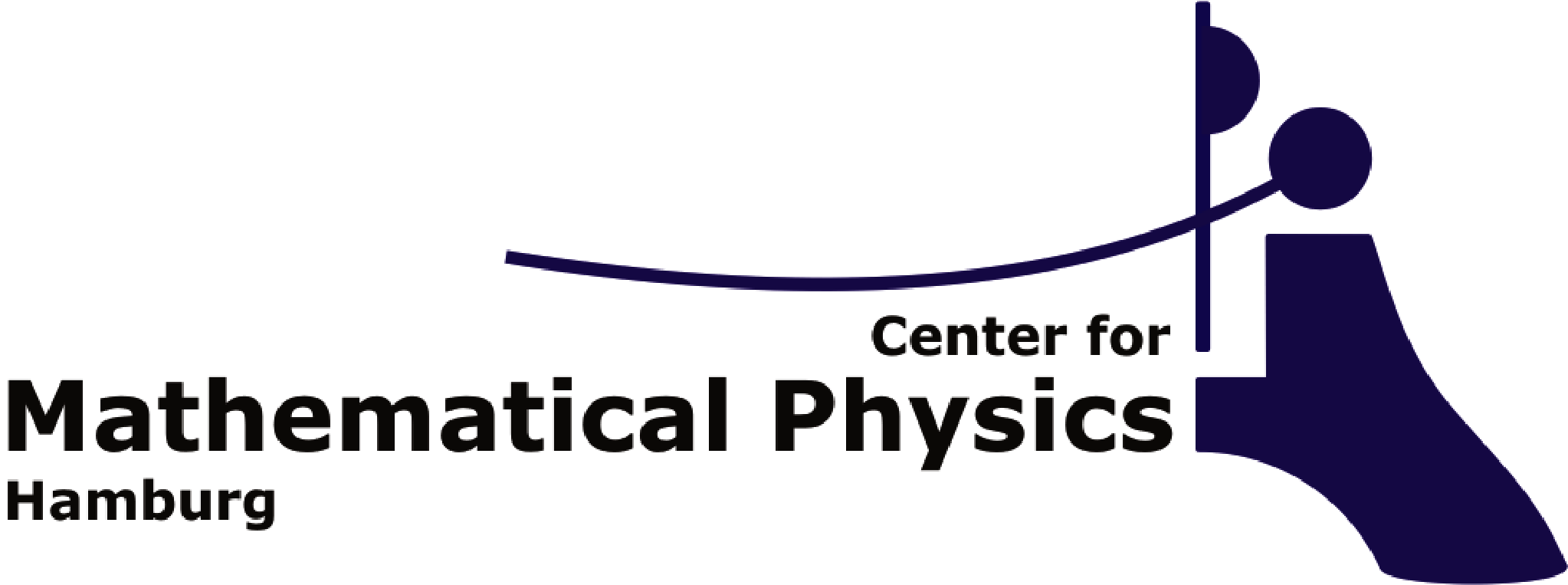}
      \end{flushleft}  
    \end{minipage}&
    \begin{minipage}[c]{0.55\textwidth}
      \begin{flushright}
      {\small\sf\@preprint}
      \end{flushright}
    \end{minipage}
  \end{tabular}
  \vskip 2cm}

\def\maketitle{
  \newpage
  \noindent
  \mytitle
  \begin{center}
    \LARGE\@ntitle
    \if!\@author!\else \vskip 0.5cm \large\@author\fi
    \if!\@adress!\else \vskip 0.5cm \normalsize\@adress\fi
  \end{center}
  \vskip 2cm\thispagestyle{empty}}

\newlength{\bibitemwidth}
\newenvironment{kobibitem}{
  \setlength{\bibitemwidth}{\textwidth}
  \addtolength{\bibitemwidth}{-\labelwidth}
  \begin{minipage}[t]{\bibitemwidth}
}{\end{minipage}}

\def\kobib{

}

\newif\if@fewtab\@fewtabtrue{
  \count255=\time\divide\count255 by 60
  \xdef\hourmin{\number\count255}
  \multiply\count255 by-60\advance\count255 by\time
  \xdef\hourmin{\hourmin:\ifnum\count255<10 0\fi\the\count255}}
\def\ps@draft{
  \let\@mkboth\@gobbletwo
  \def\@oddfoot{
    \hbox to 7 cm{\tiny \versionno\hfil}
    \hskip -7cm\hfil\rm\thepage\hfil{\tiny\draftdate}}
  \def\@oddhead{}
  \def\@evenhead{}
  \let\@evenfoot\@oddfoot}
\def\draftdate{\number\month/\number\day/\number\year\ \ \ \hourmin }
\newcommand\version[1]{
  \typeout{}\typeout{#1}\typeout{}
  \vskip-1.7cm \centerline{\fbox{{\normalsize\tt DRAFT -- #1 -- }
  {\small\draftdate}}} \vskip0.92cm}
\def\draft#1{
  \def\versionno{#1}
  \pagestyle{draft}\thispagestyle{draft}
  \gdef\@ntitle{\version\versionno \@title}
  \global\def\draftcontrol{1}}
\global\def\draftcontrol{0}
\def\todo#1{\ifnum\draftcontrol>0 \\\bigskip {\tt #1} \\ \fi}

\makeatother

\def\proof{{Proof. }}
\def\endofproof{\hfill{$\square$}\\}

\newcommand{\diffco}[3]{Z^{#2}_{#3}(#1)^{\infty}}

\newcommand{\smsp}{D^{\infty}}

\newcommand{\sm}{C^{\infty}}
\newcommand{\fu}{\mathcal{P}}
\newcommand{\fo}{\mathcal{D}}

\newcommand{\bigon}[5]{\alxydim{@C=1.2cm}{#1 \ar@/^1.5pc/[r]^{#3}="1" \ar@/_1.5pc/[r]_{#4}="2" \ar@{=>}"1";"2"|{#5} & #2}}

\newcommand{\quadrat}[9]{\alxydim{@C=1.2cm@R=1.2cm}{#1 \ar[r]^{#5} \ar[d]_{#6} & #2 \ar@{=>}[dl]|{#9} \ar[d]^{#7} \\ #3 \ar[r]_{#8} & #4}}

\newcommand{\dreieck}[7]{\alxydim{@R=1.2cm@C=0.3cm}{& #2  \ar[dr]^{#5} & \\ #1 \ar[ur]^{#4} \ar[rr]_{#6}="1" && #3 \ar@{=>}[ul];"1"|{#7}}
}

\newcommand{\kopfdreieck}[7]{\alxydim{@R=1.2cm@C=0.3cm}{#1 \ar[dr]_{#4} \ar[rr]^{#6}="1" && #3 \ar@{=>}"1";[dl]|{#7} \\ & #2  \ar[ur]_{#5} & }
}

\title{Smooth Functors vs. Differential Forms}
\author{Urs Schreiber and Konrad Waldorf}
\adress{Department Mathematik\\Schwerpunkt Algebra und Zahlentheorie\\Universität
Hamburg\\Bundesstra\ss e 55\\D--20146 Hamburg}
\preprint{Hamb. Beitr. Math. Nr. 297\\ZMP-HH/08-04}

\def\mytitle{\zmptitle}
\usepackage{amstext}
\usepackage{hyperref}

\begin{document}


\maketitle

\def\der{derivative}
\def\Der{Derivative}

\begin{abstract}
We establish a relation between smooth 2-functors defined on the path 2-groupoid of a smooth manifold and   differential forms on this manifold.   This relation can be understood as a part of a dictionary between fundamental notions from category theory and  differential geometry. We show that smooth 2-functors appear in several fields, namely as connections on  (non-abelian) gerbes,   as  \der s of smooth functors and as critical points in BF theory. We demonstrate further that our dictionary provides a powerful tool to discuss the transgression of geometric objects to loop spaces. 
\end{abstract}

\newpage

\tableofcontents
\section*{Introduction}
\addcontentsline{toc}{section}{Introduction}

The present article is the second of three articles aiming at a general and systematical approach to connections on (non-abelian) gerbes and their surface holonomy.

In the first article \cite{schreiber3} we have established an equivalence between  categories of  fibre bundles with connection over a smooth manifold $X$, and categories of certain functors, called transport functors. Let us spell out this equivalence, reduced to  \emph{trivial} principal $G$-bundles with connection. These are just $\mathfrak{g}$-valued 1-forms on $X$, where $\mathfrak{g}$ is the Lie algebra of the Lie group  $G$. Our equivalence is then a bijection 
\begin{equation*}
\Omega^1(X,\mathfrak{g}) \cong \left\lbrace \alxydim{}{\txt{Smooth functors \\ $\mathcal{P}_1(X) \rightarrow \mathcal{B} G$}} \right \rbrace
\end{equation*}
between the set of $\mathfrak{g}$-valued 1-forms on $X$ and the set of smooth functors between two groupoids $\mathcal{P}_1(X)$ and $\mathcal{B}G$. One the one hand, we have the path groupoid $\mathcal{P}_1(X)$ which is associated to the manifold $X$: its objects are the points of $X$, and the morphisms between two points are (thin homotopy classes of smooth) paths between these two points. On the other hand, we have a groupoid $\mathcal{B}G$ which is associated to the Lie group $G$: it has just one object and every group element acts as an automorphisms of this object. The notation $\mathcal{B}G$ is devoted to the fact that the geometric realization of its nerve is the classifying space $BG$ of the group $G$. 

Now, the  functors $F:\mathcal{P}_1(X) \to \mathcal{B}G$ we have on the right hand side of the above bijection assign group elements $F(\gamma)$ to  paths $\gamma$ in $X$; this assignment  is smooth in a sense that can be expressed in terms of smooth maps between smooth manifolds. For the convenience of the reader, we review this relation between  smooth functors and differential forms  in Section \ref{sec2}.

In the present article we generalize the above bijection between smooth functors and 1-forms to smooth 2-functors and 2-forms. The aim of this generalization is multiple, but for a start we want to  give  an impression how the generalized bijection looks like. The first step is the generalization of the    categories $\mathcal{P}_1(X)$ and $\mathcal{B}G$ to appropriate 2-categories. On the one hand, we introduce the path  2-groupoid $\mathcal{P}_2(X)$ of a smooth manifold $X$ by adding 2-morphisms to the path groupoid $\mathcal{P}_1(X)$. These 2-morphisms are (thin homotopy classes of) smooth homotopies between  smooth paths in $X$. On the other hand,  we infer that the   group $G$ which was present before has to be replaced by a (strict)   2-group: basically, this is a group object in  categories, i.e. a category with additional structure. The concept of  2-groups can  be refined to \emph{Lie} 2-groups; such a Lie 2-group $\mathfrak{G}$ underlies the generalized relation we are after. We form a 2-category $\mathcal{B}\mathfrak{G}$ mimicking the same idea we used for the category $\mathcal{B}G$: it has just  one object and the Hom-category of this object is the category  $\mathfrak{G}$. 

Equipped with these generalized  2-categories, we consider 2-functors
\begin{equation*}
F: \mathcal{P}_2(X) \to \mathcal{B}\mathfrak{G}\text{,}
\end{equation*}
and it even makes perfectly sense to qualify some as  \emph{smooth} 2-functors. As before, the smoothness can  be expressed in terms of smooth maps between smooth manifolds (Definition \ref{def7}). 

In order to explore which kind of differential form corresponds to a smooth 2-functor $F: \mathcal{P}_2(X) \to \mathcal{B}\mathfrak{G}$, we put the abstract concept of a  2-group in a more familiar setting. According to  Brown and Spencer \cite{brown1}, a  2-group $\mathfrak{G}$ is equivalent to a crossed module: a structure introduced by Whitehead \cite{whitehead1}  consisting of two  ordinary groups $G$ and $H$, a  group homomorphism $t:H \to G$ and a compatible action of $G$ on $H$. Similarly, a \emph{Lie} 2-group corresponds to \emph{Lie} groups $G$ and $H$ and  \emph{smooth} additional structure. We denote the Lie algebras of the two Lie groups $G$ and $H$ by $\mathfrak{g}$ and $\mathfrak{h}$, respectively. The first  result
of this article (Proposition \ref{prop1}) is that the smooth 2-functor $F$ defines a $\mathfrak{g}$-valued 1-form $A$ \emph{and} an $\mathfrak{h}$-valued 2-form $B$ on $X$ that are related by
\begin{equation*}
\mathrm{d}A + [A \wedge A] = t_{*} (B)\text{.}
\end{equation*}
Here, we have the ordinary curvature 2-form of $A$ on the left hand side, and $t_{*}$ is the Lie algebra homomorphism induced by $t$. 
The two differential forms $A$ and $B$ contain in fact all information about the 2-functor $F$ they came from: we describe an explicit procedure how to integrate two forms $A$ and $B$ (satisfying the above condition) to a smooth 2-functor $F: \mathcal{P}_2(X) \to \mathcal{B}\mathfrak{G}$. This integration involves iterated solutions of ordinary differential equations.
The main result of this article (Theorem \ref{th2}) is that we obtain a bijection
\begin{equation*}
\alxydim{}{\left\lbrace\txt{Smooth 2-functors\\$F:\mathcal{P}_2(X) \rightarrow \mathcal{B}\mathfrak{G}$} \right\rbrace} \cong \alxydim{}{\left \lbrace \txt{ $(A,B)\in \Omega^1(X,\mathfrak{g})\times\Omega^2(X,\mathfrak{h})$\\ with $\mathrm{d}A + [A \wedge A] = t_{*} (B) $} \right \rbrace\text{.}}
\end{equation*}
This is the announced generalization of the relation between smooth 1-functors and differential 1-forms from \cite{schreiber3}. Besides, we also explore the geometric structure that corresponds to morphisms (pseudonatural transformations)  and 2-morphisms (modifications) between 2-functors. The derivation of all the relations that are imposed on this structure takes a large part of this article, and  is collected in  Section  \ref{sec3}. 

\medskip

In Section \ref{sec4} we try to convince the reader that smooth 2-functors are implicitly  present in various fields, and we describe the impact of our new bijection to these fields. We give three examples. 
The first example are connections on (possibly non-abelian) gerbes. As mentioned at the beginning of this introduction, ordinary smooth functors correspond to connections on \emph{trivial} principal bundles. We claim here that smooth \emph{2-}functors correspond in the same way to connections on \emph{trivial gerbes}. 

For abelian gerbes, such connections have been studied by Brylinski on sheaves of groupoids \cite{brylinski1}, and by Murray on bundle gerbes \cite{murray}. In both cases, our claim shows to be true. Similar results for abelian gerbes have been obtained in \cite{mackaay1}. Connections on a certain class of (possibly) non-abelian gerbes have been introduced by Breen and Messing \cite{breen1}. Their connection -- considered on a trivial gerbe -- is a pair $(A,B)$ of a 1-form and a 2-form like they arise  from a smooth 2-functor in the way outlined above. Interestingly, the two forms of a Breen-Messing connection are not necessarily related to each other, in contrast to the two forms coming from a smooth 2-functor. We argue that this difference is  related to the  question,  whether such connections induce a  notion of  \emph{surface holonomy}. 

Moreover, certain higher gauge theories can be described by pairs $(A,B)$ of differential forms with values in the Lie algebras belonging to the two Lie groups of a crossed module, and even the relation between $A$ and $B$ we found here is already present in this context \cite{girelli1}.  Since higher gauge theories are naturally related to connections on gerbes, a further link between smooth 2-functors and connections on gerbes is present. 

A deeper discussion of connections on  \emph{non-}trivial  gerbes and their surface holonomy is the content of the third part \cite{schreiber2} in our series of articles.

The second example  of smooth 2-functors we want to give are \der\  2-functors. For any Lie group $G$, a smooth functor $F:\mathcal{P}_1(X) \to \mathcal{B}G$ determines a smooth 2-functor 
\begin{equation*}
\mathrm{d}F: \mathcal{P}_2(X) \to \mathcal{BE}G\text{,}
\end{equation*}
where $\mathcal{E}G$ is a  Lie 2-group  introduced by Segal as a model for the universal $G$-bundle $EG$ \cite{segal3}. Using our dictionary between smooth functors and differential forms for both $F$ and $\mathrm{d}F$ we see that the functor $F$ corresponds to a trivial principal $G$-bundle with connection $\omega$ over $X$, while  the \der\ 2-functor $\mathrm{d}F$ induces a 2-form $B \in \Omega^2(X,\mathfrak{g})$. We prove that $B$ is the curvature of the connection $\omega$, so  that the relation between $F$ and $\mathrm{d}F$ implies a  relation between the holonomy of $\omega$ and its curvature. We show that this  establishes a new proof of the so-called non-abelian Stokes' Theorem (Theorem \ref{th1}).  

\medskip

The third example where smooth 2-functors arise is a certain topological field theory which is called BF-theory due to the presence of two fields $B$ and $F$. These fields are 2-forms with values in the Lie algebra $\mathfrak{g}$ of a Lie group $G$; actually $F=\mathrm{d}A + [A \wedge A]$ is the curvature of a 1-form $A$. We prove that the critical points of the  BF action functional are those pairs $(A,B)$ satisfying the relation described above. In other words, smooth 2-functors arise as the classical solutions of the field equations of BF-theory (Proposition \ref{prop8}). 

\medskip

Section \ref{sec5} is devoted to the following observation: every element in the loop space $LX$ of a smooth manifold $X$ can be understood as a particular morphism in the path groupoid $\mathcal{P}_1(X)$, and also  as a particular 1-morphism in the path 2-groupoid $\mathcal{P}_2(X)$. This way, functors on the path groupoid, and 2-functors on the path 2-groupoid are intrinsically related to structure on the loop space of $X$. 

First we observe that the structure on the loop space which is induced by a smooth functor $F:\mathcal{P}_1(X) \to \mathcal{B}G$ is a smooth function $LX \to G$, and that this function is nothing but the holonomy of the (trivial) principal $G$-bundle with connection associated to $F$. Then we prove that the structure which is induced by a smooth 2-functor $F:\mathcal{P}_2(X) \to \mathcal{B}\mathfrak{G}$  is a smooth functor \begin{equation*}
\mathcal{P}_1(LX) \to \Lambda\mathcal{B}\mathfrak{G}\text{,}
\end{equation*}
where $\Lambda\mathcal{B}\mathfrak{G}$ is a certain  category constructed from the 2-groupoid $\mathcal{B}\mathfrak{G}$. 
In order to be able to speak about smooth functors on the loop space, we work with the canonical diffeology on $LX$: this is a structure which generalizes a smooth manifold structure and  is more suitable for mapping spaces. We extend the relation between functors and 1-forms from \cite{schreiber3} to  diffeological spaces (Theorem \ref{th4}), and prove that the above smooth functor on the path groupoid of the loop space corresponds to the following structure: a smooth function $LX \to G$, a 1-form $A_F\in\Omega^1(LX,\mathfrak{g})$ and a 1-form $\varphi_F\in\Omega^1(LX,\mathfrak{h})$.

Denoting by $(A,B)$ the differential forms that belong to the smooth 2-functor $F$ we started with, we derive an explicit relation between the differential forms $(A,B)$ on $X$ and $(A_F,\varphi_F)$ on $LX$ (Proposition \ref{prop7}). This relation involves integration along the fibre, and admits an outlook on the question, what the transgression of a non-abelian gerbe over $X$ to the loop space $LX$ is.

Finally, we have included an Appendix in which we review basic notions from 2-category theory and important definitions and examples related to Lie 2-groups and smooth crossed modules. For the convenience of the reader,  there is also a table of notation.

\section{Review: Smooth Functors and 1-Forms}

\label{sec2}

In this section we review some relevant definitions and results from \cite{schreiber3} and references therein.

\subsection{The Path Groupoid of a Smooth Manifold}

\label{sec1_1}

In the topological category, the basic idea of the path groupoid is very simple: for a topological space $X$, it is the category whose objects are the points of $X$, and whose morphisms are homotopy classes of  continuous paths in~$X$. For smooth manifolds, one considers \emph{smooth paths}: these are 
smooth maps $\gamma:[0,1] \to X$ with sitting instants, i.e.  a number $0 < \varepsilon <\frac{1}{2}$ with
$\gamma(t)=\gamma(0)$ for $0 \leq t < \varepsilon$ and $\gamma(t)=\gamma(1)$ for $1-\varepsilon
< t \leq 1$. The set of smooth paths in $X$ is denoted by $PX$. \label{not:smoothpaths} The sitting instants assure that two  paths $\gamma:x \to y$
and $\gamma':y \to z$ can be composed   to a new  path $\gamma' \circ \gamma:x \to
z$. However, the composition of paths is not associative, so that a category can only be defined using certain quotients of $PX$ as its morphisms. 

There are essentially  three ways to define such quotients. The first  is to take reparameterization classes  $P^{0}X := PX / \sim_{0}$,
where   $\gamma_1 \sim_0 \gamma_2$  if there exists an orientation-preserving diffeomorphism $\varphi$ of $[0,1]$ such
that $\gamma_2 = \gamma_1 \circ \varphi$.  The second  way is to take thin homotopy classes, $P^1X:=PX / \sim_{1}$: \label{not:thinhomotopyclasses} 
\begin{definition}
\label{def3}
Two paths $\gamma_1,\gamma_2:x \to y$ are called \emph{thin homotopy equi\-va\-lent}, denoted $\gamma_1 \sim_1 \gamma_2$,
if there exists a smooth map 
$h: [0,1]^2 \to M$
such that 
\begin{enumerate}
\item[(1)]
it has sitting instants: there exists a number $0 < \epsilon < \frac{1}{2}$ with 
\begin{enumerate}
\item[a)] 
$h(s,t)=x$ for $0 \leq t <\epsilon$ and $h(s,t)=y$ for $1-\epsilon< t
\leq 1$.
\item[b)]
$h(s,t)=\gamma_1(t)$ for $0 \leq s < \epsilon$ and $h(s,t)=\gamma_2(t)$ for
$1-\epsilon<s\leq 1$. 
\end{enumerate}

\item[(2)]
the differential of $h$ has at most rank 1. 
\end{enumerate}
\end{definition}

The third way is to take homotopy classes $P^2X := PX/ \sim_2$ just like in  Definition \ref{def3} but without condition (2). Notice that there are projections
\begin{equation}
\label{34}
\alxy{PX \ar[r] & P^{0}X \ar[r] & P^1X \ar[r] & P^2X}
\end{equation}
and that the above-mentioned composition of paths induces well-defined compositions on all $P^{i}X$. We denote  by $\id_x$ the constant path at a point $x$. In $P^{0}X$ we have 
\begin{equation}
\label{32}
\gamma \circ \id_x \;\sim_0\; \gamma\; \sim_0\; \id_y \circ \gamma
\quad\text{ and }\quad
(\gamma_3 \circ \gamma_2) \circ \gamma_1 \;\sim_0 \;\gamma_3 \circ (\gamma_2 \circ \gamma_1)\text{;}
\end{equation}
these are the axioms of a category with objects $X$ and morphisms $P^0X$. We further  denote by $\gamma^{-1}:y \to x$ the path $\gamma^{-1}(t):=
\gamma(1-t)$. In $P^1X$ we have additionally to \erf{32}
\begin{equation*}
\gamma^{-1}
\circ \gamma \;\sim_1\; \id_x
\end{equation*}
for any path $\gamma:x \to y$, so that the corresponding
category with morphisms $P^{1}X$ is even a groupoid. This groupoid is denoted by $\mathcal{P}_1(X)$  and called the
\emph{path groupoid}\label{not:pathgroupoid} of $X$. The groupoid $\Pi_1(X)$ with morphisms $P^2X$ is well-known as the \emph{fundamental groupoid} of the smooth manifold $X$. All these categories are compatible with smooth maps between smooth manifolds in the sense that any smooth map $f:X \to Y$ induces maps
$f_{*}: P^iX \to P^iY$, and that these maps  furnish functors between the respective categories.

\begin{remark}
\normalfont
The groupoids $\mathcal{P}_1(X)$ and $\Pi_1(X)$ are important for parallel transport in a fibre bundle with connection over $X$ in the sense that any such bundle defines a  functor
\begin{equation*}
\mathrm{tra}:\mathcal{P}_1(X) \to T\text{,}
\end{equation*}
where $T$ is a category in which the fibres of the bundle are objects.  If the connection is flat, this functor factors through the fundamental groupoid $\Pi_1(X)$. More on the relation between functors and connections in fibre bundles can be found in Section \ref{sec4_1} and in \cite{schreiber3}.
\end{remark}

\subsection{Diffeological Spaces}

\label{sec1_2}

A \emph{Lie category} is a category $S$ whose sets $S_{0}$ of objects and $S_1$ of morphisms are smooth manifolds, whose target and source maps are surjective submersions, whose identity map is an embedding and whose composition is smooth.   A functor $F:S \to T$ between Lie categories $S$ and $T$ is called \emph{smooth}, if its assignments $F_0:S_0 \to T_0$ and $F_1:S_1 \to T_1$ are smooth maps. The path groupoid $\mathcal{P}_1(X)$ of a smooth manifold is, however, \emph{not} a Lie category, since its set of morphisms $P^{1}X$ has not the structure of a smooth manifold.  In \cite{schreiber3} we have instead equipped  $P^{1}X$  with a \emph{diffeology}, a structure that generalizes a smooth manifold structure \cite{chen1,souriau1}. This diffeology on $P^1 X$ has  been introduced  in \cite{caetano}. 
 For the convenience of the reader let us recall the basic definitions (see also Appendix A.2 of \cite{schreiber3}). 

\begin{definition}
\label{def2}
A \emph{diffeological space} is a set $X$ together with a collection of  plots: maps
\begin{equation*}
c : U \to X
\end{equation*}
each of them defined on an open subset $U \subset \R^{k}$ for any $k\in
\N_0$, such that three axioms are satisfied:
\begin{itemize}
\item[(D1)]
for any plot $c: U \to X$ and any smooth function $f:V \to U$ also 
$c \circ f$ is a plot. 
\item[(D2)]
every constant map $c:U \to X$ is a plot.
\item[(D3)]
if $f:U \to X$ is a map defined on $U \subset \R^k$ and  $\lbrace U_i \rbrace_{i\in I}$ is an open cover of $U$ for which all restrictions $f|_{U_{i}}$ are plots of $X$, then also $f$ is a plot.

\end{itemize}
A \emph{diffeological map} between diffeological spaces  $X$ and $Y$ is a map $f:X \to Y$ such that
for every plot $c: U \to X$ of $X$ the map $f \circ c:U \to Y$ is a plot of $Y$. The set of all diffeological maps is denoted by $\smsp(X,Y)$. 
\end{definition}

Any smooth manifold  is a diffeological space, whose plots are all smooth maps defined on all open subsets of $\R^k$, for all $k$.  If $M$ and $N$ are smooth manifolds, a map $f:M \to N$ is diffeological if and only if it is smooth. In other words,  diffeological spaces and maps form a category $\smsp$ \label{not:dinfty} that contains the category $\sm$ of smooth manifolds (without boundary) as a \emph{full} subcategory. Besides from smooth manifolds, we have three further examples of sets with a canonical diffeology:
\begin{enumerate}
\item 
If $X$ and $Y$ are diffeological spaces, the set $\smsp(X,Y)$ of diffeological maps between $X$ and $Y$ is a diffeological space in the following way: a map
\begin{equation*}
c: U \to \smsp(X,Y)
\end{equation*}
is a plot if and only if for any plot $c':V \to X$ of $X$ the composite 
\begin{equation*}
\alxydim{@C=1.5cm}{U \times V \ar[r]^-{c \times c'} & \smsp(X,Y) \times X \ar[r]^-{\mathrm{ev}}
& Y}
\end{equation*}
is a plot of $Y$. Here, $\mathrm{ev}$ denotes the evaluation map $\mathrm{ev}(f,x):=f(x)$.

\item
Subsets $Y$ of a diffeological space $X$ are diffeological: its plots are those plots of $X$ whose image is contained in $Y$.

\item
If $X$ is a diffeological space, $Y$ is a set and $p:X \to Y$ is a map, $Y$ becomes  a diffeological space whose plots are those maps $c:U \to Y$ for which there exists a cover of $U$ by open sets $U_{\alpha}$ and plots $c_{\alpha}:U_{\alpha} \to X$ of $X$ such that $c|_{U_{\alpha}}=p\circ c_{\alpha}$.

\end{enumerate}

Equipped with these examples,  the sets $P^iX$ we have defined in Section \ref{sec1_1} become diffeological spaces in the following way. The set $PX$ is a subset of the diffeological space $D^{\infty}([0,1],X)$, and hence a diffeological space. Then we consider one of the projections $\mathrm{pr}^i:PX \to P^iX$ from (\ref{34}) to either reparameterization classes, thin homotopy classes or homotopy classes. According to  the third example above, all the sets $P^iX$ become diffeological spaces. We also have examples of diffeological maps: if $f:X \to Y$ is a smooth map between smooth manifolds, the induced maps $f_{*}:P^iX \to P^iY$ are all diffeological.

The most important question for us will be, whether a map $P^iX \to M$ from one of these diffeological spaces to a smooth manifold $M$ -- regarded as a diffeological space -- is diffeological. From the definitions above one can deduce the following result, and the reader is free to take it either as a result from the background of diffeological spaces, or  as a definition. 

\begin{lemma}[\cite{schreiber3}, Proposition A.7 \textit{i)}]
\label{lem2}
A map $f:P^iX \to M$ is diffeological, if and only if  for every  $k \in \N_0$, every open subset $U \subset \R^k$  and every map $c: U \to PX$ such that the composite
\begin{equation*}
\alxydim{@C=1.6cm}{U \times [0,1] \ar[r]^-{c \times \id} & PX \times [0,1] \ar[r]^-{\mathrm{ev}} & X}
\end{equation*}
is smooth, also the map
\begin{equation*}
\alxydim{}{U \ar[r]^-{c} & PX \ar[r]^-{\mathrm{pr}^i} & P^i X \ar[r]^-{f} & M}
\end{equation*} 
is smooth.
\end{lemma}

Now we can study \emph{smooth} functors $F:\mathcal{P}_1(X) \to S$ to a Lie category $S$: on objects $F: X \to S_0$ is a smooth map and on morphisms $F: P^1X \to S_1$ is a diffeological map. Similarly, if $\eta:F \to F'$ is a natural transformation between two smooth functors, 
 it is called   \emph{smooth} natural transformation, if its components $\eta(x) \in S_1$ furnish a smooth map $X \to S_1$. Smooth functors $F$ and smooth natural transformations $\eta$ form a category $\mathrm{Funct}^{\infty}(\mathcal{P}_1(X),S)$.

\begin{remark}
\normalfont
Concerning Definition \ref{def2} of a diffeological space, several different conventions for  plots are common. For example, in order to deal properly with manifolds with boundary or corners, it is more convenient to consider plots being defined on convex subsets $U \subset \R^k$ rather than open ones \cite{baez3}. Such questions do not affect the results of this article, since we consider either maps defined on manifolds without boundary or maps which are constant near the boundary of a manifold, for example paths with sitting instants. 
\end{remark}

\medskip

\subsection{Equivalence between Functors and Forms}

\label{sec1_3}

In \cite{schreiber3}  we have established an isomorphism between two categories,
\begin{equation}
\label{17}
\mathrm{Funct}^{\infty}(\mathcal{P}_1(X),\mathcal{B} G) \cong \diffco{G}{1}{X}\text{.}
\end{equation}
Both categories depend on a smooth manifold $X$ and a Lie group $G$.  On the left hand side we have the category of smooth functors from the path groupoid $\mathcal{P}_1(X)$ of $X$ to the Lie groupoid $\mathcal{B} G$.\label{not:smoothfun}\label{not:oneobject} We recall from the introduction that the Lie groupoid $\mathcal{B}G$ has one object, and its set of morphisms is the Lie group $G$. The composition is defined by $g_2 \circ g_1 := g_2g_1$.  On the right hand side we have a category $\diffco{G}{1}{X}$\label{not:gconn}  defined as follows:
its objects are 1-forms $A\in\Omega^1(X,\mathfrak{g})$ with values in the Lie algebra $\mathfrak{g}$ of $G$, and a morphism $g:A \to A'$ is a smooth function $g:X \to G$ such that 
\begin{equation}
\label{18}
A' = \mathrm{Ad}_{g}(A) - g^{*}\bar\theta\text{,}
\end{equation}
where $\bar\theta$ is the right invariant Maurer-Cartan form on $G$. 
The identity morphism is the constant function $g=1$ and the composition is the multiplication of functions, $g_2 \circ g_1 := g_2g_1$.
The equivalence (\ref{17}) can be given explicitly in both directions: there are two functors \label{not:dp}
\begin{equation*}
\alxydim{@C=2cm}{\mathrm{Funct}^{\infty}(\mathcal{P}_1(X),\mathcal{B} G) \ar@/^1.5pc/[r]^{\fo} & \diffco{G}{1}{X} \ar@/^1.5pc/[l]^{\fu}}
\end{equation*}
whose definitions  we  shall review in the following. Several details and proofs will be skipped and can be found in \cite{schreiber3}. 

\medskip

Given a smooth functor $F:\mathcal{P}_1(X) \to \mathcal{B} G$,  a 1-form $A\in\Omega^1(X,\mathfrak{g})$ is defined in the following three steps:
\begin{enumerate}
\item 
For a point $x\in X$ and a tangent vector $v\in T_xX$, we choose a smooth curve $\Gamma:\R \to X$ with $\Gamma(0)=x$ and $\dot\Gamma(0)=v$. Let $\gamma_{\R}(t_0,t)$ be the (up to thin homotopy unique) path in $\R$ that goes from $t_0$ to $t$, regarded as a map $\gamma_{\R}: \R^2 \to P^1\R$. 
\item
One can show (Appendix B.4 in \cite{schreiber3}) that the composite
\begin{equation}
\label{19}
F_{\Gamma} := F \circ \Gamma_{*} \circ \gamma_{\R}: \R^{2} \to G
\end{equation}
is a smooth map with $F_{\Gamma}(t_0,t_0)=1$ for all $t_0\in \R$. We define
\begin{equation}
\label{22}
A_x(v) :=- \left. \frac{\mathrm{d}}{\mathrm{d} t} \right|_{0} F_{\Gamma}(0,t) \in \mathfrak{g}\text{.}
\end{equation}

\item 
One can then verify that the value $A_x(v)$ is independent of the choice of $\Gamma$ (Lemma B.2 in \cite{schreiber3}), and that the assignment $A: TX \to \mathfrak{g}$ is smooth and linear (Lemma B.3 in \cite{schreiber3}). This defines the 1-form $\fo(F):=A$.
\end{enumerate}
The components of a smooth natural transformation $\rho: F \to F'$ form by definition a  smooth map $\fo(\rho) := g:X \to G$. Let again $\Gamma: \R \to X$ be a smooth curve and $F_{\Gamma}$ and $F'_{\Gamma}$ the functions (\ref{19}) associated to the functors $F$ and $F'$, and let $g_{\Gamma} := g \circ \Gamma$.  The naturality of $\rho$ implies the equation
\begin{equation*}
g_{\Gamma}(t) \cdot F_{\Gamma}(0,t) = F'_{\Gamma}(0,t) \cdot g_{\Gamma}(0)\text{,}
\end{equation*}
whose derivative evaluated at $t=0$ shows (\ref{18}) for $A$ and $A'$ the 1-forms defined by $F_{\Gamma}$ and $F_{\Gamma}'$, respectively. Hence, $\fo(\rho)$ is a morphism in $\diffco{G}{1}{X}$; this defines the functor $\fo$.

\medskip

Conversely, consider a 1-form $A\in\Omega^1(X,\mathfrak{g})$. Then, a smooth functor $F:\mathcal{P}_1(X) \to \mathcal{B} G$ is defined in the following way:
\begin{enumerate}
\item 
Let $\gamma$ be a path in $X$, which we extend trivially to $\R$ by $\gamma(t):=\gamma(0)$ for $t<0$ and $\gamma(t):=\gamma(1)$ for $t>1$. We pose the initial value problem
\begin{equation}
\label{36}
\frac{\partial }{\partial t}u_{\gamma}(t) = - \mathrm{d}r_{u(t)}|_1 \left ( A_{\gamma(t)} \left ( \frac{\mathrm{d} \gamma}{\mathrm{d} t}\right) \right )
\quad\text{ and }\quad
u(t_{0})=1
\end{equation} 
for a smooth function $u:\R \to G$ and  fixed $t_0\in \R$. Here $\mathrm{d}r_{u(t)}$ is the differential of the multiplication with $u(t)$ from the right.  

\item
The initial value problem (\ref{36}) has a unique  solution $f_{A,\gamma}(t_0,t)$, from which we define a map
\begin{equation}
\label{23}
F: PX \to G: \gamma \mapsto f_{A,\gamma}(0,1)\text{.}
\end{equation}
We remark that this map is sometimes referred to as the ,,path-ordered exponential`` 
\begin{equation}
\label{74}
F(\gamma) = \mathcal{P}\mathrm{exp}\left( \int_{\gamma}A \right )\text{.}
\end{equation}

\item
One can show that the map $F$ is independent of the thin homotopy class of $\gamma$ (Proposition 4.3 in \cite{schreiber3}), and that it factors through a smooth map $F: P^1X \to G$ (Lemma 4.5 in \cite{schreiber3}). It respects the composition of paths so that we have defined a smooth functor $\fu(A):=F$.
\end{enumerate}

For a smooth function $g:X \to G$ considered as a morphism $g:A \to A'$ between two 1-forms $A,A'\in\Omega^1(X,\mathfrak{g})$ we need to define an associated smooth natural transformation $\rho=\fo(g):F \to F'$ between the associated functors. The component of $\rho$ at $x\in X$ is defines as $g(x)$. One can then show that $g(y) \cdot f_{A,\gamma}(t_{0},t) \cdot g(x)^{-1}$ solves the initial value problem (\ref{36}) for $A'$ and $\gamma$ (Lemma 4.2 in \cite{schreiber3}), which implies the naturality of the natural transformation $\rho$. This defines the functor $\fu$.

\begin{theorem}[Proposition 4.7 in \cite{schreiber3}]
\label{th3}
Let $X$ be a smooth manifold and $G$ be a Lie group.
The two functors $\fo$ and $\fu$ satisfy
\begin{equation*}
\fo \circ \fu = \id_{\diffco{G}{1}{X}}
\quad\text{ and }\quad
\fu \circ \fo = \id_{\mathrm{Funct}^{\infty}(\mathcal{P}_1(X),\mathcal{B}G)}\text{,}
\end{equation*}
in particular, they form an isomorphism of categories.
\end{theorem}

We give a short sketch of the proof. If we start with a 1-form $A\in\Omega^1(X,\mathfrak{g})$, we shall test the 1-form $\fo(\fu(A))$ at a point $x\in X$ and a tangent vector $v\in T_xX$. Let $\Gamma:\R \to X$ be a curve in $X$ with $x=\Gamma(0)$ and $v=\dot \Gamma(0)$. If we further denote $\gamma_{\tau} := \Gamma_{*}(\gamma_{\R}(0,\tau)) \in PX$ we have
\begin{multline*}
-\fo(\fu(A))|_x(v) \stackrel{(\ref{22})}{=} \left . \frac{\partial}{\partial \tau} \right|_{0} \fu(A)_\Gamma(0,\tau)\\\stackrel{(\ref{19})}{=} \left . \frac{\partial}{\partial \tau} \right |_{0}\fu(A)(\gamma_{\tau})\stackrel{(\ref{23})}{=} \left . \frac{\partial}{\partial \tau} \right |_{0}f_{A,\gamma_{\tau}}(0,1)
\end{multline*}
Here, $f_{A,\gamma_{\tau}}$ denotes the unique solution of the initial value problem (\ref{36}) for $\gamma_{\tau}$. A uniqueness argument shows $f_{A,\gamma_\tau}( t_0,t) = f_{A,\gamma_1}(\tau t_0,\tau t)$, so that
\begin{equation*}
\left . \frac{\partial }{\partial \tau}f_{A,\gamma_\tau}(0,t) \right |_{\tau=0,t=1} =\left . \frac{\partial }{\partial t} f_{A,\gamma_1}(0,t)  \right |_{t=0} =-A_{p}(v)\text{,}
\end{equation*}
this yields $\fo(\fu(A))=A$. 

On the other hand, if $F: \mathcal{P}_1(X) \to \mathcal{B} G$ is a smooth functor, we test the functor $\fu(\fo(F))$ on a path $\gamma$ in $X$.  By (\ref{23}), 
\begin{equation*}
\fu(\fo(F))(\gamma) = f_{\fo(F),\gamma}(0,1)
\end{equation*}
 where $f_{\fo(F),\gamma}$ is the solution of the initial value problem (\ref{36}) for the 1-form $\fo(F)$ and the path $\gamma$. Due to the definition (\ref{22}) of $\fo(F)$ by the function $F_\gamma:\R^2 \to G$ we have
\begin{equation*}
(\gamma^{*}\fo(F))_{t} \left (\frac{\partial}{\partial t} \right ) =- \left . \frac{\partial}{\partial \tau} \right |_{\tau=0}F_{\gamma}(t,t+\tau)\text{.}
\end{equation*}
Since $F$ is a functor, $F_{\gamma}(x,z)=F_{\gamma}(y,z)F_{\gamma}(x,y)$. Both together show that $F_{\gamma}$ also solves the initial value problem, so that, by uniqueness, 
\begin{equation*}
f_{\fo(F),\gamma}(0,1)= F_{\gamma}(0,1) = F(\gamma)\text{.}
\end{equation*}
This shows $\fu(\fo(F)) = F$.

\medskip

Remarkably, there is not much structure that is preserved by the functors $\fu$ and $\fo$ (unless the Lie group $G$ is abelian). For example, sums and negatives of differential forms, or products and inverses of smooth functors are  all \emph{not} preserved. We only know the following fact:

\begin{proposition}
\label{prop6}
The functors $\fu$ and $\fo$ are compatible with pullbacks along a smooth map $f:X \to Y$ between smooth manifolds $X$ and $Y$, i.e.
\begin{equation*}
\fu(f^{*}A) = f^{*}\fu(A)
\quad\text{ and }\quad
\fo(f^{*}F) = f^{*}\fo(F)
\end{equation*}
for a 1-form $A\in\Omega^1(Y,\mathfrak{g})$ and a smooth functor $F:\mathcal{P}_1(Y) \to \mathcal{B}G$, and similarly for morphisms. 
\end{proposition}

Here we have used the notation $f^{*}F$ for the functor $F \circ f_{*}$, where $f_{*}$ is the induced map on path groupoids.
Proposition \ref{prop6} follows in a straightforward way from the naturality of the definitions of the functors $\mathcal{D}$ and $\mathcal{P}$.
\section{Smooth 2-Functors and Differential Forms}

\label{sec3}

In this section we generalize Theorem \ref{th3} -- the equivalence between 1-forms and smooth functors -- to 2-functors. The basic 2-categorical notions such as 2-categories, 2-functors, pseudonatural transformations and modifications are summarized in Appendix \ref{app1}; for the reader familiar with these notions it is important to notice that all 2-categories and 2-functors in this article are assumed to be strict without further notice. 

The first step in  concerns the path groupoid $\mathcal{P}_1(X)$ that was present in Theorem \ref{th3}: in Section \ref{sec2_1} we define the path 2-groupoid $\mathcal{P}_2(X)$ associated to a smooth manifold $X$. Instead of the Lie group $G$ that was present in Theorem \ref{th3} we  use a (strict) Lie 2-group $\mathfrak{G}$. In the same way that a category $\mathcal{B}G$ is associated to any Lie group $G$, a 2-category $\mathcal{B}\mathfrak{G}$ is associated to any Lie 2-group $\mathfrak{G}$, and the 2-functors we consider are of the form
\begin{equation*}
\label{57}
F: \mathcal{P}_2(X) \to \mathcal{B}\mathfrak{G}\text{.}
\end{equation*}
A convenient and  concrete way to deal with Lie 2-groups is provided by crossed modules \cite{whitehead1, brown1}. Their definition, their relation to Lie 2-groups, and the associated 2-categories $\mathcal{B}\mathfrak{G}$ are  described in Appendix \ref{app2}.

The announced generalization of Theorem \ref{th3} is worked out in three steps: in Section \ref{sec3_1} we extract differential forms from 2-functors, pseudonatural transformations and modifications. We derive conditions on the extracted differential forms that lead us straightforwardly to an appropriate generalization $\diffco{\mathfrak{G}}{2}{X}$ of the category $\diffco{G}{1}{X}$ that was present in Theorem \ref{th3}. The goal of Section \ref{sec3_1} is that extracting differential forms furnishes 2-functor
\begin{equation*}
\fo: \mathrm{Funct}^{\infty}(\mathcal{P}_2(X),\mathcal{B}\mathfrak{G}) \to \diffco{\mathfrak{G}}{2}{X}\text{.}
\end{equation*}
In Section \ref{sec3_2} we introduce a 2-functor 
\begin{equation*}
\fu: \diffco{\mathfrak{G}}{2}{X} \to \mathrm{Funct}^{\infty}(\mathcal{P}_2(X),\mathcal{B}\mathfrak{G})
\end{equation*}
in the opposite direction, that reconstructs 2-functors, pseudonatural transformations and modifications from given differential forms. Finally, we prove in Section \ref{sec3_3} the main result of this article, namely that the 2-functors
$\fo$ and $\fu$ establish an isomorphism of 2-categories.

\subsection{The Path 2-Groupoid of a Smooth Manifold}

\label{sec2_1}

As mentioned in the introduction, the path 2-groupoid is obtained by adding 2-morphisms to the path groupoid $\mathcal{P}_1(X)$. These 2-morphisms are smooth homotopies in the sense of Definition \ref{def3} without the restriction (2) on their rank, explicitly:

\begin{definition}
Let $\gamma_0,\gamma_1:x \to y$ be paths in $X$.
A \emph{bigon} $\Sigma:\gamma_0 \Rightarrow \gamma_1$ is a smooth map $\Sigma:
[0,1]^2 \to X$ such that there exists a number $0<\varepsilon <\frac{1}{2}$
with
\begin{enumerate}
\item[a)]
$\Sigma(s,t) = x$ for $0 \leq t < \varepsilon$ and $\Sigma(s,t)=y$ for $1-\varepsilon
< t \leq 1$.
\item[b)] 
$\Sigma(s,t)=\gamma_0(t)$ for $0 \leq s <\varepsilon$ and $\Sigma(s,t)=\gamma_1(t)$
for $1-\varepsilon < s \leq 1$.
\end{enumerate}
\end{definition}

We denote the set
of  bigons in $X$ by $BX$. \label{not:bx}Bigons can be composed in two  ways: If $\Sigma:\gamma_1
\Rightarrow \gamma_2$ and $\Sigma':\gamma_2 \Rightarrow \gamma_3$ are bigons,
we have a new bigon $\Sigma'\bullet\Sigma: \gamma_1 \Rightarrow \gamma_3$
defined by
\begin{equation}
\label{sec2_4}
(\Sigma' \bullet \Sigma)(s,t)  = \begin{cases}\Sigma(2s,t) & \text{ for
}0 \leq s < \frac{1}{2} \\
\Sigma'(2s-1,t) & \text{ for }\frac{1}{2}\leq s \leq 1\text{;} \\
\end{cases}
\end{equation}
and for two bigons $\Sigma_1: \gamma_1 \Rightarrow \gamma_1'$ and $\Sigma_2:\gamma_2 \Rightarrow \gamma_2'$ such that $\gamma_1(1)=\gamma_2(0)$, we have another new bigon $\Sigma_2 \circ \Sigma_1: \gamma_2 \circ \gamma_1 \Rightarrow \gamma_2' \circ \gamma_1'$
 defined by
\begin{equation}
\label{sec2_5}
(\Sigma_2 \circ \Sigma_1)(s,t)  := \begin{cases}\Sigma_1(s,2t) & \text{ for
}0 \leq t < \frac{1}{2} \\
\Sigma_2(s,2t-1) & \text{ for }\frac{1}{2}\leq t \leq 1\text{.} \\
\end{cases}
\end{equation}
Due to the sitting instants, the new maps (\ref{sec2_4}) and (\ref{sec2_5}) are again smooth and have sitting instants. 

Like in the case of paths, there are several equivalence relations on the set $BX$ of bigons in $X$, starting with reparameterization classes, and continued by  a ladder of types of homotopy classes, graded by an upper bound for the rank of the homotopies. The corresponding quotient spaces are denoted by
\begin{equation*}
\alxydim{}{BX \ar[r] & B^0X \ar[r] & B^1 X  \ar[r] & B^2X \ar[r] & B^3X\text{.}}
\end{equation*}
In this article we are only interested in $B^2X=BX/\!\sim_2$. 
\begin{definition}
\label{def1}
Two bigons $\Sigma:\gamma_0 \Rightarrow \gamma_1$ and $\Sigma':\gamma'_0 \Rightarrow \gamma'_1$
are called \emph{thin ho\-mo\-to\-py equivalent}, denoted $\Sigma \sim_2 \Sigma'$, if there exists a smooth
map $h: [0,1]^3 \to X$ such that
\begin{enumerate}
\item[(1)] 
it has sitting instants: there exists a number $0 < \varepsilon < \frac{1}{2}$ with
\begin{enumerate}
\item[a)]
$h(r,s,t)=x$  for $0 \leq t < \varepsilon$ and $h(r,s,t)=y$ for $1-\varepsilon
< t \leq 1$.
\item[b)]
$h(r,s,t)=h(r,0,t)$  for $0 \leq s < \varepsilon$ and $h(r,s,t)=h(r,1,t)$ for $1-\varepsilon
< s \leq 1$.
\item[c)]
$h(r,s,t) = \Sigma(s,t)$  for all  $0\leq r < \varepsilon$ and $h(r,s,t)=\Sigma'(s,t)$
for all $1-\varepsilon < r \leq 1$. 
\end{enumerate}

\item[(2)]
the differential of $h$ satisfies
\begin{enumerate}
\item[a)]
$\mathrm{rank}(\mathrm{d}h|_{(r,s,t)}) \leq 2$ for all $r,s,t\in [0,1]$, and
\item[b)]
$\mathrm{rank}(\mathrm{d}h|_{(r,i,t)}) \leq 1$ for $i=0,1$ fixed.
\end{enumerate}
\end{enumerate}
\end{definition}

Condition (1) assures that thin homotopy is an equivalence
relation on $BX$.  Condition (2b) asserts that two thin homotopy equivalent bigons $\Sigma:\gamma_0 \Rightarrow \gamma_1$ and $\Sigma':\gamma_0' \Rightarrow \gamma_1'$ start and end on thin homotopy equivalent paths
$\gamma_0 \sim_{1} \gamma_0'$ and $\gamma_1 \sim_{1} \gamma_1'$.
The composition $\circ$ of two bigons defined above clearly induces a well-defined composition on $B^2X$. For the composition $\bullet$ this is more involved: let $\Sigma:\gamma_1 \Rightarrow \gamma_2$ and $\Sigma': \gamma'_2 \Rightarrow \gamma_3$ be two bigons  such that $\gamma_2 \sim_1 \gamma_2'$. Let $h:[0,1] \to X$ be any thin homotopy between $\gamma_2$ and $\gamma_2'$; this is a particular bigon $h:\gamma_2 \Rightarrow \gamma_2'$. Now we define the composition of the corresponding classes in $B^2X$ by
\begin{equation*}
[\Sigma']_{\sim_2} \bullet [\Sigma]_{\sim_2} := [\Sigma' \bullet h \bullet \Sigma]_{\sim_2}\text{.}
\end{equation*}
The proof that this is independent of the choice of $h$ requires a technical computation carried out in \cite{martins1}.  Another important fact is that the two compositions $\circ$ and $\bullet$ are compatible with each other in the sense that
\begin{equation}
\label{31}
(\Sigma'_1 \bullet \Sigma'_2) \circ (\Sigma_1 \bullet \Sigma_2) \sim_2 (\Sigma_1'\circ\Sigma_1)\bullet(\Sigma_2'\circ\Sigma_2)
\end{equation}
whenever all these compositions are well-defined. 

\begin{definition}
\label{not:p2g}
The \emph{path 2-groupoid} $\mathcal{P}_2(X)$ of a smooth manifold $X$ is the
2-category whose  objects are the points of $X$, whose set of 1-morphisms is $P^{1}X$, and whose set of 2-morphisms is $B^2X$. Horizontal and vertical composition are given by $\circ$ and $\bullet$, and the identities are  the identity path $\id_x:x \to x$ and the identity bigon $\id_{\gamma}:\gamma\Rightarrow \gamma$ defined by $\id_{\gamma}(s,t) := \gamma(t)$. 
\end{definition}

The axioms of a 2-category (see Definition \ref{def6}) are satisfied: Axiom (C1) follows from the second equation in (\ref{32}), axiom (C2) follows from the first equation in (\ref{32}) and from an elementary construction of  homotopies $\Sigma \bullet \id_{\gamma_1} \sim_2 \Sigma \sim_2 \id_{\gamma_2} \bullet \Sigma$. Axiom (C3) is (\ref{31}).  It is also clear the the category $\mathcal{P}_2(X)$ is indeed a groupoid. 

For parallel transport along surfaces, the path 2-groupoid plays the same role  the path groupoid $\mathcal{P}_1(X)$  plays for parallel transport along curves (see Section \ref{sec1_1}): the corresponding geometric objects are  (weak) 2-functors
\begin{equation*}
\mathrm{tra}:\mathcal{P}_2(X) \to T
\end{equation*}
into some 2-category $T$, as outlined in Section 6 of \cite{schreiber3}. A detailed discussion of these 2-functors will follow in  \cite{schreiber2}. 
\medskip

In exactly the same way as we have diffeological maps $P^iX \to M$ we have diffeological maps from all the equivalence classes $B^{i}X$ of bigons in $X$ to smooth manifolds $M$. Analogous to Lemma \ref{lem2}, we have

\begin{lemma}
\label{lem13}
A map $f: B^iX \to M$ is diffeological if and only if for every $k\in \N_0$, every open subset $U \subset \R^k$ and every map $c: U \to BX$ such that the composite
\begin{equation*}
\alxydim{@C=1.6cm}{U \times [0,1]^2 \ar[r]^-{c \times \id} & BX \times [0,1]^2 \ar[r]^-{\mathrm{ev}} & X}
\end{equation*}
is smooth, also the map 
\begin{equation*}
\alxydim{}{U \ar[r]^{c} & BX \ar[r]^-{\mathrm{pr}^i} & B^iX \ar[r]^-{f} & M}
\end{equation*}
is smooth.
\end{lemma}

This admits to define smooth 2-functors defined on the path 2-groupoid of $X$ with values in \emph{smooth 2-categories} $S$: 2-categories for which objects $S_0$, 1-morphisms $S_1$ and 2-morphisms $S_2$ are smooth manifolds and all structure maps are smooth. 

\begin{definition}
\label{def7}
A 2-functor $F: \mathcal{P}_2(X) \to S$ from the path 2-groupoid of a smooth manifold $X$ to a smooth 2-category $S$ is called \emph{smooth}, if 
\begin{enumerate}
\item 
on objects, $F:X \to S_{0}$ is smooth.
\item
on 1-morphisms, $F: P^1X \to S_1$ is diffeological (see Lemma \ref{lem2}).
\item
on 2-morphisms, $F: B^2X \to S_2$ is diffeological (see Lemma \ref{lem13}).
\end{enumerate}
\end{definition}

For the definitions of morphisms between 2-functors, the pseudonatural transformations, and morphisms between those, the modifications, we refer the reader again to Appendix \ref{app1}. A pseudonatural transformation $\rho: F \to F'$  is called \emph{smooth}, if its components $\rho(x) \in S_1$ at objects $x\in X$ furnish a smooth map $X \to S_1$, and its components $\rho(\gamma)\in S_2$ at 1-morphisms $\gamma\in P^1X$ furnish a diffeological map $P^1 X \to S_2$. Similarly, a modification $\mathcal{A}: \rho \Rightarrow \rho'$ is called smooth, if its components $\mathcal{A}(x)\in S_2$ from a smooth map $X \to S_2$.
Summarizing, these structures form  a 2-category $\mathrm{Funct}^{\infty}(\mathcal{P}_2(X),S)$. 

\subsection{From Functors to Forms}

\label{sec3_1}

As we explain in Appendix \ref{app2} that the 2-category $\mathcal{B}\mathfrak{G}$ associated to a Lie 2-group $\mathfrak{G}$ which is represented by a smooth crossed module $(G,H,t,\alpha)$ has one object, the set of morphisms is  $G$ and the set of 2-morphisms is the semi-direct product $G \ltimes H$, where $G$ acts on  $H$ via a smooth map   $\alpha:G \times H \to H$. The guideline how to extract differential forms from smooth 2-functors is the same as reviewed in Section \ref{sec1_3}: we evaluate the Lie group-valued functors on certain paths, obtain Lie group-valued maps, and take their derivative.

\subsubsection{Extracting Forms I: 2-Functors}

\label{sec_ex1}

Here we start with a given smooth 2-functor
\begin{equation*}
F:\mathcal{P}_2(X) \to \mathcal{B}\mathfrak{G}\text{.}
\end{equation*}
Clearly, $F$ restricted to objects and 1-morphisms is just a smooth 1-functor 
$F_{0,1}:\mathcal{P}_1(X) \to \mathcal{B} G$.
By Theorem \ref{th3} it corresponds to a $\mathfrak{g}$-valued 1-form $A$ on $X$. From the remaining map $F_2: B^2X \to G \ltimes H$ we now define an $\mathfrak{h}$-valued 2-form $B$ on $X$. Its definition is pointwise: let $x\in X$ be a point and $v_1,v_2\in T_xX$  be tangent vectors. We choose a smooth map $\Gamma:\R^2 \to X$ with $x= \Gamma(0)$ and
\begin{equation}
\label{27}
v_1 = \left . \frac{\mathrm{d}}{\mathrm{d}s} \right |_{s=0} \Gamma(s,0)
\quad\text{ and }\quad
v_2 = \left . \frac{\mathrm{d}}{\mathrm{d}t} \right |_{t=0} \Gamma(0,t)\text{.}
\end{equation}
Note that in $\R^2$ there is only one thin homotopy class of bigons between each two fixed paths. In particular, we have a canonical family  $\Sigma_{\R}: \R^2 \to B^2\R^2$, where
\begin{equation}
\label{25}
\Sigma_{\R}(s,t) := \alxydim{}{(0,0) \ar[r] \ar[d] & (0,t) \ar@{=>}[dl] \ar[d] \\ (s,0) \ar[r] & (s,t)}\text{.}
\end{equation}
We use this canonical family of bigons to produce a map
\begin{equation}
\label{9}
F_\Gamma := p_H \circ F_{2} \circ \Gamma_{*} \circ \Sigma_{\R}: \R^2 \to  H
\end{equation}
where $p_H: G \ltimes H \to H$ is the projection to the second factor.
\begin{lemma}
\label{lem3}
The map $F_\Gamma: \R^2 \to H$ is smooth.  Furthermore, its second mixed derivative evaluated at $0\in \R^2$ is a well-defined element in the Lie algebra $\mathfrak{h}$ of $H$, and is independent of the choice of $\Gamma$, i.e. if $\Gamma_0,\Gamma_1:\R^2 \to X$ are smooth maps with $\Gamma_0(0)=\Gamma_1(0)=x$ and both satisfying (\ref{27}), then
\begin{equation}
\label{26}
 \left . \frac{\partial^2F_{\Gamma_0}}{\partial s \partial t}\right |_{(0,0)} = \left . \frac{\partial^2F_{\Gamma_1}}{\partial s \partial t}\right |_{(0,0)}\text{.}
\end{equation}
\end{lemma}

\proof
The smoothness of $F_\Gamma$ follows from the smoothness of the 2-functor $F$ as explained in Section \ref{sec1_2}:  the relevant evaluation map
\begin{equation*}
\alxydim{@C=1.6cm}{\R^2 \times [0,1]^2 \ar[r]^-{\Gamma_{*} \circ \Sigma_{\R} \times \id} & BX \times [0,1]^2 \ar[r]^-{\mathrm{ev}}& X}
\end{equation*} 
is smooth. Since $F$ is smooth on 2-morphisms, $F \circ \Gamma_{*} \circ \Sigma_{\R}$ is smooth.

Next we  notice that   $F_{\Gamma}(0,t)=F_{\Gamma}(s,0)=1$ for all $s,t\in \R$, so that the second mixed derivative naturally takes values  in $\mathfrak{h}$. Now we consider the two  families 
\begin{equation*}
\Sigma_k := (\Gamma_k)_{*} \circ \Sigma_{\R}: \R^2 \to B^2X
\end{equation*}
of bigons in $X$. We work in a suitable open neighborhood $V \subset \R^2$ of the origin, and are going to construct a homotopy
$H: V \times [0,1] \to BX$ with $H(x,y,k)=\Sigma_k(x,y)$ for $k=0,1$. Then we will show that $H$ factors through the smooth map
\begin{equation*}
f: V \times [0,1] \to Z: (x,y,\alpha) \mapsto (x,y,(x^2+y^2)\alpha)\text{,}
\end{equation*}
where $Z:=\left \lbrace (x,y,z) \in V \times [0,1] \;|\; 0 \leq z \leq x^2+y^2 \right \rbrace$, and another smooth map $B: Z \to BX$. Applying the chain rule to the decomposition  $H=B \circ f$ gives
\begin{eqnarray*}
 \left . \frac{\partial^2F_{\Gamma_k}}{\partial s \partial t}\right |_{(0,0)} &=& \left . \frac{\partial^2}{\partial s\partial t} \right |_{(0,0)} (F \circ \Sigma_k)
\;\;=\;\; \left . \frac{\partial^2}{\partial s\partial t} \right |_{(0,0,k)} (F \circ B \circ f)
\\&=& \mathrm{H}(F \circ B)|_{f(0,0,k)} \left ( \left.\frac{\partial f}{\partial x} \right|_{(0,0,k)} , \left. \frac{\partial f}{\partial y} \right|_{(0,0,k)} \right ) \\ && \hspace{2cm}+\; \mathrm{D}(F \circ B)|_{f(0,0,k)} \left ( \left . \frac{\partial^2 f}{ \partial x \partial y} \right|_{(0,0,k)} \right )\text{,}
\end{eqnarray*}
where $\mathrm{H}(F \circ B)$ denotes the Hesse matrix of $F \circ B$, considered as a symmetric, fibre-wise bilinear form $TZ \times_Z TZ \to \mathfrak{h}$.
By construction of the map $f$,  the latter expression is independent of $k$.

In order to construct $H$ and $B$, we
work in a chart that identifies $V$ with an open neighborhood of $x$, and form the ,,linear interpolation`` 
\begin{equation*}
h: V \times [0,1] \times [0,1]^2 \to X: (x,y,\alpha,s,t) \mapsto \Sigma_0(x,y)(s,t) + \alpha \cdot d(x,y,s,t)\text{,}
\end{equation*}
where the ,,difference`` $d$ is given by
\begin{equation*}
d: V \times [0,1]^2 \to X: (x,y,s,t) \mapsto \Sigma_1(x,y)(s,t) - \Sigma_0(x,y)(s,t)\text{.}
\end{equation*}
Then, we set $H(x,y,\alpha)(s,t) := h(x,y,s,t,\alpha)$. Next we construct the map $B$. The coincidence of the values and the first derivatives of the maps $\Gamma_0$, $\Gamma_1$ at $(0,0)$ imply via Hadamard's lemma that there exist smooth maps $a,b,c: V \times [0,1]^2 \to X$ such that
\begin{equation*}
d(x,y,s,t) = x^2 \cdot a(x,y,s,t) + y^2 \cdot b(x,y,s,t) + 2xy \cdot c(x,y,s,t)\text{.}
\end{equation*}

Now we change to polar coordinates. We denote by $U \subset \R_{\geq 0} \times [0,2\pi)$ a suitable  open neighborhood of $(0,0)$ so that the coordinate transformation $\tau: (r,\phi) \mapsto (r \cdot \cos\phi,r\cdot \sin\phi)$ is a map $\tau:U \to V$. Then, we get
\begin{equation*}
d_{p}(r,\phi,s,t) := d(\tau(r,\phi),s,t) = r^2\cdot \tilde d_p(r,\phi,s,t)
\end{equation*}
with
\begin{multline*}
\tilde d_p(r,\phi,s,t) := \cos^2\phi \cdot a(\tau(r,\phi),s,t) + \sin\phi\cos\phi\cdot b(\tau(r,\phi),s,t) \\ + \sin^2\phi\cdot c(\tau(r,\phi),s,t)
\end{multline*}
defining a smooth map $\tilde d_p: U \times [0,1]^2 \to X$. We look at the smooth map
\begin{equation*}
b_p: U \times [0,1] \times [0,1]^2 \to X: (r,\phi,\alpha,s,t) \mapsto \Sigma_0(\tau(r,\phi),s,t) + \alpha \cdot \tilde d_p(r,\phi,s,t)\text{.}
\end{equation*}
We consider $Z_p := \left \lbrace (r,\phi,z) \in U \times  [0,1] \;\;\right | \;\left . 0\leq z \leq r^{2} \; \right \rbrace$, and claim that the restriction of $b_p$ to
$Z_p \times [0,1]^2$ 
transforms back into a smooth map in \emph{cartesian} coordinates. To see this, it suffices to notice that the term $\alpha \cdot \tilde d_p(r,\phi,s,t)$ and all its $r$-derivatives vanish at $r=0$, since then also $\alpha=0$. This way we obtain a smooth map $b: Z \times [0,1]^2 \to X$ satisfying 
\begin{equation*}
b(\tau(r,\phi),\alpha,s,t) = b_p(r,\phi,\alpha,s,t)
\end{equation*}
for all $(r,\phi) \in U$.

Finally, we set $B(x,y,\alpha)(s,t) := b(x,y,s,t,\alpha)$. It is straightforward to see that $B(x,y,\alpha)$ is a bigon, and the smoothness of $b$ implies the smoothness of the map  $B: Z \to BX$. We calculate for $(r,\phi) \in U$:
\begin{eqnarray*}
(B \circ f)(\tau(r,\phi),\alpha)(s,t) 
&=& b_p(r,\phi,r^{2}\alpha,s,t) 
\\&=& \Sigma_0(\tau(r,\phi),s,t) + \alpha \cdot r^2 \cdot \tilde d_p(r,\phi,s,t)
\\&=& \Sigma_0(\tau(r,\phi),s,t) + \alpha\cdot d(\tau(r,\phi),s,t)
\\ &=& h(\tau(r,\phi),s,t,\alpha)\text{.}
\end{eqnarray*}
Since the coordinate transformation is surjective, this shows the claimed coincidence $B \circ f = H$. 
\endofproof

With Lemma \ref{lem3} we have now extracted a well-defined map
\begin{equation}
\label{3}
\alpha_{F}: TX \times_X TX \to \mathfrak{h} : (x,v_1,v_2) \mapsto -  \left . \frac{\partial^2F_\Gamma}{\partial s \partial t}\right |_{(0,0)} 
\end{equation}
from the 2-functor $F$.
\begin{lemma}
\label{lem1}
The map $\alpha_{F}$ has the following properties:
\begin{enumerate}
\item[(a)]
for fixed $x\in X$, it is antisymmetric and bilinear.
\item[(b)]
it is smooth.
\end{enumerate}
\end{lemma}

\proof
To prove (a), let $\bar\Gamma(s,t) := \Gamma(t,s)$, and let $F_\Gamma$ and $F_{\bar\Gamma}$ the corresponding smooth maps (\ref{9}). Due to the permutation, the derivatives (\ref{26}) yield the values for $\alpha_F(x,v_1,v_2)$ and $\alpha_F(x,v_2,v_1)$, respectively. Note that $\bar\Gamma_{*} \circ \Sigma_{\R} = \Gamma_{*} \circ \Sigma_{\R}^{-1}$, where $\Sigma_{\R}^{-1}$ is the 2-morphism inverse to the 2-morphism (\ref{25}) under vertical composition. Since the 2-functor $F$ sends inverse 2-morphisms to inverse group elements, we have $F_{\bar\Gamma} = F^{-1}_{\Gamma}$. Hence, by taking derivatives, we get $\alpha_F(x,v_2,v_1)=-\alpha_F(x,v_1,v_2)$. 

It remains to show that $\alpha_F(v_1 + \lambda v_1',v_2) = \alpha_F(v_1,v_2) + \lambda\alpha_F(v_1',v_2)$. If $\Gamma$ and $\Gamma'$ are smooth functions for the tangent vectors $(v_1,v_2)$ and $(v_1',v_2)$, respectively, we use a chart $\phi: U \to X$ of a neighbourhood of $x$ with $\phi(0)=x$ and construct a smooth function $\tilde\Gamma:(-\epsilon,\epsilon)^2 \to X$ by
\begin{equation*}
\tilde\Gamma (s,t) :=\phi( \phi^{-1}(\Gamma(s,t)) + \lambda \phi^{-1}\Gamma' (s,t))
\end{equation*}
where $\epsilon$ has to chosen small enough.  It is easy to see that $\tilde\Gamma(0,0)=p$ and that
\begin{equation}
\label{82}
v_1 + \lambda v_1' = \left . \frac{\mathrm{d}}{\mathrm{d}s} \right |_{s=0} \tilde \Gamma(s,0)
\quad\text{ and }\quad
v_2 = \left . \frac{\mathrm{d}}{\mathrm{d}t} \right |_{t=0} \tilde \Gamma(0,t)\text{.}
\end{equation}
On the other hand,
\begin{equation*}
\tilde \Gamma_{*}(\Sigma_{\R}(s,t)) = \phi_{*}(\phi^{-1}_{*}(\Gamma_{*}(\Sigma_{\R}(s,t))) + \lambda\phi^{-1}_{*}(\Gamma'_{*}(\Sigma_{\R}(s,t))))\text{,}
\end{equation*}
where we have used that bigons in $\R^2$ can be added and multiplied with scalars. This shows that
\begin{equation*}
 \left . \frac{\partial^2F_{\tilde\Gamma}}{\partial s \partial t}\right |_{(0,0)} =\left . \frac{\partial^2F_\Gamma}{\partial s \partial t}\right |_{(0,0)} +\lambda\left . \frac{\partial^2F_{\Gamma'}}{\partial s \partial t}\right |_{(0,0)}\text{;}
\end{equation*}
together with (\ref{82}) this proves that $\alpha_F$ is bilinear.

To prove (b), let $\phi:U \to X$ be a chart of $X$ with an open subset $U \subset \R^n$. The induced chart  $\phi_{TX}: U \times \R^n \to TX$ of the tangent bundle sends a point $(u,v) \in U \times \R^n$ to $\mathrm{d}\phi|_u(v)\in T_{\phi(u)}X$. We show that
\begin{equation*}
\alxydim{@C=1.5cm@R=0cm}{U \times \R^n \times \R^n \ar[r]^{\phi_{TX}^{[2]}} & TX \times_X TX \ar[r]^-{\alpha_{F}}  & \mathfrak{h}}
\end{equation*}
is a smooth map. For this purpose, let $c:U \times \R^n \times \R^n \times \R^2 \to BX$ be defined by $c(x,v_1,v_2,\sigma,\tau)(s,t) := \phi(u + \beta(s\sigma) v_1 + \beta(t \tau) v_2)$, where $\beta:[0,1] \to [0,1]$ is a smooth map with $\beta(0)=0$ and $\beta(1)=1$ and with sitting instants.  The map $c$ depends obviously smoothly on all parameters, so that $f_c := p_{2} \circ F \circ c: U \times \R^n \times \R^n \times \R^2 \to H$ is a smooth function. Notice that $\Gamma_{u,v_1,v_2}(s,t) := c(u,v_1,v_2,s,t)(1,1)$ defines a smooth map with the properties
\begin{equation*}
\Gamma(0,0) = \phi(u)
\quad\text{, }\quad
\left . \frac{\partial \Gamma}{ \partial s} \right |_{(0,0)}=\mathrm{d}\phi|_u(v_1)
\quad\text{ and }\quad
\left . \frac{\partial \Gamma}{ \partial t} \right |_{(0,0)}=\mathrm{d}\phi|_u(v_2)\text{.}
\end{equation*}
It is furthermore still related to $c$ by
\begin{equation*}
(\Gamma_{u,v_1,v_2})_{*}(\Sigma_{\R}(s,t)) = c(u,v_1,v_2,s,t)\text{.}
\end{equation*}
Now, 
\begin{eqnarray*}
(\alpha_F \circ \phi_{TX}^{[2]})(x,v_1,v_2) &=& \alpha_F(\phi(u),\mathrm{d}\phi|_u(v_1),\mathrm{d}\phi|_u(v_2))
\nonumber \\ &=&- \left . \frac{\partial^2 }{\partial s \partial t} \right |_{(0,0)} (p_2 \circ F \circ (\Gamma_{u,v_1,v_2})_{*} \circ \Sigma_{\R})(s,t)
\nonumber \\ &=&- \left . \frac{\partial ^2}{\partial s \partial t} \right |_{(0,0)} f_c(x,v_1,v_2,s,t)\text{.}
\end{eqnarray*}
The last expression is, in particular,  smooth in $x$, $v_1$ and $v_2$. 
\endofproof

All together, 
\begin{equation*}
B_x(v_1,v_2) := \alpha_{F}(x,v_1,v_2)
\end{equation*}
defines an $\mathfrak{h}$-valued 2-form $B\in \Omega^2(X,\mathfrak{h})$ on $X$, which is canonically associated to the smooth 2-functor $F$.

\begin{proposition}
\label{prop1}
Let $F:\mathcal{P}_2(X) \to \mathcal{B} \mathfrak{G}$ be a smooth  2-functor, and  $A\in \Omega^1(X,\mathfrak{g})$ and $B \in \Omega^2(X,\mathfrak{h})$ be the corresponding differential forms. Then,
\begin{equation}
\label{5}
\mathrm{d}A + [A
\wedge A] = t_{*} (B)\text{,}
\end{equation}
where $t_{*} :=\ \mathrm{d}t|_1: \mathfrak{h} \to \mathfrak{g}$ is the  Lie algebra homomorphism induced from the Lie group homomorphism $t$ which is part of the crossed module corresponding to the Lie 2-group $\mathfrak{G}$.
\end{proposition}

\proof
We consider again the bigon $\Sigma_\R(s,t)$ and the associated smooth map $F_\Gamma: \R^2 \to H$ from (\ref{9}). If we denote by $\gamma_1(s,t)$ the source path and by $\gamma_2(s,t)$ the target path of $\Sigma_{\R}(s,t)$, we obtain further smooth maps
\begin{equation*}
f_i := F \circ \Gamma_{*} \circ \gamma_i: \R^2 \to G\text{.}
\end{equation*}
Note that each of these two paths can be decomposed into  horizontal and  vertical paths, $\gamma_1(s,t)=\gamma^{v}_1(s,t) \circ \gamma_1^h(t)$ and $\gamma_2(s,t) = \gamma_2^h(s,t) \circ  \gamma_2^v(s)$, and that this decomposition induces accordant decompositions of the functions $f_i$, namely $f_1(s,t)= f^v_1(s,t) \cdot f^h_1(t)$ and $f_2(s,t)=f_2^h(s,t) \cdot f_2^v(s)$. We recall from (\ref{22}) that the 1-form $A \in \Omega^1(X,\mathfrak{g})$ is related to these functions by
\begin{equation*}
A_x(v_1) =-\left .  \frac{\partial f_i^v}{\partial s} \right |_{(0,0)}
\quad\text{ and }\quad
A_x(v_2) =-\left .  \frac{\partial f_i^h}{\partial t} \right |_{(0,0)}
\end{equation*}
for $i=1,2$. Now we employ the target-matching-condition  (\ref{10}) for the  2-morphism $F(\Sigma)$: 
\begin{equation}
\label{11}
f_2 = (t \circ F_{\Gamma}) \cdot f_1
\end{equation}
as functions from $\R^2$ to $G$. The second partial derivatives of the functions $f_i$ are at $(0,0)$:
\begin{eqnarray*}
\left . \frac{\partial^2f_1}{\partial s\partial t} \right |_{(0,0)}&=& \left . \frac{\partial ^2 f^v_1}{\partial s\partial t} \right|_{(0,0)} + A_x(v_1)A_x(v_2) \\
\left . \frac{\partial^2f_2}{\partial s\partial t} \right |_{(0,0)}&=& \left .  \frac{\partial^2 f^h_2}{\partial s\partial t} \right|_{(0,0)} +A_x(v_2)A_x(v_1)
\end{eqnarray*}
Here notice that $f_k(0,t)=1$ and $f_k(s,0)=1$ so that the second derivatives naturally take values in $\mathfrak{g}$. Also, we write $XY := \mathrm{d}m|_{(1,1)}(X,Y) \in \mathfrak{g}$ for $X,Y \in \mathfrak{g}$, which is -- in a faithful matrix representation -- just the product of matrices. The  first derivatives vanish,
\begin{equation}
\label{14}
 \left . \frac{\partial F_{\Gamma}}{\partial t} \right|_{(0,0)}
=
 \left . \frac{\partial F_{\Gamma}}{\partial s} \right|_{(0,0)}
=0\text{,}
\end{equation}
because $F$ is constant on families of identity bigons. Summarizing, equation (\ref{11})  becomes
\begin{equation*}
\left .  \frac{\partial ^2 f^h_2}{\partial s\partial t} \right|_{(0,0)}  +A_x(v_2) A_x(v_1)=-  \mathrm{d}t|_1 (B_x(v_1,v_2))+ \left . \frac{\partial ^2 f^v_1}{\partial s\partial t} \right|_{(0,0)}+A_x(v_1) A_x(v_2) \text{,}
\end{equation*}
this implies the claimed equality. 
\endofproof

\subsubsection{Extracting Forms II:  Pseudonatural Transformations}

\label{sec_ex2}

Now we discuss a smooth pseudonatural transformation
\begin{equation*}
\rho: F \to F'
\end{equation*}
between two smooth  2-functors $F,F': \mathcal{P}_2(X) \to \mathcal{B} \mathfrak{G}$. The components of $\rho$ are a smooth map $g: X \to G$ and a diffeological map  $\rho_1: P^{1}X \to G \ltimes H$. Notice that $\rho_1$ is not functorial, i.e. $\rho_1(\gamma_2 \circ \gamma_1)$ is in general not the product of $\rho_1(\gamma_2)$ and $\rho_1(\gamma_1)$. In order to remedy this problem, we construct another map $\tilde\rho:P^1X \to G \ltimes H$ from $\rho$ that will be functorial. We denote the projection of $\rho_1$ to $H$ by $\rho_H:= p_H \circ \rho_1: P^1X \to H$. Then we define
\begin{equation}
\label{61}
\tilde\rho(\gamma) := (F'(\gamma),\rho_{H}(\gamma)^{-1})\text{.}
\end{equation}
\begin{lemma}
\label{lem5}
$\tilde\rho$
defines a smooth functor $\tilde\rho:\mathcal{P}_1(X) \to \mathcal{B}(G \ltimes H)$.
\end{lemma}

\proof
For our convention concerning the semi-direct product, we refer the reader the equation (\ref{2}) in Appendix \ref{app2}. With this convention, axiom (T1) of the pseudonatural transformation $\rho$ infers for two composable paths $\gamma_1$ and $\gamma_2$ that
\begin{equation}
\label{49}
 \alpha(F'(\gamma_2),\rho_H(\gamma_1))\rho_H(\gamma_{2}) = \rho_H(\gamma_2\circ \gamma_1)\text{.}
\end{equation}
Then, the product of $\tilde\rho(\gamma_2)$ with $\tilde\rho(\gamma_1)$ in the semi-direct product $G \ltimes H$ is
\begin{eqnarray*}
&&\hspace{-3cm}(F'(\gamma_2),\rho_H(\gamma_2)^{-1}) \cdot (F'(\gamma_1),\rho_H(\gamma_1)^{-1}) \\&\stackrel{\text{(\ref{2})}}{=}& (F'(\gamma_2)F'(\gamma_1),\rho_H(\gamma_2)^{-1}\alpha(F'(\gamma_2),\rho_H(\gamma_1)^{-1}))
\\&\stackrel{\text{(\ref{49})}}{=}& (F'(\gamma_2\circ\gamma_1),\rho_H(\gamma_2 \circ\gamma_1)^{-1})\text{,}
\end{eqnarray*}
and thus equal to $\tilde\rho(\gamma_2 \circ \gamma_1)$. Since $F'(\id_x)=1$, equation (\ref{49}) also shows that $\tilde\rho(\id_x) = (1,1)$. Thus, $\tilde\rho$ is a functor. Its smoothness is clear from the definition. 
\endofproof

By Theorem \ref{th3}, the smooth functor $\tilde\rho$ corresponds to a 1-form with values in $\mathfrak{g}\ltimes \mathfrak{h}$, which in turn gives by projection into the two summands  an $\mathfrak{h}$-valued 1-form $\varphi\in\Omega^1(X,\mathfrak{h})$ and a $\mathfrak{g}$-valued 1-form. The ladder identifies (due to the definition of $\tilde\rho$) with the 1-form $A'$ that corresponds to the functor $F'$. Summarizing, the smooth pseudonatural transformation $\rho$ defines a smooth function $g:X\to G$ and a 1-form $\varphi\in\Omega^1(X,\mathfrak{h})$.

\begin{proposition}
\label{prop2}
Let $F,F': \mathcal{P}_2(X) \to \mathcal{B} \mathfrak{G}$ be smooth 2-functors with associated 1-forms $A,A'\in \Omega^1(X,\mathfrak{g})$ and 2-forms $B,B'\in \Omega^2(X,\mathfrak{h})$ respectively. The  smooth function $g:X \to G$ and the 1-form $\varphi\in \Omega^1(X,\mathfrak{h})$ extracted from a smooth pseudonatural transformation $\rho:F \to F'$ satisfy the relations 
\begin{eqnarray}
\label{6}
A' +t_{*} ( \varphi)&=& \mathrm{Ad}_g(A) - g^{*}\bar\theta    
\\
\label{7}
B' + \alpha_{*}(A' \wedge \varphi) + \mathrm{d}\varphi + [\varphi \wedge \varphi]&=& (\alpha_g)_{*} ( B)  \text{.}
\end{eqnarray}
In (\ref{6}), $\bar\theta$ is the right invariant Maurer-Cartan form  on $G$.
In (\ref{7}), $A' \wedge \varphi$ is a 2-form with values in $\mathfrak{h} \oplus \mathfrak{g}$, which is sent by the linear map $\alpha_{*}$ to a 2-form with values in $\mathfrak{h}$.
\end{proposition}

\proof
Like in the proof of Proposition \ref{prop1} we employ the target-matching condition  (\ref{10}) for the component
\begin{equation*}
\quadrat{\ast}{\ast}{\ast}{\ast}{F(\gamma)}{g(x)}{g(y)}{F'(\gamma)}{\rho(\gamma)}
\end{equation*}
 of the pseudonatural transformation $\rho$ at 1-morphism $\gamma:x \to y$ in $\mathcal{P}_2(X)$. For this purpose we choose a smooth curve $\Gamma: \R \to X$ through a point $x:=\Gamma(0)$ and consider the associated tangent vector $v\in T_xX$. With the standard path $\gamma_{\R}(t)$ in the real line from $0$ to $t$ we form from the 2-functors the smooth maps
\begin{equation*}
f := F \circ \Gamma_{*} \circ \gamma_{\R}:\R \to G
\quad\text{ and }\quad
f' := F' \circ \Gamma_{*} \circ \gamma_{\R}:\R \to G
\end{equation*}
and from the pseudonatural transformation the smooth maps 
\begin{equation}
\label{15}
\tilde g := \rho \circ \Gamma: \R \to G
\quad\text{ and }\quad
h := \rho_H \circ \Gamma_{*} \circ \gamma_{\R}:\R \to H\text{.}
\end{equation}
The condition we want to employ then becomes
\begin{equation}
\label{eq:hadtoshow}
f'(t) \cdot \tilde g(0) =t(h(t)) \cdot \tilde g(t) \cdot f(t)\text{.}
\end{equation}
If we take the definition of the function $g: X \to G$ and the 1-forms $A$, $A'$ and $\varphi$ into account, namely $g(x)=\tilde g(0)$ and
\begin{equation*}
A_x(v) =- \left . \frac{\partial f}{\partial t} \right|_{t=0}
\quad\text{, }\quad
A'_x(v) =- \left . \frac{\partial f'}{\partial t} \right|_{t=0}
\quad\text{ and }\quad
\varphi_x(v) =- \left . \frac{\partial}{\partial t} \right|_{t=0}h^{-1}\text{,}
\end{equation*}
the derivative of this equation at zero yields
\begin{equation*}
-A'_x(v) \cdot g(x) =  \mathrm{d}t|_1 ( \varphi_x(v)) \cdot g(x) +  \mathrm{d}g_x(v) - g(x) \cdot A_x(v) \text{,}
\end{equation*}
this implies  equation \erf{eq:hadtoshow} that we had to show. Here, the symbol $\cdot$ stands for the derivatives of left or right multiplication.

To prove the second equation we use axiom (T2) of the pseudonatural transformation $\rho$, namely the compatibility with 2-morphisms. For a 2-morphism $\Sigma$ in $\mathcal{P}_2(X)$, that we take of the form
\begin{equation*}
\quadrat{x_1}{y_1}{x_2}{y_2}{\gamma_1^h}{\gamma_2^v}{\gamma_1^v}{\gamma_2^h}{\Sigma}
\end{equation*}
this axiom requires
\begin{equation*}
\alxydim{@C=0.2cm@R=0.7cm}{& F(x_1) \ar[rr]^{F(\gamma_1^h)} \ar[dl]_{\rho(x_1)} && F(y_1) \ar@{=>}[dlll]|{\rho(\gamma_1^h)} \ar[dd]^{F(\gamma_1^v)} \ar[dl]|{\rho(y_1)} \\ F'(x_1) \ar[rr]|{F'(\gamma_1^h)} \ar[dd]_{F'(\gamma_2^v)}& & F'(y_1) \ar@{=>}[ddll]|{F'(\Sigma)} \ar[dd]|{F'(\gamma_1^v)} & \\ &&& F(y_2)\ar@{=>}[ul]|{\rho(\gamma_1^v)} \ar[dl]^{\rho(y_2)} \\ F'(x_2) \ar[rr]_{F'(\gamma_2^h)} && F'(y_2)}
=
\alxydim{@C=0.2cm@R=0.7cm}{& F(x_1) \ar[dd]|{F(\gamma_2^v)} \ar[rr]^{F(\gamma_1^h)} \ar[dl]_{\rho(x_1)} && F(y_1) \ar@{=>}[ddll]|{F(\Sigma)}  \ar[dd]^{F(\gamma_1^v)}  \\ F'(x_1)  \ar[dd]_{F'(\gamma_2^v)}& &    & \\ &F(x_2) \ar@{=>}[ul]|{\rho(\gamma_2^v)} \ar[dl]|{\rho(x_2)} \ar[rr]|{F(\gamma_2^h)}&& F(y_2) \ar@{=>}[dlll]|{\rho(\gamma_2^h)} \ar[dl]^{\rho(y_2)} \\ F'(x_2) \ar[rr]_{F'(\gamma_2^h)} && F'(y_2)}
\end{equation*}
With a choice of a smooth map $\Gamma: \R^2 \to X$ we can pullback these diagrams to $\R^2$ and use the standard bigon $\Sigma_{\R}(s,t)$. We use  the smooth functions $F_{\Gamma}$, $f_1$ and $f_2$ defined by the 2-functor $F$ as described in the proof of Proposition \ref{prop1}, and the analogous functions $F'_{\Gamma}$, $f_1'$ and $f_2'$ for the 2-functor $F'$. From the pseudonatural transformation $\rho$ we further obtain a function $\tilde g := \rho \circ \Gamma: \R^2 \to X$ and functions $h_i^h := \rho_H \circ \Gamma_{*} \circ \gamma_i^h$ and $h_i^v := \rho_H \circ \Gamma_{*} \circ \gamma_i^v$. Now we have
\begin{equation*}
\alxydim{@C=0.1cm@R=0.7cm}{& f(0,0) \ar[rr]^{f_1^h(t)} \ar[dl]_{\tilde g(0,0)} && f(0,t) \ar@{=>}[dlll]|{h_1^h(t)} \ar[dd]^{f_1^v(s,t)} \ar[dl]|{\tilde g(0,t)} \\ f'(0,0) \ar[rr]|{f_1^{\prime h}(t)} \ar[dd]_{f_2^{\prime v}(t)}& & f'(0,t) \ar@{=>}[ddll]|{F_{\Gamma}'(s,t)} \ar[dd]|{f_1^{\prime v}(s,t)} & \\ &&& f(s,t)\ar@{=>}[ul]|{h_1^v(s,t)} \ar[dl]^{\tilde g(s,t)} \\ f'(s,0) \ar[rr]_{f_2^{\prime h}(s, t)} && f'(s,t)}
=
\alxydim{@C=0.1cm@R=0.7cm}{& f(0,0) \ar[dd]|{f_2^v(s)} \ar[rr]^{f_1^h(t)} \ar[dl]_{g(0,0)} && f(0,t) \ar@{=>}[ddll]|{F_{\Gamma}(s,t)}  \ar[dd]^{f_1^v(s,t)}  \\ f'(0,0)  \ar[dd]_{f_2^{\prime v}(s)}& &    & \\ &f(s,0) \ar@{=>}[ul]|{h_2^v(s)} \ar[dl]|{\tilde g(s,0)} \ar[rr]|{f_2^h(s,t)}&& f(s,t) \ar@{=>}[dlll]|{h_2^h(s,t)} \ar[dl]^{\tilde g(s,t)} \\ f'(s,0) \ar[rr]_{f_2^{\prime h}(s,t)} && f'(s,t)}
\end{equation*}
Using the rules (\ref{12}) and (\ref{13}) for vertical and horizonal composition in $\mathcal{B} \mathfrak{G}$,  the above diagram boils down to  the equation
\begin{multline*}
F_{\Gamma}'(s,t)\cdot \alpha(f_1^{\prime v}(s,t),h_1^h(t)) \cdot h_1^v(s,t)
 \\= \alpha(f_2^{\prime h}(s,t),h_2^v(s)) \cdot h_2^h(s,t) \cdot \alpha(\tilde g(s,t),F_{\Gamma}(s,t))\text{.}
\end{multline*}
We now take the second mixed derivative and evaluate at $(0,0)$. 

For the evaluation we use the properties of the 2-functor $F$ that imply -- on the level of 2-morphisms -- $f(0,0)=1$ and -- on the level of 1-morphisms -- $f_1^h(0)=f_1^v(0,0)=f_2^v(0)=f_2^h(0,0)=1$. The same rules hold for  $F'$. Similarly, the properties of the functor $\tilde\rho$ give additionally $h_1^h(0)=h_1^v(0,0)=1$ and $h_2^v(0)=h_2^h(0,0)=1$. To compute the derivative of the terms that contain $\alpha$, it is convenient to use the rule
\begin{equation}
\label{33}
\mathrm{d}\alpha|_{(g,h)} (X,Y) =\mathrm{d}\alpha_{h}|_{g}(X) + \mathrm{d}\alpha_g|_{h}(Y)  
\end{equation}
where $\alpha_g: H \to H$ and $\alpha_h:G \to H$ are obtained from $\alpha$ by fixing one of the two parameters, and the differentials on the right hand side are taken only with respect to the remaining parameter. Finally, we use (\ref{14}). The result of the computation is (in notation introduced in the proof of Proposition \ref{prop1})
\begin{multline*}
\frac{\partial^2 F'_{\Gamma}}{\partial s\partial t}  + \mathrm{d}\alpha|_{(1,1)} \left ( \frac{\partial f_1'^v}{\partial s}, \frac{\partial h_1^h}{\partial t} \right )+ \frac{\partial h_1^h }{\partial t} \cdot \frac{\partial h_1^v}{\partial s} + \frac{\partial^2 h_1^v}{\partial s\partial t}
\\= \mathrm{d}\alpha|_{(1,1)}\left ( \frac{\partial f_2'^h}{\partial t}, \frac{\partial h_2^v}{\partial s} \right )  + \frac{\partial h_2^v}{\partial s} \cdot \frac{\partial h_2^h}{\partial t} + \frac{\partial^2 h_2^h}{\partial s\partial t} + \mathrm{d}\alpha_{\tilde g(0,0)}|_{1}\left ( \frac{\partial^2 F_{\Gamma}}{\partial s\partial t} \right )
\end{multline*}
Expressed by differential forms, this gives
\begin{multline*}
-B'_x(v_1,v_2) - \alpha_{*}(A'_x(v_1),\varphi_x(v_2))+\varphi_x(v_2)\varphi_x(v_1)+ \frac{\partial^2 h_1^v}{\partial s\partial t}
\\= -\alpha_{*} (A'_x(v_2),\varphi_x(v_1))+\varphi_x(v_1)\varphi_x(v_2) + \frac{\partial^2 h_2^h}{\partial s\partial t} - (\alpha_g)_{*} (B)\text{,}
\end{multline*}
which yields the second equality. 
\endofproof

\subsubsection{Extracting Forms III:  Modifications}

\label{sec_ex3}

Let us now consider a smooth modification
\begin{equation*}
\mathcal{A}:\rho_1 \Rightarrow \rho_2
\end{equation*}
between smooth pseudonatural transformations $\rho_1,\rho_2:F \to F'$ between two smooth 2-functors $F,F': \mathcal{P}_2(X) \to \mathcal{B}\mathfrak{G}$. Its components furnish a smooth map $X \to G \ltimes H$. We denote its projection on the second factor by $a:X \to H$.

\begin{proposition}
\label{prop3}
Let $F,F': \mathcal{P}_2(X) \to \mathcal{B} \mathfrak{G}$ be smooth 2-functors with associated 1-forms $A,A'\in\Omega^1(X,\mathfrak{g})$, let $\rho_1,\rho_2:F \to F'$ be smooth pseudonatural transformations with associated smooth functions $g_1,g_2: X \to G$ and 1-forms $\varphi_1,\varphi_2\in\Omega^1(X,\mathfrak{h})$. Then, the smooth map $a:X \to H$ associated to a smooth modification $\mathcal{A}: \rho_1 \Rightarrow \rho_2$ satisfies
\begin{equation}
\label{8}
g_2 = (t \circ a) \cdot g_1
\quad\text{ and }\quad
\varphi_2 +(r_{a}^{-1} \circ \alpha_{a})_{*}(A') =  \mathrm{Ad}_a (\varphi_1) -a^{*}\bar\theta\text{,}
\end{equation}
where $r_{a(x)}:H \to H$ is the multiplication with $a(x)$ from the right.
\end{proposition}

\proof
In the same way as before we choose a smooth map $\Gamma: \R \to X$ with $\Gamma(0)=: x$ and $\dot\Gamma(0)=:v \in T_{x}X$ and consider the smooth functions $f_{\Gamma},f_{\Gamma}': \R \to G$ from (\ref{9}), the smooth functions $\tilde g_1, \tilde g_2: \R \to G$ and $h_1,h_2: \R \to H$ from (\ref{15}), and define an additional smooth function $a_{\Gamma} := a \circ \Gamma: \R \to H$ with $a_{\Gamma}(0)=a(x)$. The target-matching condition (\ref{10}) for the 2-morphism
\begin{equation*}
\bigon{f_{\Gamma}(0)}{f'_{\Gamma}(0)}{\tilde g_1(0)}{\tilde g_2(0)}{a_{\Gamma}(0)}
\end{equation*}
in $\mathcal{B} \mathfrak{G}$ obviously gives us the first equation.  The axiom for the modification $\mathcal{A}$ implies
\begin{equation*}
\alpha(f'_{\gamma}(t),a_{\Gamma}(0)) \cdot h_1(t) = h_2(t) \cdot a_{\Gamma}(t)\text{.}
\end{equation*}
The first derivative evaluated at $0$ gives
\begin{equation*}
(\alpha_{a_{\Gamma}(0)})_{*}\left(\left .\frac{\partial f'_{\gamma}}{\partial t} \right|_0 \right)h_1(0) + \alpha(f'^{\gamma}(0),a_{\Gamma}(0)) \cdot \left . \frac{\partial h_1}{\partial t} \right |_{0} = \left . \frac{\partial h_2}{\partial t} \right |_0 \cdot a_{\Gamma}(0) + h_2(0) \cdot \left . \frac{\partial a_{\Gamma}}{\partial t} \right |_0
\end{equation*}
With $f'^{\gamma}(0)=h_{1}(0)=h_2(0)=1$ this yields
\begin{equation*}
(\alpha_{a(x)})_{*}(-A') + a(x) \cdot \varphi_1|_x(v) =  \varphi_2|_x(v) \cdot a(x) + \mathrm{d}a|_x(v)
\end{equation*}
which is the second equation we had to prove. 
\endofproof

\subsubsection{Summary of Section \ref{sec3_1}}

\label{sec2_2_4}

In order to obtain a precise relation between smooth 2-functors and differential forms, we define a 2-category which is adapted to the relations we have found in Propositions \ref{prop1}, \ref{prop2} and \ref{prop3}. 

\begin{definition}
\label{def8}
Let $\mathfrak{G}$ be a Lie 2-group, $(G,H,t,\alpha)$ the corresponding smooth
crossed module, and $X$ a smooth manifold. We define the following 2-category $\diffco{\mathfrak{G}}{2}{X}$:
\begin{enumerate}

\item 
 An object is a pair $(A,B)$ of a 1-form $A\in\Omega^1(X,\mathfrak{g})$
and a 2-form $B \in \Omega^2(X,\mathfrak{h})$ which satisfy the relation (\ref{5}):
\begin{equation*}
\mathrm{d}A + [A
\wedge A] = t_{*} ( B)\text{.}
\end{equation*}

\item
A 1-morphism
$(g,\varphi): (A,B) \to (A',B')$
is a smooth map $g:X \to G$ and a 1-form $\varphi\in \Omega^1(X,\mathfrak{h})$
that satisfy the relations (\ref{6}) and (\ref{7}):
\begin{eqnarray*}
A' + t_{*} (\varphi) &=& \mathrm{Ad}_g(A) - g^{*}\bar\theta   
\\
B' + \alpha_{*}(A' \wedge \varphi) + \mathrm{d}\varphi + [\varphi \wedge \varphi]&=& (\alpha_g)_{*} ( B)  \text{.}
\end{eqnarray*}
The composition of 1-morphisms
\begin{equation*}
\alxydim{@C=1.5cm}{(A,B) \ar[r]^-{(g_1,\varphi_1)} & (A',B') \ar[r]^-{g_2,\varphi_2} & (A'',B'')}
\end{equation*}
is given by the map $g_2g_1:X \to G$ and the 1-form $(\alpha_{g_{2}})_{*} (\varphi_1) +  \varphi_2$, where $\alpha_g:H \to H$ is the action of $G$ on $H$ with fixed $g$. The identity 1-morphism is given by $g=1$ and $\varphi=0$.

\item

A 2-morphism $a:(g_1,\varphi_1) \Rightarrow (g_2,\varphi_2)$ is a smooth map $a:X \to H$ that satisfies (\ref{8}):
\begin{equation*}
g_2 = (t \circ a) \cdot g_1
\quad\text{ and }\quad
\varphi_2 +(r_{a}^{-1} \circ \alpha_{a})_{*}(A') =  \mathrm{Ad}_a (\varphi_1) -a^{*}\bar\theta\text{.}
\end{equation*}
The vertical composition 
\begin{equation*}
\alxydim{}{(g,\varphi) \ar@{=>}[r]^-{a_1} & (g',\varphi') \ar@{=>}[r]^{a_2} & (g'',\varphi'')}
\end{equation*}
is given by $a_2a_1$. The horizontal composition is
\begin{equation*}
\alxydim{}{(A,B) \ar@/^2pc/[r]^{(g_1,\varphi_1)}="1" \ar@/_2pc/[r]_{(g_1',\varphi_1')}="2" \ar@{=>}"1";"2"|{a_1} & (A',B') \ar@/^2pc/[r]^{(g_2,\varphi_2)}="1" \ar@/_2pc/[r]_{(g_2',\varphi_2')}="2" \ar@{=>}"1";"2"|{a_2} & (A'',B'')}
=
\alxydim{@C=2.3cm}{(A,B) \ar@/^2pc/[r]^{(g_2g_1, (\alpha_{g_2})_{*}(\varphi_1) + \varphi_2)}="1" \ar@/_2pc/[r]_{(g_2'g_1', (\alpha_{g_2'})_{*}(\varphi'_1) +\varphi'_2)}="2" \ar@{=>}"1";"2"|{a_2\alpha(g_2,a_1)} & (A'',B'')\text{,}}
\end{equation*}
and the identity 2-morphism is given by $a=1$.
\end{enumerate}
\end{definition}
It is straightforward to check that this definition gives indeed a 2-category. In the Sections \ref{sec_ex1}, \ref{sec_ex2} and \ref{sec_ex3} above  we have collected the structure of a 2-functor \label{not:d}
\begin{equation*}
\fo: \mathrm{Funct}^{\infty}(\mathcal{P}_2(X),\mathcal{B}
\mathfrak{G}) \to \diffco{\mathfrak{G}}{2}{X}\text{.}
\end{equation*}  
Let us check that the axioms of a 2-functor are satisfied. Horizontal and vertical compositions of 2-morphisms are respected because these are just smooth maps $a: X \to H$ which become multiplied in exactly the same way in both 2-categories. It remains to check the compatibility with the composition of 1-morphisms, i.e. we have to show  that
\begin{equation*}
\fo(\rho_2 \circ \rho_1) = \fo(\rho_2) \circ \fo(\rho_1)
\end{equation*}
for smooth pseudonatural transformations $\rho_1:F \to F'$ and $\rho_2:F' \to F''$. Let $(g_i,\varphi_i) := \fo(\rho_i)$ for $i=1,2$. According to the definition (\ref{28}) of the  composition of pseudonatural transformations, the component of $\rho_2 \circ \rho_1$ at an object  $x\in X$ is $g_2(x)g_1(x)\in G$, and its component at a 1-morphism $\gamma:x \to y$ is $\rho_2(\gamma)\cdot \alpha(g_2(y),\rho_1(\gamma)) \in H$. If we consider the smooth functions $\tilde g_1, \tilde g_2:\R \to G$ and $h_1,h_2: \R \to H$ associated to $\rho_1$ and $\rho_2$ like in (\ref{15}), the 1-form associated to $\rho_{2} \circ \rho_1$ is, at $x:=\Gamma(0)$ and $v:=\dot\Gamma(0)$ and using (\ref{33}),
\begin{eqnarray*}
\nonumber
&&\hspace{-0.5cm}-\left . \frac{\mathrm{d}}{\mathrm{d}t} \right |_{0} \alpha(\tilde g_2(t),h_1(t)^{-1})h_2(t)^{-1}
\\ &&\hspace{1cm}= - \mathrm{d}\alpha_{\tilde g_2(0)}|_{h_1(0))} \left (  \left . \frac{\partial h^{-1}_{1}}{\partial t} \right|_0 \right )  h_2(0)^{-1}-\alpha(\tilde g_2(0),h_1(0)^{-1}) \left . \frac{\partial h^{-1}_{2}}{\partial t} \right |_{0}  \nonumber
 \\&&\hspace{1cm}=  (\alpha_{g_2(x)})_{*} ( \varphi_1|_x(v))+ \varphi_2|_x(v) \text{,}
\end{eqnarray*}
this is exactly the rule for horizontal composition of 1-morphisms in $\diffco{\mathfrak{G}}{2}{X}$.    

\subsection{From Forms to Functors}

\label{sec3_2}

In this section we introduce a 2-functor
\begin{equation*}
\fu: \diffco{\mathfrak{G}}{2}{X} \to \mathrm{Funct}^{\infty}(\mathcal{P}_2(X),\mathcal{B}
\mathfrak{G})
\end{equation*}
that goes in the  direction opposite to the 2-functor $\fo$ defined in Section \ref{sec3_1}. The principle here is to pose initial value problems governed by differential forms. Their unique  solutions  define smooth 2-functors, smooth pseudonatural transformations and smooth modifications.

\subsubsection{Reconstruction I: 2-Functors}

\label{sec3_2_1}

Here we consider a given 1-form $A\in \Omega^1(X,\mathfrak{g})$ and a given 2-form $B\in \Omega^2(X,\mathfrak{h})$ that satisfy the condition from Proposition \ref{prop1},
\begin{equation}
\label{39}
\mathrm{d}A + [A \wedge A]=t_{*}( B)\text{.} 
\end{equation}
By Theorem \ref{th3} the 1-form $A$ defines a smooth functor $F_A: \mathcal{P}_1(X) \to \mathcal{B} G$. Our aim is now to define a map $k_{A,B}:B^2X \to  H$ such that $F_A$ and $k_{A,B}$ together define a smooth 2-functor $F:\mathcal{P}_2(X) \to \mathcal{B} \mathfrak{G}$, which is dedicated to be the image  of the pair $(A,B)$ under the 2-functor $\fu$ we want to define.
For our convention concerning the semi-direct product, we refer the reader again to equation (\ref{2}) in Appendix \ref{app2}.

In order to find the correct definition of $k_{A,B}$ we look at the target-matching condition
\begin{equation}
\label{38}
F_A(\gamma_1) =  t(k_{A,B}(\Sigma)) \cdot F_A(\gamma_{0})
\end{equation}
that has to be satisfied for any bigon $\Sigma:\gamma_0 \Rightarrow \gamma_1$. For technical reasons we consider $\Sigma: [0,1]^2 \to X$ to be extended trivially over all of $\R^2$, i.e. 
\begin{equation*}
\Sigma(s,t)=\begin{cases}\gamma_0(0)=\gamma_1(0) & \text{ for }t<0 \\
\gamma_0(1)=\gamma_1(1) & \text{ for }t>1 \\
\gamma_0(t) & \text{ for }s<0\text{ and }0\leq t\leq1 \\
\gamma_1(t) & \text{ for }s>1\text{ and }0\leq t\leq1\text{.}
\end{cases}
\end{equation*}
Let $\tau_{s_0}(s,t)$ be the closed path in $\R^2$ that runs counter-clockwise around the rectangle spanned by $(s_0,0)$ and $(s_0 + s,t)$, and let the smooth function  $u_{A,s_0}: \R^2 \to G$ be defined by $u_{A,s_0}(s,t) := F_{A}(\Sigma_{*}(\tau_{s_0}(s,t)))$. For this function, we recall

\begin{lemma}[Lemma B.1 in \cite{schreiber3}]\
\label{lem7}
\begin{enumerate}
\item[(a)]
$u_{A,0}(1,1)=F_A(\gamma_{0}^{-1}
\circ \gamma_{1})$

\item[(b)]
$u_{A,s_0}(s,1)=u_{A,s_0}(s',1) \cdot u_{A,s_0+s'}(s-s',1)$

\item[(c)]
$\displaystyle\left .    \frac{\partial}{\partial s}
    \frac{\partial}{\partial t}
    u_{A,s_0} \right |_{(0,t)}
    =-
    \mathrm{Ad}_{F_A(\gamma_{s_0,t})}^{-1} 
     \left (\Sigma^* K \right )_{(s_0,t)} \left (\frac{\partial}{\partial s},\frac{\partial}{\partial t} \right)$
\end{enumerate}
with $\gamma_{s,t}$ the path defined by  $\gamma_{s,t}(\tau) := \Sigma(s,\tau t)$ and
 $K$ the curvature 2-form $K := \mathrm{d}A +[A \wedge A] \in \Omega^2(X,\mathfrak{g})$.
\end{lemma}

The function $u_{A,s_0}$ is interesting for us because by (a)  $u_{A,0}(1,1)$ coincides up to conjugation  with the image of the group element $k_{A,B}(\Sigma)\in  H$ we want to determine under the homomorphism $t$. The multiplicative property (b) shows that the smooth function $f:\R \to G$ defined by $f(\sigma) :=  u_{A,0}(\sigma,1)$ solves the initial value problem
\begin{equation}
\label{30}
\frac{\partial}{\partial \sigma} f(\sigma) = \mathrm{d}l_{f(\sigma)}|_1  \left ( \left . \frac{\partial}{\partial s} \right|_0 u_{A,\sigma}(s,1) \right )
\quad\text{and}\quad
f(0)=1\text{.}
\end{equation}
This initial value problem is governed by the 1-form
\begin{multline}
\label{35}
\left . \frac{\partial}{\partial s} \right|_0 u_{A,\sigma}(s,1) = \int_0^1 \mathrm{d}t \left \lbrace \frac{\partial}{\partial s} \left .\frac{\partial}{\partial t} u_{A,\sigma} \right|_{(0,t)} \right \rbrace \\\stackrel{\mathrm{(c)}}{=}-  \int_0^1 \mathrm{d}t \; \mathrm{Ad}^{-1}_{F_A(\gamma_{\sigma,t})} \left ( (\Sigma^{*}K)_{(\sigma,t)} \left ( \frac{\partial}{\partial s},\frac{\partial}{\partial t} \right ) \right )\text{.}
\end{multline}
Here, $\mathrm{Ad}^{-1}_{F_A(\gamma_{-,-})} \circ \Sigma^{*}K$ is a $\mathfrak{g}$-valued 2-form on $[0,1]^2$, and we have just performed a fibre integration over the second factor $[0,1]$. The result is a $\mathfrak{g}$-valued 1-form on $[0,1]$. This form actually lies in the image of $t_{*}$, 
\begin{equation}
\label{37}
(t \circ (\alpha_{F_A(\gamma_{-,-})^{-1}})_{*} (\Sigma^{*}B)\; \stackrel{\mathrm{(\ref{39})}}{=}\;\mathrm{Ad}^{-1}_{F_A(\gamma_{-,-})} (\Sigma^{*}K) \text{.}
\end{equation}
We are thus forced to consider the 1-form
\begin{equation}
\label{44}
\mathcal{A}_{\Sigma} :=- \int_{[0,1]} (\alpha_{F_A(\gamma_{-,-})^{-1}})_{*} (\Sigma^{*}B) \;\;\in \Omega^1([0,1],\mathfrak{h}) \text{.}
\end{equation}
Due to the sitting instants of $\Sigma$, we can equivalently speak of a 1-form on $\R$ which vanishes outside of $[0,1]$.
Now we use again Theorem \ref{th3} and obtain a smooth functor $F_{\mathcal{A}_{\Sigma}}:P^1\R \to H$. Since $P^1\R$ can be identified with $\R \times \R$ (compare Lemma 4.1 in \cite{schreiber3})  this is just a smooth function $f_{\Sigma}:\R^2 \to H$. The purpose of these definitions is, that by  (\ref{35}) and (\ref{37}) the smooth function
\begin{equation*}
f:\R \to G : \sigma \mapsto t(f_{\Sigma}(0,\sigma))^{-1}
\end{equation*}
 solves the initial value problem (\ref{30}). Thus, by uniqueness $t (f_{\Sigma}(0,\sigma))^{-1}=  u_{A,0}(\sigma,1)$. If we now define
\begin{equation}
\label{40}
k_{A,B}:BX \to H : \Sigma \mapsto \alpha(F_A(\gamma_0),f_{\Sigma}(0,1)^{-1})
\end{equation}
for $\gamma_0$ the source path of the bigon $\Sigma$ we have achieved
\begin{equation}
t(k_{A,B}(\Sigma)) = F_A(\gamma_0)\cdot t(f_{\Sigma}(0,1))^{-1} \cdot F_A(\gamma_0)^{-1}\stackrel{\mathrm{(a)}}{=} F_A( \gamma_{1})\cdot F_A(\gamma_0)^{-1}\text{;}
\end{equation}
this is the required target-matching condition (\ref{38}). 
Another indication that the map $k_{A,B}$ is we have found is the correct one is the following

\begin{proposition}
\label{prop4}
The map $k_{A,B}:BX \to H$ is diffeological. For any smooth map $\Gamma:\R^2 \to X$ with $x:=\Gamma(0)$, $v_1:=\frac{\partial\Gamma}{\partial s}$ and $v_2:=\frac{\partial\Gamma}{\partial t}$, we have
\begin{equation*}
 -  \left . \frac{\partial^2}{\partial s \partial t}\right |_{(0,0)} k_{A,B}(\Gamma_{*}\Sigma_{\R}(s,t)) = B_x(v_1,v_2)\text{.}
\end{equation*}
\end{proposition}

\proof
Assume that $c: U \to BX$ is a map from an open subset $U \subset \R^n$ to $BX$ such that the evaluation $c(u)(s,t) \in X$ is smooth on all of $U \times [0,1]^2$. Hence, the differential form $\mathcal{A}_{c(u)}$ from (\ref{44}) depends smoothly on $u\in U$, and so does the solution $f_{c(u)}:\R^2 \to H$ of the differential equation governed by $\mathcal{A}_{c(u)}$. This implies that $k_{A,B} \circ \mathrm{pr} \circ c: U \to H$ is smooth, so that $k_{A,B}$ is diffeological by Lemma \ref{lem2}. 

Now we consider $U=\R^2$ and $c := \Gamma_{*} \circ \Sigma_{\R}$ the standard bigon (\ref{25}), so that $k_{A,B} \circ c:\R^2 \to H$ is a smooth map. In order to compute the derivative
\begin{equation*}
\left . \frac{\partial}{\partial s} \right|_{0} k_{A,B}(c(s,t)) = \left . \frac{\partial}{\partial s} \right |_{0} \alpha(F_A(\Gamma_{*}\gamma_0(s,t)), f_{c(s,t)}(0,1)^{-1})
\end{equation*}
we observe that $\mathcal{A}_{c(s,t)}|_\sigma = \sigma \mathcal{A}_{c(1,t)}|_{\sigma s}$. For the solutions of the corresponding differential equations we obtain by a uniqueness argument $f_{c(s,t)}(0,\sigma) = f_{c(1,t)}(0,s\sigma)$. We  compute
\begin{eqnarray}
\hspace{-0.7cm}\left . \frac{\partial}{\partial s} \right |_0 f_{c(s,t)}(0,1)^{-1} &=& - \left . \frac{\partial}{\partial s}\right|_0 f_{c(1,t)}(0,s)
\nonumber\\&=& \mathcal{A}_{c(1,t)}|_{0}\left ( \frac{\partial}{\partial s} \right ) \nonumber\\&=& - \int_0^1 \mathrm{d}\tau \left ( \alpha_{F_A(\Gamma_{*}\gamma_{0,t\tau})} \right )_{*} (c(1,t)^{*}B)_{(0,\tau)}\left ( \frac{\partial}{\partial s}, \frac{\partial}{\partial \tau} \right ) \nonumber\\&=& -\int_0^t \mathrm{d}\tau'   \left ( \alpha_{F_A(\Gamma_{*}\gamma_{0,\tau'})} \right )_{*} (\Gamma^{*}B)_{(0,\tau')}\left ( \frac{\partial}{\partial s}, \frac{\partial}{\partial \tau'} \right )  
\label{47}
\end{eqnarray}
In the last step we have performed an integral transformation and used that $c(1,1)=\Gamma$. Finally
\begin{eqnarray*}
\left . \frac{\partial^2}{\partial s\partial t} \right |_0 k_{A,B}(c(s,t)) &=& \left . \frac{\partial}{\partial t} \right |_0 (\alpha_{F_A(\Gamma_{*}\gamma_0(0,t))})_{*} \left ( \left . \frac{\partial}{\partial s} \right |_0 f_{c(s,t)}(0,1)^{-1} \right )\nonumber \\&=& \left . \frac{\partial^2}{\partial s\partial t} \right |_0 f_{c(s,t)}(0,1)^{-1}  
\nonumber\\&\stackrel{\mathrm{(\ref{47})}}{=}&- (\Gamma^{*}B)|_0\left ( \frac{\partial}{\partial s}, \frac{\partial}{\partial t} \right ) 
\nonumber\\&=&- B_x(v_1,v_2)\text{.}
\end{eqnarray*}
In the first line we have used (\ref{33}) and that $f_{c(0,t)}(0,1)=1 \in H$, so that the differential of $\alpha_1: G \to H$ is the zero map. 
\endofproof

The next thing we would like to know about the map $k_{A,B}$ is its compatibility with the horizontal and vertical composition of bigons in $X$. Concerning the vertical composition, this will be straightforward, but for the horizontal composition we have to introduce firstly an auxiliary horizontal composition and to check the compatibility of $k_{A,B}$ with this one. 

\begin{figure}[h]
\begin{center}
\includegraphics{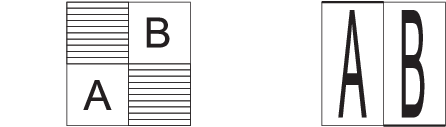}\setlength{\unitlength}{1pt}\begin{picture}(0,0)(245,729)\put(31.21455,754.01546){$p:$}\put(146.41455,754.01546){$\mapsto$}\end{picture}
\end{center}
\caption{A useful deformation of the unit square. }
\label{fig1}
\end{figure}

To define this auxiliary horizontal composition, we consider two bigons $\Sigma_1: \gamma_1 \Rightarrow \gamma_1'$ and $\Sigma_2:\gamma_2 \Rightarrow \gamma_2'$, with $\gamma_1,\gamma_1':x \to y$ and $\gamma_2,\gamma_2':y \to z$. The result will be a bigon
\begin{equation*}
\Sigma_2 \ast \Sigma_1 : \gamma_2 \circ \gamma_1 \Rightarrow \gamma_2' \circ \gamma_1'\text{.}
\end{equation*}
We define a map $p: [0,1]^{2} \to [0,1]^2$ by
\begin{equation*}
p(s,t) := \begin{cases} (0,t) & \text{ for } 0 \leq t < \frac{1}{2}\text{ and }0\leq s < \frac{1}{2} \\
(2s,t) & \text{ for }\frac{1}{2}\leq t \leq 1\text{ and }0\leq s < \frac{1}{2} \\
(2s-1,t) & \text{ for }0\leq t< \frac{1}{2}\text{ and }\frac{1}{2}\leq s \leq 1 \\
(1,t) & \text{ for }\frac{1}{2}\leq t \leq 1\text{ and }\frac{1}{2}\leq s\leq 1 \\
\end{cases}
\end{equation*}
see Figure \ref{fig1}. 
This map $p$ is not smooth, but its composition with $\Sigma_2 \circ \Sigma_1$ is, due to the sitting instants of the bigons $\Sigma_1$ and $\Sigma_2$. We define
\begin{equation*}
\Sigma_2 \ast \Sigma_1 := (\Sigma_2 \circ \Sigma_1) \circ p
\end{equation*}  
to be this smooth map, this defines the auxiliary horizontal composition of $\Sigma_1$ and $\Sigma_2$. 

\begin{lemma}
\label{lem6}
The map $k_{A,B}:BX \to H$ respects the vertical composition of bigons in the sense that
\begin{equation*}
k_{A,B}(\id_{\gamma})=1
\quad\text{ and }\quad
k_{A,B}(\Sigma_2\bullet\Sigma_1)= k_{A,B}(\Sigma_2) \cdot k_{A,B}(\Sigma_1)
\end{equation*}
for any path $\gamma\in PX$ and any two vertically composable bigons $\Sigma_1$ and $\Sigma_2$.
It respects the auxiliary horizontal composition $\ast$ in the sense that
\begin{equation*}
k_{A,B}(\Sigma_2 \ast \Sigma_1)=k_{A,B}(\Sigma_2) \cdot \alpha(F_A(\gamma_2),k_{A,B}(\Sigma_1))
\end{equation*}
for any two horizontally composable bigons $\Sigma_1:\gamma_1\Rightarrow\gamma_1'$ and $\Sigma_2:\gamma_2 \Rightarrow \gamma_2'$.
\end{lemma}

\proof
Concerning the vertical composition, the identity bigon $\id_{\gamma}:\gamma\Rightarrow\gamma$ has the 1-form $\mathcal{A}_{\id_{\gamma}} = 0$, so that $f_{\Sigma}(0,\sigma)$ is constant. Hence, $k_{A,B}(\id_{\gamma})=1$. Now let $\Sigma_1:\gamma_0 \Rightarrow \gamma_1$ and $\Sigma_2:\gamma_1 \Rightarrow \gamma_2$ be two bigons. For the 1-form (\ref{44}) associated to the bigon $\Sigma_2 \bullet \Sigma_1$ we find 
\begin{equation*}
\frac{1}{2}\mathcal{A}_{\Sigma_2 \bullet \Sigma_1}|_{\sigma} = \begin{cases} \mathcal{A}_{\Sigma_1}|_{2\sigma} & \text{ for } 0 \leq \sigma \leq \frac{1}{2} \\
\mathcal{A}_{\Sigma_2}|_{2\sigma-1} & \text{ for } \frac{1}{2}\leq \sigma \leq 1 \\
\end{cases}
\end{equation*} 
and accordingly
\begin{equation}
\label{41}
f_{\Sigma_2\bullet\Sigma_1}(0,1)=f_{\Sigma_2\bullet\Sigma_1}\left (\frac{1}{2},1 \right ) \cdot f_{\Sigma_2\bullet\Sigma_1}\left (0,\frac{1}{2} \right)= f_{\Sigma_2}(0,1) \cdot f_{\Sigma_1}(0,1)\text{.} 
\end{equation}
A short calculation then shows that
\begin{eqnarray*}
k_{A,B}(\Sigma_2\bullet\Sigma_1) &\stackrel{\mathrm{(\ref{40})}}{=}& \alpha(F_A(\gamma_0),f_{\Sigma_2\bullet\Sigma_1}(0,1)^{-1}) \nonumber
\\ &\stackrel{(\mathrm{\ref{41}})}{=}& \alpha(F_A(\gamma_0),f_{\Sigma_1}(0,1)^{-1} \cdot f_{\Sigma_2}(0,1)^{-1})\nonumber 
\\
&=& \alpha(F_A(\gamma_0) \cdot t(f_{\Sigma_1}(0,1))^{-1},f_{\Sigma_2}(0,1)^{-1} ) \\&&\hspace{5cm}\cdot \alpha(F_A(\gamma_0),f_{\Sigma_1}(0,1)^{-1} )
\\
&\stackrel{\mathrm{(a)}}{=}& \alpha(F_A(\gamma_1),f_{\Sigma_2}(0,1)^{-1} ) \cdot \alpha(F_A(\gamma_0),f_{\Sigma_1}(0,1)^{-1} )
\\&\stackrel{\text{\erf{40}}}{=}&  k_{A,B}(\Sigma_2) \cdot k_{A,B}(\Sigma_1)\text{.}
\end{eqnarray*}
In the step in the middle we have used the axioms of a crossed module, namely that
\begin{equation*}
\alpha(g,h_1h_
2)  = \alpha(g,\alpha(t(h_1),h_2) \cdot h_1) = \alpha(g\cdot t(h_1),h_2) \cdot \alpha(g,h_1)
\end{equation*}
for all $g\in G$ and $h_1,h_2\in H$.

Concerning the auxiliary horizontal composition, we obtain for the 1-form (\ref{44}) associated to the bigon $\Sigma_2 \ast \Sigma_1$
\begin{equation*}
\frac{1}{2}\mathcal{A}_{\Sigma_2 \ast \Sigma_1}|_{\sigma} = \begin{cases} \mathrm{Ad}^{-1}_{F_A(\gamma_1)}\left ( \mathcal{A}_{\Sigma_2}|_{2\sigma} \right ) & \text{ for } 0 \leq \sigma \leq \frac{1}{2} \\
\mathcal{A}_{\Sigma_1}|_{2\sigma-1} & \text{ for }\frac{1}{2}\leq \sigma\leq 1 \\
\end{cases}
\end{equation*}
and accordingly
\begin{multline}
\label{43}
f_{\Sigma_2 \ast \Sigma_1}(0,1)= f_{\Sigma_2 \ast \Sigma_1}\left(\frac{1}{2},1\right) \cdot f_{\Sigma_2\ast\Sigma_1}\left(0,\frac{1}{2}\right)\\= f_{\Sigma_1}(0,1) \cdot \alpha(F_A(\gamma_1)^{-1},f_{\Sigma_2}(0,1))\text{.}
\end{multline}
Then we obtain
\begin{eqnarray*}
k_{A,B}(\Sigma_2 \ast \Sigma_1)
&\stackrel{\mathrm{(\ref{40})}}{=}&
\alpha(F_A(\gamma_2 \circ \gamma_1), f_{\Sigma_2\ast\Sigma_1}(0,1)^{-1})
\nonumber
\\
&\stackrel{\mathrm{(\ref{43})}}{=}&
\alpha(F_A(\gamma_2 \circ \gamma_1), \alpha(F_A(\gamma_1)^{-1},f_{\Sigma_2}(0,1))^{-1} \cdot f_{\Sigma_1}(0,1)^{-1})
\nonumber
\\
&\stackrel{\mathrm{(\ref{40})}}{=}&
k_{A,B}(\Sigma_2) \cdot \alpha(F_A(\gamma_2),k_{A,B}(\Sigma_1))
\end{eqnarray*}
this yields the required identity.
\endofproof

Before we come to the original horizontal composition of bigons it is convenient to show first the following
\begin{lemma}
\label{lem4}
For thin homotopy equivalent bigons $\Sigma \sim_2 \Sigma'$ we have
\begin{equation*}
k_{A,B}(\Sigma) = k_{A,B}(\Sigma')\text{.}
\end{equation*} 
\end{lemma}

We have moved the proof of this lemma to Appendix \ref{app3}. Then it follows that $k_{A,B}$ factors through $B^2X$,
\begin{equation*}
\alxydim{}{BX \ar[r]^{\mathrm{pr}^2} & B^2X \ar[r] & H}
\end{equation*}
Since $\mathrm{pr}^2$ is surjective, the map $B^2X \to H$ is uniquely determined, and by Proposition \ref{prop4}, it is  diffeological. We denote this unique diffeological map also by $k_{A,B}:B^2X \to H$.

\begin{proposition}
\label{prop5}
The assignment
\begin{equation}
\label{48}
F \quad:\quad \bigon{x}{y}{\gamma}{\gamma'}{\Sigma} \quad\mapsto\quad \bigon{\ast}{\ast}{F_A(\gamma)}{F_A(\gamma')}{k_{A,B}(\Sigma)}
\end{equation}
defines a smooth 2-functor $F: \mathcal{P}_2(X) \to \mathcal{B} \mathfrak{G}$.
\end{proposition}

\proof
Since $F_A$ is a smooth functor, we have nothing to show for 1-morphisms. On 2-morphisms, the assignment $k_{A,B}$ is smooth by Proposition \ref{prop4}. By Lemma \ref{lem6} it further respects the vertical composition. Concerning the horizontal composition, notice that
\begin{equation*}
h:[0,1] \times [0,1]^2 \to X : (r,s,t) \mapsto (\Sigma_2 \circ \Sigma_2)(rp + (1-r)\id_{[0,1]^2})(s,t)\text{,}
\end{equation*}
defines a homotopy between $\Sigma_2 \ast \Sigma_1$ and $\Sigma_2 \circ \Sigma_1$, and since its rank is limited by dimensional reasons to 2, this  homotopy is thin.
Then, by Lemmata \ref{lem6} and \ref{lem4} we have
\begin{equation}
\label{42}
k_{A,B}(\Sigma_2 \circ \Sigma_1) =k_{A,B}(\Sigma_2 \ast \Sigma_1)= k_{A,B}(\Sigma_2) \cdot \alpha(F_A(\gamma_2),k_{A,B}(\Sigma_1))\text{.}
\end{equation}
Thus, the 2-functor $F$ respects the horizontal composition.
\endofproof

\subsubsection{Reconstruction II: Pseudonatural Transformations}

Here we consider a 1-morphism
\begin{equation*}
(g,\varphi): (A,B) \to (A',B')
\end{equation*}
in the 2-category $\diffco{\mathfrak{G}}{2}{X}$, i.e.
a smooth map $g:X \to G$ and a 1-form $\varphi\in \Omega^1(X,\mathfrak{h})$
that satisfy the relations from Proposition \ref{prop2},
\begin{eqnarray}
\label{45}
A'+ t_{*} ( \varphi)    &=& \mathrm{Ad}_g(A) - g^{*}\bar\theta 
\\
\label{51}
B' + \alpha_{*}(A' \wedge \varphi) + \mathrm{d}\varphi + [\varphi \wedge \varphi]&=& (\alpha_g)_{*} ( B)  \text{.}
\end{eqnarray}
The 1-forms $A'$ and $\varphi$ define a 1-form $(A',\varphi)\in \Omega^1(X,\mathfrak{g}\ltimes \mathfrak{h})$, and thus by Theorem \ref{th3} a smooth functor $\tilde\rho:\mathcal{P}_1(X) \to \mathcal{B}(G \ltimes H)$. We denote its projection to $H$ by $h: P^1X \to H$. 
We want to define a smooth pseudonatural transformation $\rho:F \to F'$ between the 2-functors $F := \fu(A,B)$ and $F' := \fu(A',B')$  by
\begin{equation}
\label{62}
\rho \quad:\quad
\alxy{ x \ar[r]^{\gamma} & y}
\quad\longmapsto\quad
\alxydim{@=1.2cm}{\ast \ar[r]^{F(\gamma)} \ar[d]_{g(x)} & \ast \ar[d]^{g(y)}
\ar@{=>}[dl]|{h(\gamma)^{-1}} \\ \ast \ar[r]_{F'(\gamma)} & \ast}\text{.}
\end{equation}
We have to show

\begin{lemma}
\label{lem8}
The target-matching condition
\begin{equation}
\label{50}
F'(\gamma)\cdot g(x) = t(h(\gamma)^{-1}) \cdot g(y) \cdot F(\gamma)
\end{equation}
for the 2-morphism $h(\gamma)^{-1}$ is satisfied.
\end{lemma}

\proof
We recall that $F(\gamma)$, $F'(\gamma)$ and $h(\gamma)$ are  values of solutions $f_{\gamma}, f_{\gamma}': \R \to G$ and $h_{\gamma}:\R \to H$ of  initial value problems. We show that
\begin{equation*}
f_{\gamma}'(0,t)= t(h_{\gamma}(t)^{-1})\cdot g(\gamma(t))\cdot f_{\gamma}(0,t) \cdot g(\gamma(0))^{-1} =: \beta(t)
\end{equation*}
which gives for $t=1$
 equation (\ref{45}). For this purpose, we show that $\beta(t)$ satisfies the initial value problem for $f_{\gamma}'$. The initial condition $\beta(0)=1$ is satisfied. Notice that  with $p:=\gamma(t)$ and $v := \dot \gamma(t)$
\begin{equation}
\label{55}
\frac{\partial h_{\gamma}(t)}{\partial t} = - \mathrm{d}r_{h_{\gamma}(t)}|_1 (\varphi_p(v)) -  (\alpha_{h_{\gamma}(t)})_{*}(A'_p(v))
\end{equation} 
so that -- using Axiom 2a) of the crossed module --
\begin{multline*}
\frac{\partial}{\partial t}t(h_{\gamma}(t)^{-1})=\mathrm{d}t|_{h_{\gamma}(t)^{-1}}\left( \frac{\partial h_{\gamma}(t)^{-1}}{\partial t} \right ) \\=  \left( \mathrm{Ad}^{-1}_{t(h_{\gamma}(t))} \left ( t_{*} (\varphi_{p}(v)) + A'_{p}(v) \right )- A'_{p}(v) \right ) \cdot t(h_{\gamma}(t)^{-1}) \text{.}
\end{multline*}
Then we compute
\begin{equation*}
\frac{\partial\beta}{\partial t} = \left ( \mathrm{Ad}^{-1}_{t(h_{\gamma}(t))} \left ( t_{*} (\varphi_{p}(v)) + A'_{p}(v)+ g^{*}\bar\theta|_p(v) - \mathrm{Ad}_{g}(A_p(v)) \right ) - A'_p(v) \right ) \cdot \beta(t)\text{.}
\end{equation*}
Using equation (\ref{45}), the right hand side becomes $-A'_p(v)\beta(t)$. Hence, $\beta(t)$ solves the same initial value problem as $f_{\gamma}'(0,t)$. By uniqueness, both functions coincide.
\endofproof

It remains to check that the axioms of a pseudonatural transformation are satisfied. Axiom (T1) follows from the fact that $\tilde\rho$ is a functor by the same arguments as given in the proof of Lemma \ref{lem5}. For axiom (T2) we have to prove

\begin{lemma}
The 2-morphism $h(\gamma)$ satisfies
\begin{equation*}
F'(\Sigma)\cdot h^{-1}(\gamma_0) = h^{-1}(\gamma_1)\cdot \alpha(g(y),F(\Sigma))
\end{equation*}
for any bigon $\Sigma:\gamma_0 \Rightarrow \gamma_1$.
\end{lemma}

\proof
 We recall that $F(\Sigma) = k_{A,B}(\Sigma) = \alpha(F(\gamma),f_{\Sigma}(0,1)^{-1})$, where $f_{\Sigma}(0,s)$ is the solution of a initial value problem governed by a 1-form $\mathcal{A}_{\Sigma}$. For $F(\gamma')$ the same is true with primed quantities. We define the notion $\gamma_s(t):= \Sigma(s,t)$ consistent with $\gamma_0$ and $\gamma_1$. Then, the equation
\begin{equation*}
f'_{\Sigma}(0,s)=\alpha(F'(\gamma_{0})^{-1},h(\gamma_{0})^{-1}\cdot \alpha(g(y)\cdot F(\gamma_{0}),f_{\Sigma}(0,s))\cdot h(\gamma_s)) := \kappa(s) \text{,}
\end{equation*} 
evaluated for $s=1$, is the equation we have to prove. Like in the proof of Lemma \ref{lem8} we show that $\kappa(s)$ solves the initial value problem for $f_{\Sigma}'(0,s)$. In a first step, the derivative $\partial\kappa/\partial s$ can be written as  $\mathrm{d}r_{\kappa(s)}|_1 X(s)$ where $X(s)\in \mathfrak{h}$ is
\begin{eqnarray*}
X(s)&=&(\alpha_{F'(\gamma_0)^{-1}})_{*} \left (- \mathrm{Ad}^{-1}_{h(\gamma_0)} (\alpha_{g(y)F(\gamma_0)})_{*} \mathcal{A}_{\Sigma}|_s\left(\frac{\partial}{\partial s} \right) \right . \\ && + \left . \mathrm{Ad}_{h(\gamma_0)^{-1}\alpha(g(y)F(\gamma_0),f_{\Sigma}(0,s))}\left ( \frac{\partial h(\gamma_s)}{\partial s}h(\gamma_s)^{-1} \right ) \right ) \\&=& -(\alpha_{g(x)})_{*} \left ( \mathcal{A}_{\Sigma}|_s \left ( \frac{\partial }{\partial s} \right ) \right ) +(\alpha_{F'(\gamma_s^{-1})} )_{*} h(\gamma_s)^{-1} \frac{\partial h(\gamma_s)}{\partial s}\text{.}
\end{eqnarray*}
In the second line we have used the target matching conditions (\ref{38}) and (\ref{50}). With the definition (\ref{44}) and again (\ref{50}), the first summand becomes
\begin{equation*}
-(\alpha_{g(x)})_{*}\mathcal{A}_{\Sigma}|_s\left ( \frac{\partial}{\partial s} \right ) = \int_{0}^1 \mathrm{d}t \,(\alpha_{F'(\gamma_{s,t})^{-1}})_{*} \mathrm{Ad}^{-1}_{h(\gamma_{s,t})} W(s,t)\text{.}
\end{equation*}
where we have written 
\begin{equation*}
W(s,t) := \Sigma^{*}((\alpha_g)_{*} (B))_{(s,t)}\left (\frac{\partial}{\partial s},\frac{\partial}{\partial t} \right ) \in \mathfrak{h}\text{.}
\end{equation*}
To compute the second summand, we recall from Section \ref{sec3_2_1} the definition of the path $\tau_{s_0}(s,t)$ that runs counter-clockwise around the rectangle spanned by $(s_0,0)$ and $(s_0 + s,t)$. We consider the smooth function  $u_{s_0}: \R^2 \to G \ltimes H$ be defined by $u_{s_0}(s,t) := \tilde\rho(\Sigma_{*}(\tau_{s_0}(s,t)))$, where $\tilde\rho$ is the smooth functor corresponding to the 1-form $(A',\varphi)\in \Omega^1(X,\mathfrak{g}\ltimes \mathfrak{h})$ we started with. For this smooth function, we recall Lemma \ref{lem7} (a), here $h(\gamma_{0}^{-1} \circ \gamma_s) =p_H(u_0(s,1))$. Furthermore, we have
\begin{eqnarray}
\frac{\partial}{\partial s}u_0(s,1) &\stackrel{\mathrm{(b)}}{=}& u_0(s,1) \cdot \left . \frac{\partial}{\partial \sigma} \right |_0 u_s(\sigma,1) \nonumber
\\ &=& u_0(s,1) \cdot  \int_0^1 \mathrm{d}t \left . \frac{\partial}{\partial \sigma} \frac{\partial}{\partial t} \right |_{(0,t)}  u_s(\sigma,t) \nonumber
\\ &\stackrel{\mathrm{(c)}}{=}&-u_0(s,1) \cdot \int_0^1 \mathrm{d}t\; \mathrm{Ad}^{-1}_{\tilde\rho(\gamma_{s,t})} (\Sigma^{*}K)_{(s,t)} \left (\frac{\partial}{\partial s},\frac{\partial}{\partial t} \right ) \label{52} \text{.}
\end{eqnarray}
In the last line, $K=(K_{A'},K_{\varphi})$ is the curvature 2-form of the 1-form $(A',\varphi)$, consisting of
\begin{equation*}
K_{A'} = \mathrm{d}A' + [A' \wedge A']\stackrel{\mathrm{(\ref{39})}}{=}t_{*}\circ B'
\quad\text{ and }\quad
K_{\varphi}= \alpha_{*}(A' \wedge \varphi) + \mathrm{d}\varphi + [\varphi \wedge \varphi]\text{.}
\end{equation*}
If we write $Y(s,t)$ for $\Sigma^{*}K_{\varphi}$ evaluated at $(s,t)$, and similarly $Z(s,t)$ for $\Sigma^{*}B'$, the  adjoint action in (\ref{52}) on the semidirect product $\mathfrak{g}\ltimes \mathfrak{h}$ is
\begin{equation*}
\mathrm{Ad}^{-1}_{(g,h)}(t_{*}(Z),Y)=\left(\mathrm{Ad}_{g}^{-1}(t_{*}(Z)),(\alpha_{g^{-1}})_{*}\left ( \mathrm{Ad}^{-1}_h(Y+Z) - Z \right) \right)\text{.}
\end{equation*}
With $Y + Z=W$ from (\ref{51}), the projection of (\ref{52}) to $\mathfrak{h}$  becomes
\begin{multline}
\label{53}
\frac{\partial h(\gamma_{0}^{-1} \circ \gamma_s)}{\partial s}=-h(\gamma_{0}^{-1} \circ \gamma_s)\cdot (\alpha_{F'(\gamma_{0}^{-1} \circ \gamma_s)} )_{*}\\\left( \int_0^1 \mathrm{d}t \, (\alpha_{F'(\gamma_{s,t})^{-1}})_{*}\left ( \mathrm{Ad}^{-1}_{h(\gamma_{s,t})}(W(s,t)) - Z(s,t)\right)  \right )
\end{multline}
Then, with $h(\gamma_0^{-1} \circ \gamma_s)= h(\gamma_0^{-1})\alpha(F'(\gamma_0)^{-1},h(\gamma_s))$, we have summarizing
\begin{equation*}
X(s) =\int_{0}^1 \mathrm{d}t (\alpha_{F'(\gamma_{s,t})^{-1}})_{*}  Z(s,t)  \\ \stackrel{\mathrm{(\ref{44})}}{=} -\mathcal{A}_{\Sigma}'|_s \left ( \frac{\partial}{\partial s} \right )\text{.}
\end{equation*}
This shows $\kappa(s) = f_{\Sigma}'(0,s)$.
\endofproof

\subsubsection{Reconstruction III: Modifications}

We consider a 2-morphism
\begin{equation*}
a:(g,\varphi) \Rightarrow (g',\varphi')
\end{equation*}
in the 2-category $\diffco{G}{2}{X}$, between two 1-morphisms $(g,\varphi)$ and $(g',\varphi')$ from $(A,B)$ to $(A',B')$. This is a smooth map $a:X \to H$ that satisfies (\ref{8}):
\begin{equation}
\label{54}
g_2 = (t \circ a) \cdot g_1
\quad\text{ and }\quad
\varphi_2 +(r_{a}^{-1} \circ \alpha_{a})_{*}(A') =  \mathrm{Ad}_a (\varphi_1) -a^{*}\bar\theta\text{.}
\end{equation}
We want to define a smooth modification $\mathcal{A}:\rho \Rightarrow \rho'$ between the pseudonatural transformations $\rho := \fu(g,\varphi)$ and $\rho ':= \fu(g',\varphi')$. We define
\begin{equation*}
\mathcal{A} \quad: \quad x \quad \mapsto \bigon{\ast}{\ast}{g(x)}{g'(x)}{a(x)}\text{.} \end{equation*}
The target-matching condition for the 2-morphism $f(x)$ is obviously satisfied due to the first equation in (\ref{54}). The axiom for the modification $\mathcal{A}$ is

\begin{lemma}
\label{lem11}
The 2-morphism $a(x)$ satisfies
\begin{equation*}
\alpha(F'(\gamma),a(x))\cdot h(\gamma)^{-1} = h'(\gamma)^{-1}\cdot a(y)
\end{equation*}
for all paths $\gamma \in PX$.
\end{lemma}

\proof
We rewrite the equation as
\begin{equation*}
h'_\gamma(t) = a(\gamma(t))\cdot h_{\gamma}(t) \cdot \alpha(f_{\gamma}'(0,t),a(x)^{-1}) := \lambda(t)
\end{equation*}
which we will  prove by showing that $\lambda(t)$ satisfies the same initial value problem as $h_{\gamma}'(t)$, namely (\ref{55}):
\begin{equation}
\label{56}
\frac{\partial h_{\gamma}'(t)}{\partial t} = - \mathrm{d}r_{h'_{\gamma}(t)}|_1 (\varphi'_p(v)) -  (\alpha_{h'_{\gamma}(t)})_{*}(A'_p(v))
\end{equation} 
for $p:=\gamma(t)$ and $v:=\dot \gamma(t)$. A straightforward calculation shows that
\begin{multline*}
\frac{\partial\lambda}{\partial t} =- \mathrm{d}r_{\lambda(t)}|_1\left (-(a^{*}\bar\theta)_p(v) + \mathrm{Ad}_{a(p)}(\varphi_{1}|_{p}(v)) \right . \\ \left . -(r_{a(p)}^{-1} \circ \alpha_{a(p)})_{*}(A'_p(v))  \right )  
 - (\alpha_{\lambda(t)})_{*}(A'_p(v))\text{.}
\end{multline*}
For this calculation, one twice has to use  the identity
\begin{equation}
\label{73}
(\alpha_{h_1h_2})_{*}(X)=
\mathrm{d}r_{h_2}|_{h_1}(\alpha_{h_1})_{*}(X) + \mathrm{d}l_{h_1}|_{h_2}(\alpha_{h_2})_{*}(X)\text{.}
\end{equation}
Using then the second equation of (\ref{54}), we have shown that $\lambda(t)$ satisfies the differential equation (\ref{56}). Thus, $\lambda(t)=h'_\gamma(t)$.
\endofproof

\subsubsection{Summary of Section \ref{sec3_1}}

Above we have collected the structure of a 2-functor
\begin{equation*}
\fu:\diffco{\mathfrak{G}}{2}{X} \to \mathrm{Funct}^{\infty}(\mathcal{P}_2(X),\mathcal{B}
\mathfrak{G}) \text{.}
\end{equation*}
Let us  check that the axioms of a 2-functor are satisfied. Like in Section \ref{sec2_2_4}, horizontal and vertical composition of 2-morphisms is respected because they are defined on both sides in the same way for the same $H$-valued functions. It remains to check the compatibility with the composition of 1-morphisms,
\begin{equation*}
\fu((g_2,\varphi_2) \circ (g_1,\varphi_1)) := \fu(g_2,g_1,(\alpha_{g_{2}})_{*} (\varphi_1) +  \varphi_2) = \fu(g_2,\varphi_2) \circ \fu(g_1,\varphi_1)
\end{equation*}
for 1-morphisms $(g_i,\varphi_i):(A_i,B_i) \to (A_{i+1},B_{i+1})$. For the components at objects $x\in X$, this equality is clear. We recall the component of $\fu(g_i,\varphi_i)$ at a path $\gamma\in PX$ is a 2-morphism in $\mathcal{B}\mathfrak{G}$ given according to (\ref{62}) by a group element $h(\gamma)^{-1}$, where $h(\gamma) = h_i(1)$ for $h_i(t)$ the solution of the initial value problem (\ref{55}). Similar, the component of $\fu(g_2,g_1,  \tilde\varphi)$ at $\gamma$ with $\tilde \varphi :=(\alpha_{g_{2}})_{*} \circ\varphi_1 +  \varphi_2$ is $\tilde h(\gamma)^{-1}$, where $\tilde h(\gamma)=\tilde h(1)$ for $\tilde h(t)$ the solution of the initial value problem
\begin{equation}
\label{71}
\frac{\partial \tilde h(t)}{\partial t} = - \mathrm{d}r_{\tilde h(t)}|_1 (\tilde\varphi_2|_{\gamma(t)}(v_{t})) -  (\alpha_{\tilde h(t)})_{*}(A_3|_{\gamma(t)}(v_{t}))
\end{equation}
with $v_t := \dot\gamma(t)$. According to the definition of the composition of pseudonatural transformations,
the equation we have to prove now follows from
\begin{equation}
\label{70}
\tilde h(t) = \alpha(g_2(\gamma(t)),h_1(t))\cdot h_2(t) =: \zeta(t)
\end{equation}
evaluated at $t=1$, and we prove (\ref{70}) by showing that $\zeta(t)$ solves (\ref{71}). A straightforward calculation similar to the one performed in the proof of Lemma \ref{lem11},  using (\ref{73}) and (\ref{45}) for $(g_2,\varphi_2)$, shows that this is indeed the case.

\subsection{Main Theorem}

\label{sec3_3}

We have so far defined two 2-functors $\fo$ and $\fu$ which go from smooth  2-functors to differential forms, and from differential forms back to smooth 2-functors.
Here we prove the main theorem of this article:

\begin{theorem}
\label{th2}
The 2-functors 
\begin{equation*}
\fo: \mathrm{Funct}^{\infty}(\mathcal{P}_2(X),\mathcal{B}
\mathfrak{G}) \to \diffco{\mathfrak{G}}{2}{X}
\end{equation*}
from Section \ref{sec3_1} and
\begin{equation*}
\fu: \diffco{\mathfrak{G}}{2}{X} \to \mathrm{Funct}^{\infty}(\mathcal{P}_2(X),\mathcal{B}
\mathfrak{G})
\end{equation*}
from Section \ref{sec3_2} satisfy
\begin{equation}
\label{58}
\fo \circ \fu = \id_{\diffco{G}{2}{X}}
\quad\text{ and }\quad
\fu \circ \fo = \id_{\mathrm{Funct}^{\infty}(\mathcal{P}_2(X),\mathcal{B}\mathfrak{G})}
\end{equation}
and form hence an isomorphism  of 2-categories.
\end{theorem}

\proof
We start with an object $(A,B)$ in $\diffco{\mathfrak{G}}{2}{X}$, i.e. a 1-form $A\in\Omega^1(X,\mathfrak{g})$ and a 2-form $B\in \Omega^2(X,\mathfrak{h})$ such that $\mathrm{d}A + [A
\wedge A] = t_{*} \circ B$. We let $(A',B'):=\fo(\fu(A,B))$ be the differential forms extracted from the reconstructed 2-functor $F:=\fu(A,B)$. By Theorem \ref{th3} we have $A'=A$. Now we  test the 2-form $B'$ at a point $x\in X$ and at tangent vectors $v_1,v_2\in T_xX$. Let $\Gamma:\R^2 \to X$ be a smooth map with $x=\Gamma(0)$, $v_1=\left . \frac{\partial\Gamma}{\partial s} \right |_0$ and $v_2=\left . \frac{\partial\Gamma}{\partial t} \right |_0$. We only have to summarize
\begin{multline*}
B'_x(v_1,v_2) \stackrel{\mathrm{(\ref{3})}}{=} - \left . \frac{\partial^2}{\partial s\partial t} \right |_{0} \fu(A,B)(\Gamma_{*}\Sigma_{\R}(s,t))\\\stackrel{\mathrm{(\ref{48})}}{=}-\left . \frac{\partial^2}{\partial s\partial t} \right |_{0} k_{A,B}(\Gamma_{*}\Sigma_{\R}(s,t)) = B_x(v_1,v_2)\text{,}
\end{multline*}
where the last equality  has been shown in Proposition \ref{prop4}.

Conversely, let $F:\mathcal{P}_2(X) \to \mathcal{B}\mathfrak{G}$ be a smooth 2-functor, and let $F' := \fu(\fo(F))$. By Theorem \ref{th3} it is clear that $F'(x) = F(x)$ and $F'(\gamma)=F(\gamma)$ for every point $x\in X$ and every path $\gamma\in PX$. For a bigon $\Sigma \in B^2X$ we recall that
\begin{equation}
\label{59}
F'(\Sigma) = k_{\fo(F)}(\Sigma)=\alpha(F(\gamma_0),f'_{\Sigma}(0,1)^{-1})\text{,}
\end{equation}
where $f'_{\Sigma}$ is the solution of the initial value problem
\begin{equation}
\label{60}
\frac{\partial f'_{\Sigma}(0,s)}{\partial s} = - \mathrm{d}r_{f'_{\Sigma}(0,s)} \left (X(s) \right)
\quad\text{ and }\quad
f'_{\Sigma}(0,0)=1\text{.}
\end{equation}
This initial value problem is governed by $X(s)\in \mathfrak{h}$, which is given by the 1-form $\mathcal{A}_{\Sigma}$ from (\ref{44}), namely
\begin{equation*}
X(s):=\mathcal{A}_{\Sigma}|_s\left(\frac{\partial}{\partial s}\right ) :=- \int_{0}^1 \mathrm{d}t\,(\alpha_{F(\gamma_{s,t})^{-1}})_{*}  (\Sigma^{*}B)_{(s,t)}\left( \frac{\partial}{\partial s},\frac{\partial}{\partial t}\right) \text{,}
\end{equation*}
and $B$ is the 2-form in $(A,B)=\fo(F)$.

We define a bigon $\Sigma_{s,t}(\sigma,\tau)$ by $\Sigma_{s,t}(\sigma,\tau)(s',t') := \Sigma(s+\beta(\sigma s'),t+\beta(\tau t'))$, where $\beta$ is some fixed smooth map $\beta: [0,1] \to [0,1]$ with $\beta(0)=0$ and $\beta(1)=1$ and with sitting instants. We notice from (\ref{59}) and (\ref{60}) that $F'(\Sigma_{0,0}(s,1))$ is the unique solution of the initial value problem
\begin{equation*}
\frac{\partial}{\partial s} F'(\Sigma_{0,0}(s,1))=F'(\Sigma_{0,0}(s,1)) \cdot  \mathrm{d}\alpha_{F(\gamma_0)}(X(s))
 \quad\text{ and }\quad F'(\Sigma_{0,0}(0,1))=1\text{.}
\end{equation*}
In the following we prove that $F(\Sigma_{0,0}(s,1))$ also solves this initial value problem, so that in particular 
\begin{equation*}
F'(\Sigma)=F'(\Sigma_{0,0}(1,1))=F(\Sigma_{0,0}(1,1))=F(\Sigma)
\end{equation*}
follows, and we have $\fu(\fo(F))=F$. To show that $F(\Sigma_{0,0}(s,1))$ is a solution we compute
\begin{equation*}
\frac{\partial}{\partial s}F(\Sigma_{0,0}(s,1)) = F(\Sigma_{0,0}(s,1)) \mathrm{d}\alpha_{F(\gamma_0)} \alpha_{F(\gamma_s)^{-1}} \left( \left .\frac{\partial}{\partial \sigma} \right|_0 F(\Sigma_{s,0}(\sigma,1)) \right)
\end{equation*}
and then
\begin{multline*}
\alpha_{F(\gamma_s)^{-1}}\left(\left .\frac{\partial}{\partial \sigma} \right|_0 F(\Sigma_{s,0}(\sigma,1)) \right)\\ = \int_0^1 \mathrm{d}t \,( \alpha_{F(\gamma_{s,t})^{-1}})_{*} \left . \frac{\partial^2}{\partial \sigma \partial \tau} \right|_0 \alpha( F(\gamma_{s,t+\tau})^{-1} , F(\Sigma_{s,0}(\sigma,t+\tau)))\text{.}
\end{multline*}
To compute the derivative we decompose $\Sigma_{s,0}(\sigma,t+\tau)$ in two bigons $\Sigma_{s,0}(\sigma,t)$ and $\Sigma_{s,t}(\sigma,\tau)$ and obtain
\begin{equation*}
\left . \frac{\partial^2}{\partial \sigma \partial \tau} \right|_0 \alpha( F(\gamma_{s,t+\tau})^{-1} , F(\Sigma_{s,0}(\sigma,t+\tau))) =  -  (\Sigma^{*}B)_{s,t} \left( \frac{\partial}{\partial s}, \frac{\partial}{\partial t} \right) \text{.}
\end{equation*}
Now, the three last equations show that $F(\Sigma_{0,0}(s,1))$ solves the above initial value problem.

So far we have proved equations (\ref{58}) on the level of objects. On the level of 1-morphisms, it is a consequence of Theorem \ref{th3}: for a pseudonatural transformation $\rho:F \to F'$ with components $g:X \to G$ and $\rho_H: P^1X \to H$ we have
\begin{equation*}
\fu(\fo(\rho)) \stackrel{\mathrm{(\ref{61})}}{=} \fu(g,\fo(F',\rho_H^{-1})) \stackrel{\mathrm{(\ref{62})}}{=} (g,\fu(\fo(F',\rho_H^{-1}))^{-1})\stackrel{\mathrm{Th. \ref{th3}}}{=}(g,\rho_H)=\rho\text{,}
\end{equation*}
and conversely, for a 1-morphism $(g,\varphi):(A,B) \to (A',B')$ in $\diffco{\mathfrak{G}}{2}{X}$,
\begin{equation*}
\fo(\fu(g,\varphi)) \stackrel{\mathrm{(\ref{62})}}{=} \fo(g,\fu(A',\varphi)^{-1}) \stackrel{\mathrm{(\ref{61})}}{=} \left(g,\fo \left( \left(\fu(A',\varphi)^{-1} \right)^{-1} \right)\right) \stackrel{\mathrm{Th. \ref{th3}}}{=}(g,\varphi)\text{.}
\end{equation*}
Finally, on the level of 2-morphisms, which are on both sides just the same $H$-valued functions on $X$, there is nothing to show.
\endofproof

\numberwithin{equation}{section}
\renewcommand\theequation{\thesection.\arabic{equation}}

\section{Examples of Smooth 2-Functors}

\label{sec4}

We give three  examples of situations where smooth 2-functors are present.

\subsection{Connections on (non-abelian) Gerbes}

\label{sec4_1}

Let us first recall from \cite{schreiber3} what connections on ordinary principal bundles have to do with ordinary functors. For $G$  a Lie group, we denote by $G\text{-}\mathrm{Tor}$ the category whose objects are smooth manifolds with transitive, free and smooth $G$-action from the right, and whose morphisms are $G$-equivariant smooth maps. The functor which regards $G$ itself as a $G$-space is denoted by $i_G:\mathcal{B}G \to G\text{-}\mathrm{Tor}$. If $\gamma:x \to y$ is a path in $X$, any principal $G$-bundle $P$ provides us with objects $P_x$ and $P_y$ of $G\text{-}\mathrm{Tor}$, namely its fibres over the endpoints of $\gamma$. Furthermore, a connection $\nabla$ on $P$ defines a morphism
\begin{equation*}
\tau_{\gamma}: P_x \to P_y
\end{equation*}
in $G\text{-}\mathrm{Tor}$, namely the parallel transport along $\gamma$.  Well-known properties of  parallel transport  assure that the assignments $x \mapsto P_x$ and $\gamma\mapsto \tau_{\gamma}$ define a functor
\begin{equation*}
\mathrm{tra}_{P,\nabla}: \mathcal{P}_1(X) \to G\text{-}\mathrm{Tor}\text{.}
\end{equation*}
The main result of \cite{schreiber3} is the characterization of functors obtained like this among all functors  $F:\mathcal{P}_1(X) \to G\text{-}\mathrm{Tor}$. They are characterized by the following defining property of a \emph{transport functor}: there exists a surjective submersion $\pi:Y \to M$ and a smooth functor $\mathrm{triv}:\mathcal{P}_1(Y) \to \mathcal{B}G$ such that the functors $\pi^{*}F$ and $i_G \circ \mathrm{triv}$ are (with additional conditions we skip here) naturally equivalent. In other words, transport functors are \emph{locally} smooth functors. These transport functors form a category $\mathrm{Trans}^1_{\mathcal{B}G}(X,G\text{-}\mathrm{Tor})$, and we have

\begin{theorem}[\cite{schreiber3}, Theorem 5.4]
\label{th5}
The assignment of a functor $\mathrm{tra}_P$ to a principal $G$-bundle $P$ with connection over $X$ defines a surjective equivalence of categories
\begin{equation*}
\mathfrak{Bun}_G^{\nabla}(X) \cong \mathrm{Trans}^1_{\mathcal{B}G}(X,G\text{-}\mathrm{Tor})\text{.}
\end{equation*}
\end{theorem} 

Under this equivalence, trivial principal $G$-bundles with connection correspond to \emph{globally} smooth functors, i.e. functors $\mathrm{tra}:\mathcal{P}_1(X) \to G\text{-}\mathrm{Tor}$ with $\mathrm{tra}=i_G \circ \mathrm{triv}$ for a smooth functor $\mathrm{triv}:\mathcal{P}_1(X) \to \mathcal{B}G$. Trivial\emph{izable} principal $G$-bundles with connection correspond to functors which are naturally equivalent to globally smooth functors (again with additional assumptions on the natural equivalence). 

\medskip

We think that the concept of transport functors is adequate to be categorified and to capture all aspects of connections on 2-bundles, in particular gerbes. 
We anticipate the following results of  \cite{schreiber2}:
\begin{enumerate}
\item
Gerbes with connection over $X$ have structure 2-groups $\mathfrak{G}$.
\item 
A trivial $\mathfrak{G}$-gerbe with connection over $X$ is a smooth 2-functor \begin{equation*}
F: \mathcal{P}_2(X) \to \mathcal{B}\mathfrak{G}\text{.}
\end{equation*}
\end{enumerate}
Let us test these assertions in two examples.

\begin{example}
\label{ex5}
\normalfont
We consider the Lie 2-group
 $\mathfrak{G}=\mathcal{B}U(1)$ from Example \ref{ex1}. The corresponding $\mathcal{B}U(1)$-gerbes are also known as abelian gerbes, or $U(1)$-gerbes.  Now, a trivial $\mathcal{B}U(1)$-gerbe with connection over $X$ is by the above assertion and Theorem \ref{th2} nothing but a 2-form $B \in \Omega^2(X)$. 
\end{example}

Abelian  gerbes with connection can be realized conveniently by \emph{bundle gerbes} \cite{murray}. In this context it is well-known that a connection on a trivial bundle gerbe is indeed just a 2-form, see, e.g., \cite{waldorf1}.

\begin{example}
\normalfont
\label{ex3}
Let $H$ be a connected Lie group. We denote by $\mathfrak{aut}(H)$  the Lie algebra of the Lie group $\mathrm{Aut}(H)$ of Lie group automorphisms of $H$. We consider the Lie 2-group $\mathfrak{G}=\mathrm{AUT}(H)$ from Example \ref{ex4}.
By the above assertion and Theorem \ref{th2}, a trivial  $\mathrm{AUT}(G)$-gerbe with connection over $X$ is a pair $(A,B)$ of a 1-form $A\in\Omega^1(X,\mathfrak{aut}(H))$ and a 2-form $B\in \Omega^2(X,\mathfrak{h})$ such that
\begin{equation}
\label{69}
\mathrm{d}A + [A \wedge A] = \mathrm{ad} ( B)\text{,}
\end{equation}
where $\mathrm{ad}: \mathfrak{h} \to \mathfrak{aut}(H): X \mapsto \mathrm{ad}_X$.
\end{example}

$\mathrm{AUT}(H)$-gerbes are also known as $H$-gerbes\footnote{We have to remark that a $U(1)$-gerbe in the sense of Example \ref{ex5} is not the same as an $H$-gerbe for $H=U(1)$ in the sense of Breen and Messing. The difference becomes clear if one uses the classification of gerbes by Lie \emph{2}-groups we have proposed here: we have $\mathcal{B}U(1)$-gerbes on one side but $\mathrm{AUT}(U(1))$-gerbes on the other. Indeed, $\mathcal{B}U(1)$ is only a sub-2-group of $\mathrm{AUT}(U(1))$. } in the sense of  Breen and Messing \cite{breen1}. There, a connection on a trivial $H$-gerbe is a pair $(A,B)$ just as in Example \ref{ex3} but without the condition (\ref{69}). This difference lies at the heart of a question N. Hitchin posed at the VBAC-meeting in Bad Honnef in June 2007 after a talk by L. Breen, namely if it is possible to define a surface holonomy from a connection on an $H$-gerbe. Let us presume that \tql a surface holonomy\tqr\ is at least a 2-functorial assignment, i.e. is described by a 2-functor on the path 2-groupoid. This assumption is supported by the approach via \tql 2-holonomies\tqr\  \cite{martins1}, as well as by  \emph{transport 2-functors}  \cite{schreiber2}.  Similar considerations have also been made for ordinary holonomy \cite{caetano,schreiber3}. Then, the following three statements on a Breen-Messing connection $(A,B)$ on a trivial $H$-gerbe over $X$  are equivalent:
\begin{enumerate}
\item[(a)] 
it defines a  surface holonomy.
\item[(b)]
it satisfies condition (\ref{69}), $\mathrm{d}A + [A \wedge A] = \mathrm{ad} (B)$.
\item[(c)]
there exists a smooth 2-functor 
\begin{equation*}
F:\mathcal{P}_2(X) \to \mathcal{B}\mathrm{AUT}(H)
\end{equation*}
such that $(A,B)=\fo(F)$.
\end{enumerate}
A detailed discussion of surface holonomies that also covers non-trivial $H$-gerbes, is postponed to \cite{schreiber2}.

\subsection{\Der s of Smooth Functors}

In Appendix \ref{app2} we have reviewed the Lie 2-group $\mathcal{E}G$ associated to any Lie group $G$.
For any functor $F: \mathcal{P}_1(X) \to \mathcal{B} G$ there is an associated 2-functor
\begin{equation*}
\mathrm{d}F: \mathcal{P}_2(X) \to \mathcal{BE}G
\end{equation*} 
that we call the \emph{\der\ 2-functor} of $F$. It sends a 1-morphism $\gamma \in P^1X$ to the image $F(\gamma)\in G$ of  $\gamma$ under the functor $F$. This determines $\mathrm{d}F$ completely, since the Lie 2-groupoid $\mathcal{BE}G$ has only one object, and precisely only one 2-morphism between any two fixed 1-morphisms. 
It will be interesting to determine the unique 2-morphism  $\mathrm{d}F(\Sigma)$ associated to a bigon $\Sigma: \gamma_1 \Rightarrow \gamma_2$ explicitly. For this purpose, we denote by $\partial\Sigma$ the 1-morphism $\gamma_1^{-1} \circ \gamma_2$. Then, we obtain directly from the definitions:

\begin{theorem}[The non-abelian Stokes' Theorem for functors]
\label{th1}
Let $G$ be a Lie group and let  $F: \mathcal{P}_1(X) \to \mathcal{B} G$ be a functor. Then,
\begin{equation*}
\mathrm{d}F(\Sigma)= F(\partial \Sigma)
\end{equation*}
for any bigon $\Sigma \in B^2X$. 
\end{theorem}

In order to understand why we call this identity Stokes' Theorem, notice that if the functor $F$ is smooth, also its \der\ 2-functor $\mathrm{d}F$ is smooth. Then, we have associated differential forms:

\begin{lemma}
\label{lem12}
 Let $A\in \Omega^1(X,\mathfrak{g})$ be the 1-form associated to the smooth functor $F$, and let  $B\in \Omega^2(X,\mathfrak{g})$ be the 2-form associated to its \der\ 2-functor $\mathrm{d}F$ by Theorem \ref{th2}. Then,
\begin{equation*}
B=[A \wedge A] + \mathrm{d}A \text{.}
\end{equation*}
\end{lemma}

\proof
We recall that there is also the 1-form $A'\in\Omega^1(X,\mathfrak{g})$ associated to the 2-functor $\mathrm{d}F$, and that by Proposition \ref{prop1} $B=[A' \wedge A'] + \mathrm{d}A'$, since $t$ is the identity in the crossed module that defines $\mathcal{E}G$. Furthermore, since $\mathrm{d}F(\gamma) = F(\gamma)$, we have $A=A'$. 
\endofproof

We have reviewed in Section \ref{sec4_1} that a smooth functor $F:\mathcal{P}_1(X) \to \mathcal{B}G$ corresponds to trivial
principal $G$-bundle $P$ with connection $\omega$ in the sense that $\mathrm{tra}_{P,\omega}=i_G \circ F$. By Lemma \ref{lem12}, the 2-form $B$ determined by the 2-functor $\mathrm{d}F$ is the curvature of this connection $\omega$. Moreover, the holonomy of $\omega$ around any closed path $\gamma$ (identified with a group element) is given by $F(\gamma)$. If $\gamma$ is of the form $\gamma=\partial\Sigma$ for any bigon $\Sigma$, Theorem \ref{th1} implies
\begin{equation*}
\mathrm{Hol}_{\nabla}(\partial\Sigma)=F(\partial\Sigma)= \mathrm{d}F(\Sigma)\text{;}
\end{equation*}
this is a relation between the \emph{holonomy} and the \emph{curvature} of a connection on a (trivial) principal $G$-bundle. Further restricted to the case that the bigon $\Sigma$ is of the form $\Sigma:\id_x \Rightarrow \gamma$ for a closed path $\gamma:x \to x$, we have

\begin{corollary}
\label{co1}
Let $\omega$ be a connection on a trivial principal $G$-bundle of curvature $K$, and let $\gamma$ be a contractible loop at $x\in X$. Then,
\begin{equation*}
\mathrm{Hol}_\omega(\gamma)=\mathcal{P}\mathrm{exp} \int_0^1 \mathcal{A}_{\Sigma} =\mathcal{P}\exp \int_0^1 \mathrm{d}s \left ( \int_0^1 \mathrm{d}t \; \mathrm{Ad}_{\tau(\gamma_{s,t})}^{-1}K|_{\Sigma(s,t)}\left ( \frac{\partial\Sigma}{\partial s},\frac{\partial\Sigma}{\partial t} \right ) \right )
\end{equation*}
where  $\Sigma:\id_x \Rightarrow \gamma$ is any choice of a smooth contraction of $\gamma$ to its base point,  the group element   $\tau(\gamma_{s,t})\in G$ is the parallel transport of the connection $\omega$ along the path $\gamma_{s,t}$, and the path-ordered exponential $\mathcal{P}\mathrm{exp}$ indicates the unique solution of the respective initial value problem, like in (\ref{74}).  
\end{corollary}
Exactly the same formula can been found in \cite{alvarez2}, derived by  completely different methods. In the abelian case of $G=U(1)$ Corollary \ref{co1} boils down to the well-known identity
\begin{equation*}
\mathrm{Hol}_{\omega}(\gamma) = \mathrm{exp} \left (  \mathrm{i}\int_{\Sigma} K \right )
\end{equation*}
for a surface $\Sigma$ with boundary $\gamma=\partial\Sigma$.

\subsection{Classical Solutions in BF-Theory}

Four-dimensional BF theory is a topological field theory  on a four-dimensional, compact, oriented smooth manifold $X$, see e.g. \cite{baez4}. Usually, it involves a symmetric, non-degenerate invariant bilinear form $\left \langle -,-  \right \rangle$ on the Lie algebra $\mathfrak{g}$ of a Lie group $G$, and the fields are pairs $(A,B)$ of a 1-form $A\in \Omega^1(X,\mathfrak{g})$ and a 2-form $B \in \Omega^2(X,\mathfrak{g})$. 

We infer  that a naturally generalized setup in which BF theory should be considered, is a Lie 2-group $\mathfrak{G}$, i.e. a smooth crossed module $(G,H,t,\alpha)$, together with the invariant form $\left \langle  -,- \right \rangle$ on the Lie algebra of $G$. Other generalizations have been proposed in \cite{pfeiffer1}. The fields are now pairs $(A,B)$ of a 1-form $A\in \Omega^1(X,\mathfrak{g})$ and a 2-form $B\in\Omega^2(X,\mathfrak{h})$, and the action is, with $F_A := \mathrm{d}A + [A\wedge A]$ and  $\beta_{A,B} :=F_{A} - t_{*}B$,
\begin{equation}
\label{72}
S(A,B) := \frac{1}{2}\int_X \left \langle \beta_{A,B} \wedge \beta_{A,B}  \right \rangle\text{.} \end{equation}
Expressed in terms of $A$ and $B$, this is
\begin{equation}
S(A,B) =\frac{1}{2} \int_X \left \langle F_A \wedge F_A  \right \rangle - \int_X \left \langle t_{*}B \wedge    F_A  \right \rangle + \frac{1}{2}\int_X \left \langle t_{*}B \wedge t_{*}B  \right \rangle\text{;}
\end{equation}
these terms can be identified as: a topological Yang-Mills term, the \tql original\tqr\ BF-term and a so-called cosmological term. The variation of this action gives 
\begin{eqnarray*}
\frac{\delta S}{\delta A} = 0 &\Leftrightarrow& t_{*}\mathrm{d}B + A \wedge t_{*}B = 0
\\
\frac{\delta S}{\delta B} =0 &\Leftrightarrow& \beta_{A,B}=0\text{.}
\end{eqnarray*}
We notice that the second equation implies  the first, so that the critical point of $S(A,B)$ are exactly those with $\beta_{A,B}=0$. It follows further that the topological Yang-Mills term, which is usually not present in BF theory, has no influence on the critical points. Since pairs $(A,B)$ with $\beta_{A,B}=0$ correspond by Theorem \ref{th2} to smooth 2-functors $F: \mathcal{P}_2(X) \to \mathcal{B}\mathfrak{G}$, we have

\begin{proposition}
\label{prop8}
The critical points of the BF action (\ref{72}) are  exactly the smooth 2-functors $F:\mathcal{P}_2(X) \to \mathcal{B}\mathfrak{G}$.
\end{proposition}

\section{Transgression to Loop Spaces}

\label{sec5}

In this section we use the  observation that  2-functors defined on  the path 2-groupoid $\mathcal{P}_2(X)$ of a smooth manifold $X$ induce structure on the loop space $LX:=C^{\infty}(S^1,X)$, since loops are particular 1-morphisms in  $\mathcal{P}_2(X)$. In order to understand this structure properly, we equip $LX$ with the canonical diffeology of $LX=D^{\infty}(S^1,X)$, see Section \ref{sec1_2}. In Section \ref{sec3_1a} we generalize the relation between smooth functors and differential forms (Theorem \ref{th3}) from smooth manifolds $X$ to arbitrary diffeological spaces, in particular to $LX$. In Section \ref{sec4_2} we combine this generalized statement on $LX$ with Theorem \ref{th2} on  $X$.

\subsection{Generalization to Diffeological Spaces}

\label{sec3_1a}

In order to describe a generalization of Theorem \ref{th3} from smooth manifolds to diffeological spaces, we first  have to define the path groupoid $\mathcal{P}_1(X)$ of a diffeological space $X$. We will see that almost all definitions we gave for $X$ a smooth manifold pass through; only the notion of  thin homotopy has to be adapted. 

So, a \textit{path} in $X$ is a diffeological  map $\gamma:[0,1] \to X$ with sitting instants. As described in Section \ref{sec1_2} the set $PX$ of paths  can be considered as a subset of the diffeological space $D^{\infty}((0,1),X)$, and  is hence itself a diffeological space. By axiom (D2) for diffeological spaces, the constant path $\id_x$ at a point $x\in X$ is diffeological.  Let us further show exemplarily

\begin{lemma}
\label{lem10}
The composition $\gamma_2 \circ \gamma_1$ of two paths $\gamma_1:x \to y$ and $\gamma_2:y \to z$ is again a path.  
\end{lemma}

\proof
Notice that if $\gamma:[0,1] \to X$ is a path and $U \subset [0,1]$ is open, then $\gamma|_U:U \to X$ is a plot of $X$. To see that the composition $\gamma_{2} \circ \gamma_1$ (which is defined in the same way as for smooth manifolds) is diffeological, let $U\subset[0,1]$ be  open, let $\epsilon_i$ be a sitting instant of $\gamma_i$, and let \begin{equation*}
U_1 := U \cap (0,\textstyle\frac{1}{2})
\quad\text{, }\quad
U_2:= U \cap (\frac{1}{2}-\epsilon_1,\frac{1}{2}+\epsilon_2)
\quad\text{ and }\quad
U_3:= U \cap (\frac{1}{2},1)\text{.}
\end{equation*}
These are open sets that cover $U$, furthermore, $(\gamma_2 \circ \gamma_1)|_{U_i} = \gamma_i|_{U_i}$ for $i=1,3$ are plots of $X$ and $(\gamma_2 \circ \gamma_1)|_{U_2}$ is constant and hence also a plot of $X$ by axiom (D2). Hence, $(\gamma_2 \circ \gamma_1)|_U$ is a plot of $X$ by axiom (D3).
\endofproof

We leave it to the reader to prove  that the inverse $\gamma^{-1}$ of a path $\gamma$ is again a path. Next we have to define thin homotopy for paths in a diffeological space. For this purpose, we first give a reformulation of a thin homotopy  on smooth manifolds, which generalizes better to diffeological spaces. 

\begin{lemma}
\label{lem9}
Let $X$ and $Y$ be smooth manifolds and $f:X \to Y$  be a smooth map. The rank of the differential of $f$ is bounded above by a number $k\in\N$ if and only if the pullback of every $(k+1)$-form $\omega\in \Omega^{k+1}(Y)$ along $f$ vanishes.
\end{lemma}

\proof
Assume that the rank of the differential of $f$ is at most $k$ everywhere. Then, $f^{*}\omega=0$ for all $\omega\in \Omega^{k+1}(Y)$. Conversely, assume that $f^{*}\omega=0$ for all $\omega\in \Omega^{k+1}(Y)$. Assume further that there exists a point $p\in X$ such that $\mathrm{d}f|_p$ has rank $k'>k$. Then, there exist vectors $v_1,...,v_{k'}\in T_pX$ such that their images $w_i := \mathrm{d}f|_p(v_i)$ are linear independent. Using a chart of a neighbourhood of $f(p)$ one can construct a $k'$-form $\omega\in\Omega^{k'}(Y)$ such that $\omega_{f(p)}(w_1,...,w_{k'})$ is non-zero. Since this is equal to $(f^{*}\omega)_p(v_1,...,v_{k'})$, we have a contradiction to the assumption that $f^{*}\omega=0$. 
\endofproof

We thus have reformulated restrictions on the rank of the differential of a smooth function in terms of pullbacks of differential forms. Now we generalize  to diffeological spaces.  

\begin{definition}
\label{def5}
Let $X$ be a diffeological space. A \emph{differential $k$-form} on $X$ is a family of $k$-forms $\omega_c \in \Omega^k(U)$ for every plot $c:U \to X$, such that
\begin{equation*}
\omega_{c_1}= f^{*}\omega_{c_2}
\end{equation*}
for every smooth map $f: U_1 \to U_2$ with $c_2 \circ f = c_1$.
\end{definition}

Notice that the $k$-forms on a diffeological space $X$ form a vector space $\Omega^k(X)$, and that the wedge product and the exterior derivative generalize naturally to differential forms on diffeological spaces. Furthermore, it is clear that a differential form $\omega$ on a smooth manifold $X$ induces a differential form on $X$ regarded as a diffeological space: for a chart $\phi:U \to X$ of $X$ one takes $\omega_\phi:=\phi^{*}\omega$. 
 We have also a very simple definition of pullbacks of differential forms on diffeological spaces along diffeological maps  $f:X \to Y$ between diffeological spaces $X$ and $Y$: the pullback $f^{*}\omega$ of a $k$-form $\omega=\lbrace \omega_c \rbrace$ on $Y$ is the $k$-form on $X$ defined by
\begin{equation*}
(f^{*}\omega)_c := \omega_{f \circ c}
\end{equation*}
for every plot $c$ of $X$. Here it is important that $f \circ c$, since $f$ was supposed to be diffeological, is a plot of $Y$. 
In particular, if $Y$ is a smooth manifold, $f \circ c : U \to Y$ is a smooth map and $(f^{*}\omega)_c = (f \circ c)^{*}\omega$.

\begin{definition}
\label{def4}
Two paths $\gamma_0:x \to y$ and $\gamma_1:x \to y$ in a diffeological space $X$ are called \emph{thin homotopy equivalent}, if there exists a diffeological map $h:[0,1]^2\to X$ with sitting instants as described in  (1) of Definition \ref{def3}, such that the pullback $h^{*}\omega$ of every 2-form $\omega\in\Omega^2(X)$ vanishes.
\end{definition}

By Lemma \ref{lem9} it is clear that for $X$ a smooth manifold Definition \ref{def4} is equivalent to Definition \ref{def3}. By arguments similar to those given in the proof of Lemma \ref{lem10} one can show that Definition \ref{def4} defines an equivalence relation $\sim_1$ on the diffeological space $PX$ of paths in $X$, so that the set of equivalence classes $P^1X := PX / \sim_1$ is again a diffeological space. This will be the set of morphisms of the path groupoid $\mathcal{P}_1(X)$ we are going to define. In the following lemma we prove that the axioms of a groupoid are satisfied. 

\begin{lemma}
Let $X$ be a diffeological space. For a path $\gamma:x \to y$ we have
\begin{equation*}
\gamma^{-1} \circ \gamma\; \sim_1 \;\id_x
\quad\text{ and }\quad
\id_y \circ \gamma \;\sim_1 \;\gamma\;\sim_1\; \gamma \circ \id_x\text{.}
\end{equation*}
For three paths $\gamma_1:x \to y$, $\gamma_2:y \to z$ and $\gamma_3: z \to w$ we have
\begin{equation*}
\gamma_1 \circ  (\gamma_2 \circ \gamma_3) \;\sim_1 \;(\gamma_1 \circ \gamma_2) \circ \gamma_3\text{.}
\end{equation*}
\end{lemma}

\proof
We prove $\gamma^{-1} \circ \gamma \sim_1 \id_x$; the remaining equivalences can be shown analogously. We choose the standard homotopy: this is, for some smooth map $\beta:[0,1] \to [0,1]$ with $\beta(0)=0$ and $\beta(1)=1$ and with sitting instants, the map
\begin{equation*}
h:[0,1]^2 \to X:(s,t) \mapsto \begin{cases}\gamma(2\beta(s)t) & 0 \leq t \leq\frac{1}{2} \\
\gamma(2\beta(s)(1-t)) & \frac{1}{2} < t\leq 1\text{.}
\end{cases}
\end{equation*}
This map has sitting instants. To see that it is diffeological, we use the same trick as in the proof of Lemma \ref{lem10}, i.e. we cover $(0,1)^2$ with
\begin{equation*}
V_{\gamma} := (0,1) \times (0,\textstyle\frac{1}{2})
\hspace{0.2cm}\text{, }\hspace{0.2cm}
V_{\gamma^{-1}} := (0,1) \times (\textstyle\frac{1}{2},1)
\hspace{0.2cm}\text{ and }\hspace{0.2cm}
V_{\epsilon}:=(0,1) \times (\frac{1}{2}-\epsilon,\frac{1}{2}+\epsilon)
\end{equation*} 
for $\epsilon$ a sitting instant $\gamma$, and accordingly any open subset $U\subset [0,1]^2$ by $V_{\gamma}\cap U$, $V_{\gamma^{-1}}\cap U$ and $V_{\epsilon}\cap U$. Now, 
\begin{equation}
\label{63}
h|_{V_{\gamma}\cap U} = (\gamma \circ m_{\beta})|_{V_{\gamma}\cap U}
\end{equation}
with $m_{\beta}(s,t):= 2\beta(s)t$; it is thus the composition of a plot with a smooth map and hence by axiom (D1) a plot. Similarly $h|_{V_{\gamma^{-1}}\cap U}$ and $h_{V_{\epsilon}\cap U}$ are plots. This shows that $h|_U$ is covered by plots and thus itself a plot. This implies that $h$ is diffeological. It remains to show that the pullback $h^{*}\omega$ of every $2$-form $\omega\in\Omega^2(X)$ vanishes. This follows from the fact that $h$ restricted to each of the subsets $V_{\gamma}$, $V_{\gamma^{-1}}$ and $V_{\epsilon}$ is either constant or factors as in (\ref{63}) through the one-dimensional manifold $[0,1]$ via $\gamma$ or $\gamma^{-1}$, respectively. 
\endofproof

This finishes the definitions of the path groupoid $\mathcal{P}_1(X)$ of a diffeological space $X$. It is clear that one now can consider smooth functors
\begin{equation*}
F:\mathcal{P}_1(X) \to S
\end{equation*}
into any Lie category $S$ like before: the maps $F_0: X \to S_0$ on objects and $F_1: P^1X \to S_1$ on morphisms have to be diffeological maps. 

\medskip

Further towards a generalization of Theorem \ref{th3} we have to generalize the category $\diffco{G}{1}{X}$ introduced in Definition \ref{def8} from a smooth manifold $X$ to a diffeological space. Notice that Definition \ref{def5}  extends naturally to $\mathfrak{g}$-valued differential forms on diffeological spaces. Now, for a diffeological space $X$ an object in $\diffco{G}{1}{X}$ is a $\mathfrak{g}$-valued 1-form $A=\lbrace A_c \rbrace$ on $X$. A morphism $g: A \to A'$ is a diffeological map $g:X \to G$ such that for any plot $c: U \to X$ and the associated smooth map $g_c := g \circ c:U \to X$
\begin{equation}
\label{64}
A'_{c} = \mathrm{Ad}_{g_c}(A_{c}) - g_c^{*}\bar\theta\text{.}
\end{equation}
The functor $\fo$ from Section \ref{sec1_3} generalizes straightforwardly to a functor
\begin{equation*}
\fo: \mathrm{Funct}^{\infty}(\mathcal{P}_1(X),\mathcal{B}G) \to \diffco{G}{1}{X}
\end{equation*}
for any diffeological space $X$: 
\begin{itemize}
\item 
Let $F:\mathcal{P}_1(X) \to \mathcal{B}G$ be a smooth functor. For any plot  $c:U \to X$ of $X$ (which is itself a diffeological map), the pullback $c^{*}F$ is a smooth functor $c^{*}F:\mathcal{P}_1(U) \to \mathcal{B}G$ defined on the path groupoid of the smooth manifold $U$. Hence $A_c := \fo(c^{*}F)\in\Omega^1(U,\mathfrak{g})$ is a 1-form. If $c':U'\to X$ is another plot and $f:U \to U'$ is a smooth map with $c=c' \circ f$, we have by Proposition \ref{prop6} $A_c =\fo(c^{*}F)=f^{*} \fo(c'^{*}F)=f^{*}A_{c'}$. 

\item
Let $\rho:F \to F'$ be a smooth natural transformation. Its components furnish a diffeological map $g:X \to G$. For any plot $c:U \to X$, we have
$g_c:=\rho \circ c = \mathcal{D}(c^{*}\rho):U \to X$, hence, since $\fo$ is a functor,  (\ref{64}) is satisfied. 
\end{itemize}
The extension of the inverse functor $\fu$ to diffeological spaces is slightly more involved. Let $A=\lbrace A_c \rbrace \in \Omega^1(X,\mathfrak{g})$ be a 1-form on the diffeological space $X$. For every plot $c:U \to X$ we obtain a smooth functor $F_c := \fu(A_c)$. In particular, since every path $\gamma:[0,1] \to X$ defines a plot $\gamma|_{(0,1)}$, we have functors $F_{\gamma}$ defined on the path groupoid of the open interval  $(0,1)$. Let $\epsilon_{s,t}\in P^1((0,1))$ be the path in $(0,1)$ that goes from $s+\epsilon$ to $t-\epsilon$, where $\epsilon$ is a sitting instant of $\gamma$. 
Then, we define a map
\begin{equation*}
F: P^1X \to G: \gamma \mapsto F_{\gamma}(\epsilon_{0,1})\text{.}
\end{equation*}

\begin{lemma}
This defines a smooth functor $F:\mathcal{P}_1(X) \to \mathcal{B}G$.
\end{lemma}

\proof
To see that $F: P^1X \to G$ is diffeological, we have to show that for every plot $c:U \to P^1X$ the composite $F \circ c:U \to G$ is a smooth map. Since we can check smoothness locally, we may assume that $c=\mathrm{pr} \circ c'$ for a plot $c': U \to PX$ and the projection $\mathrm{pr}:PX \to P^1X$. The relevant evaluation map  $\tilde c: U \times (0,1) \to X$ given by
\begin{equation*}
\alxydim{}{U \times (0,1) \ar[r]^-{c' \times \id} & PX \times (0,1) \ar[r]^-{\mathrm{ev}} & X}
\end{equation*}
is  a plot of $X$. Hence, we have a smooth functor $F_{\tilde c}:\mathcal{P}_1(U \times (0,1)) \to \mathcal{B}G$. With the map $i_u: (0,1) \to U \times (0,1): t \mapsto (u,t)$ we have a plot $\tilde c \circ i_u$ and accordingly $A_{c'(u)}= i_u^{*}A_{\tilde c}$ for all $u\in U$. Then, by Proposition \ref{prop6},
\begin{equation*}
F(c(u)) = i_u^{*}F_{\tilde c}(\epsilon_{0,1})= F_{\tilde c}((i_u)_{*}\epsilon_{0,1})\text{.}
\end{equation*}
Since $U \to P^1(U \times (0,1)): u \mapsto (i_u)_{*}\epsilon_{1,2}$ is a diffeological map, and $F_{\tilde c}$ is diffeological, we have shown that $F \circ c$ is smooth. The compatibility of $F$ with the composition of paths follows from
\begin{equation}
\label{eq1}
F_{\gamma' \circ \gamma}(\epsilon_{0,1})= F_{\gamma'\circ\gamma}(\epsilon_{1\!/\!2,1}) \cdot F_{\gamma' \circ \gamma}(\epsilon_{0,1\!/\!2}) = F_{\gamma'}(\epsilon_{0,1}) \cdot F_{\gamma}(\epsilon_{0,1}) \text{.}
\end{equation} 
For the last step of \erf{eq1}, we show the equality $F_{\gamma}(\epsilon_{0,1}) = F_{\gamma' \circ \gamma}(\epsilon_{0,1\!/\!2})$; the  one for $\gamma'$ goes analogously.  Indeed, consider the inclusion $\iota_1: [0,1] \to [0,1]$ defined by $\iota_1(t) := \frac{1}{2}t$. Pulling back $A_{\gamma' \circ \gamma}$ along $\iota_1$ and using Proposition \ref{prop6} we get $F_{\gamma} = \iota_1^{*}F_{\gamma' \circ \gamma}$. Evaluating this on the path $\varepsilon_{0,1}$, and using that $(\iota_1)_{*}(\varepsilon_{0,1}) = \varepsilon_{0,1/2}$ in $P^1((0,1))$ shows the claim. 
\endofproof

Now the following theorem  follows from Theorem \ref{th3} applied to functors and forms on the codomain $U$ of each plot $c:U \to X$ of $X$.

\begin{theorem}
\label{th4}
Let $X$ be a diffeological space and $G$ a Lie group.
The functors
\begin{equation*}
\fo: \mathrm{Funct}^{\infty}(\mathcal{P}_1(X),\mathcal{B} G) \to \diffco{G}{1}{X}
\end{equation*}
and
\begin{equation*}
\fu: \diffco{G}{1}{X} \to \mathrm{Funct}^{\infty}(\mathcal{P}_1(X),\mathcal{B} G)
\end{equation*}
satisfy
\begin{equation*}
\fo \circ \fu = \id_{\diffco{G}{1}{X}}
\quad\text{ and }\quad
\fu \circ \fo = \id_{\mathrm{Funct}(\mathcal{P}_1(X),\mathcal{B}G)}\text{,}
\end{equation*}
and are hence isomorphisms  of categories.
\end{theorem}

\subsection{Induced Structure on the Loop Space}

\label{sec4_2}

In this section we discuss the diffeological space $LX=D^\infty(S^1,X)$, where $X$ is a smooth manifold. In order to formalize the relation between functors defined on the path groupoid $\mathcal{P}_1(X)$ and structure on $LX$ we are going to explore we introduce two constructions.

Firstly, we denote for any category $T$ by $\Lambda T := T_1$ set of morphisms in $T$.\label{not:lambda} Accordingly, for a functor $F:S \to T$, we call its induced map on morphisms $\Lambda F: \Lambda S \to \Lambda T$. Clearly, if $F$ was a diffeological functor, $\Lambda F$ is a diffeological map.
Secondly, we introduce a diffeological map 
\begin{equation}
\label{65}
\ell : LX \to \Lambda\mathcal{P}_1(X) \text{.} 
\end{equation}
Its definition is not completely obvious since loops have no sitting instants. We fix some smooth map $\beta:[0,1] \to [0,1]$ with $\beta(0)=0$ and $\beta(1)=1$ and with sitting instants. We have a smooth map  
$e_{\beta} : [0,1] \to S^1$ defined by $e_{\beta}(t) := \mathrm{e}^{2\pi \mathrm{i} \beta(t)}$ and accordingly a diffeological map
\begin{equation*}
\ell_{\beta}:LX \to PX: \tau \mapsto \tau \circ e_{\beta}\text{.}
\end{equation*}
We define $\ell := \mathrm{pr} \circ \ell_{\beta}$, where $\mathrm{pr}:PX \to P^1X$ is the projection to thin homotopy classes. This map $\ell$ is diffeological and indeed independent of the choice of $\beta$: for another choice $\beta'$ and some $\tau\in LX$ we find a thin homotopy $\ell_{\beta}(\tau) \sim_1 \ell_{\beta'}(\tau)$ for example by 
\begin{equation}
\label{75}
h:[0,1]^2 \to X: (s,t) \mapsto \tau \left(\mathrm{e}^{2\pi\mathrm{i}(\beta(s) \beta'(t) + (1-\beta(s))\beta(t))} \right)\text{;}
\end{equation}
this map is diffeological, has sitting instants and is evidently thin, since it factors through $S^1$. 

Now, having the two definitions $\Lambda$ and $\ell$ at hand, for $F: \mathcal{P}_1(X) \to T$ a smooth functor, 
\begin{equation}
\label{68}
\Lambda F \circ \ell : LX \to \Lambda T
\end{equation}
is a diffeological map on the loop space. A particular situation arises if the category $T=\mathcal{B}G$ for a Lie group $G$. In this case $\Lambda \mathcal{B}G=G$. We have now obtained a map
\begin{equation}
\label{67}
\mathrm{H}_1: \alxydim{}{\left \lbrace \txt{ Smooth functors\\$F:\mathcal{P}_1(X) \rightarrow \mathcal{B}G$} \right \rbrace} \to D^{\infty}(LX,G)\text{.}
\end{equation}
This map is of course well-known: as mentioned in Section \ref{sec4_1}, a smooth functor $F: \mathcal{P}_1(X) \to \mathcal{B}G$ corresponds to a (trivial) principal $G$-bundle $P$ with connection $\omega$ over $X$, in such a way that the parallel transport along a path $\gamma$ in $X$ is given by multiplication with $F(\gamma)$. For a loop $\tau\in LX$, understood as a path $\ell(\tau)$, this means
\begin{equation*}
\mathrm{H}_1(F)(\tau) = F(\ell(\tau)) = \mathrm{Hol}_\omega(\tau)\text{,}
\end{equation*}
so that $\mathrm{H}_1(F)$ is nothing but the holonomy of the connection $\omega$ around $\gamma$.
\medskip

In the following we  explore which structure on the loop space $LX$ is induced from a smooth \emph{2-}functor $F: \mathcal{P}_2(X) \to \mathcal{B}\mathfrak{G}$. To start with, we  generalize the two constructions $\Lambda$ and $\ell$ we have described before, to 2-categories. 

\begin{definition}
Let $T$ be a 2-category. We define a category $\Lambda T$ as follows: the objects are the 1-morphisms $T_{1}$ of $T$, and  the morphisms between two objects $f:X_f \to Y_f$ and $g:X_g \to Y_g$ are triples $(x,y,\varphi)\in T_1 \times T_1 \times T_2$ of 1-morphisms $x:X_f \to X_g$ and $y:Y_f \to Y_g$ and of a 2-morphism
\begin{equation*}
\alxydim{@=1.2cm}{X_f \ar[d]_{f} \ar[r]^{x} & X_g \ar[d]^{g} \ar@{=>}[dl]|{\varphi} \\ Y_f \ar[r]_{y} & Y_g}\text{.}
\end{equation*}
The composition in $\Lambda T$ is  putting these squares next to each other, and the identity of an object $f:X \to Y$ is the triple $(\id_X,\id_Y,\id_f)$.
\end{definition}

Clearly, if the sets $T_{1}$ and $T_2$ of the 2-category $T$ are diffeological spaces, the objects and morphisms of $\Lambda T$ form also diffeological spaces. For $F:S \to T$ a 2-functor, we have an associated functor 
\begin{equation*}
\Lambda F: \Lambda S \to \Lambda T\text{,}
\end{equation*}
which just acts as $F$ on 1-morphisms and 2-morphisms of $S$. If the 2-functor $F$ is diffeological, the functor $\Lambda F$ is also diffeological.
Next we  generalize the diffeological map $\ell$ introduced above to a diffeological functor \label{not:ell}
\begin{equation*}
\ell: \mathcal{P}_1(LX) \to \Lambda \mathcal{P}_2(X)\text{.}
\end{equation*}
On objects, it is just the map  $\ell$ from (\ref{65}),  regarding a loop $\tau\in LX$ as a particular path in $X$, i.e. as an object in $\Lambda \mathcal{P}_2(X)$. To define $\ell$ on morphisms, let $\gamma$ be a path in $LX$, i.e. a diffeological map $\gamma:[0,1] \to D^{\infty}(S^1,X)$ with sitting instants. We have an associated smooth map $m^{\gamma}:\R^2 \to X$ defined by $m^{\gamma}(s,t) := \gamma(t)(\mathrm{e}^{2\pi\mathrm{i}s})$, where we assume $\gamma$ to be trivially extended to $\R$ in the usual way.
Using the standard bigon $\Sigma_{\R}(s,t) \in B^2\R^2$ from (\ref{25}), we have a bigon
\begin{equation}
\label{78}
m(\gamma) := m^{\gamma}_{*}(\Sigma_{\R}(1,1)) \in B^2X
\end{equation}
associated to the path $\gamma$, and thus a well-defined map $m: PLX \to B^2X$.
\begin{lemma}
The map $m: PLX \to B^2X$ is diffeological.
\end{lemma}

\proof
We have to show that for any plot $c: U \to PLX$  the composite $m \circ c$ is a plot of $B^2X$. This means that, for a fixed representative $\Sigma \in BX$ of $\Sigma_{\R}(1,1)$, and any open subset $W \subset [0,1]^2$, the associated map
\begin{equation}
\label{76}
\alxydim{@C=1.5cm}{U \times W \ar[r]^-{c \times \id} & PLX \times W \ar[r]^{m_{*}\Sigma \times \id} &  BX \times W \ar[r]^-{\mathrm{ev}} & X}
\end{equation}
has to be smooth. Let us define the open intervals $V:=p_2(\Sigma(W))$ and $V':=p_1(\Sigma(W))$ for $p_i:[0,1]^2 \to [0,1]$ the canonical projections, and consider the chart $\varphi:V' \to S^1:s \mapsto \mathrm{e}^{2\pi\mathrm{i}s}$ of $S^1$.
Going through all involved definitions shows that the map (\ref{76}) coincides with the composite
\begin{equation}
\label{77}
\alxydim{@C=1.3cm}{U \times W \ar[r]^-{\id \times \Sigma} & U \times V \times V' \ar[r]^-{c'} & X}
\end{equation}
where $c'$ is given by
\begin{equation*}
\alxydim{@C=1.5cm}{U \times V \times V' \ar[r]^-{c \times \id \times \id} & D^{\infty}([0,1],X) \times V \times V' \ar[r]^-{\mathrm{ev} \times \varphi} & LX \times S^1 \ar[r]^-{\mathrm{ev}} & X}
\end{equation*}
Now, (\ref{77}) is the composition of two smooth maps, where $c'$ is smooth because $c$ was supposed to be a plot of $PLX \subset D^{\infty}([0,1],X)$.
\endofproof

We show next that the bigon $m(\gamma) \in B^2X$  does not depend on the thin homotopy class of the path $\gamma$. For this purpose, let $h:[0,1] \to LX$ be a thin homotopy between two paths $\gamma,\gamma'\in PLX$, and let $\Sigma \in BX$ be a representative for the bigon $\Sigma_{\R}(1,1)$.  We have an associated map $m^{h}: [0,1]^3 \to X$ defined by $m^{h}(r,s,t) := h(r,t)(\mathrm{e}^{2\pi\mathrm{i}s})$. Then, the map
\begin{equation*}
H: [0,1]^3 \to X : (r,s,t) \mapsto m^{h}(r,\Sigma(s,t))
\end{equation*} 
is a thin homotopy between  $m^{\gamma}_{*}\Sigma$ and $m^{\gamma'}_{*}(\Sigma)$. Hence, we have obtained a diffeological map $m: P^1LX \to B^2X$. 
Going through the definitions, one finds that this bigon $m(\gamma)$ has (up to thin homotopy) the following target and source paths:
\begin{equation}
\label{66}
\alxydim{@C=1.7cm@R=1.2cm}{\gamma(0)(1) \ar[d]_{\ell(\gamma(0))} \ar[r]^{\mathrm{pr}(b \circ \gamma)} & \gamma(1)(1) \ar[d]^{\ell(\gamma(1))} \ar@{=>}[dl]|{m(\gamma)} \\ \gamma(0)(1) \ar[r]_{\mathrm{pr}(b \circ \gamma)} & \gamma(1)(1)\text{,}}
\end{equation}
where $b:LX \to X$ is the base point evaluation. Hence, the triple
\begin{equation*}
\ell(\gamma) := (\mathrm{pr}(b \circ \gamma),\mathrm{pr}(b \circ \gamma),m(\gamma))
\end{equation*}
is a morphism in $\Lambda\mathcal{P}_2(X)$. The composition of paths in $LX$ is respected in the sense that $\ell(\gamma_2 \circ \gamma_1) = \ell(\gamma_2) \circ \ell(\gamma_1)$ where the latter is the composition in $\Lambda\mathcal{P}_2(X)$. Thus, we have completely  defined the diffeological functor
\begin{equation*}
\ell: \mathcal{P}_1(LX) \to \Lambda \mathcal{P}_2(X)\text{.}
\end{equation*}

\medskip

If now $F:\mathcal{P}_2(X) \to T$ is a smooth 2-functor, we obtain an associated smooth functor
\begin{equation*}
\Lambda F \circ \ell: \mathcal{P}_1(LX) \to \Lambda T\text{,}
\end{equation*}
generalizing the map (\ref{68}).
In the important case that $T=\mathcal{B}\mathfrak{G}$ for $\mathfrak{G}$ a Lie 2-group, the groupoid $\Lambda \mathcal{B}\mathfrak{G}$ is -- following a notion of \cite{mackenzie} -- trivializable: a groupoid is called \emph{trivializable}, if it is equivalent to a groupoid of the form  $\mathrm{Gr}_{S,N}=S \times \mathcal{B}N$ for a set $S$ regarded as a category with only identity morphisms, and a group $N$; these groupoids are called \emph{trivial}. Explicitly, the objects of $\mathrm{Gr}_{S,N}$ are the elements of $S$,  the Hom-set between two objects $s_1$ and $s_2$ is $N$ if $s_1=s_2$ and empty else,  and the  composition is  multiplication in $N$.
In our case we find  such an equivalence $\mathrm{tr}:\Lambda \mathcal{B}\mathfrak{G} \to \mathrm{Gr}_{G,G \ltimes H}$  as follows: on objects, $\mathrm{tr}$ is just the identity on $G$. A morphism \begin{equation*}
\alxydim{@=1.2cm}{\ast \ar[d]_{g_1} \ar[r]^{x} & \ast \ar[d]^{g_2} \ar@{=>}[dl]|{h} \\ \ast \ar[r]_{y} & \ast}\text{.}
\end{equation*}
in $\Lambda \mathcal{B}\mathfrak{G}$ is sent to the morphism $(y,h^{-1})$ in $\mathrm{Gr}_{G,G\ltimes H}$,  the inverse being necessary in order to respect the composition. The functor $\mathrm{tr}: \Lambda\mathcal{B}\mathfrak{G} \to \mathrm{Gr}_{G,G \ltimes H}$ is an equivalence of categories: the inverse functor sends a morphism $(y,h):g_1 \to g_2$ in $\mathrm{Gr}_{G,G\ltimes H}$ to the triple $(g_2^{-1}t(h)yg_1,y,h^{-1})$. 
Now, we have constructed a smooth functor
\begin{equation*}
\mathrm{tr} \circ \Lambda F \circ \ell : \mathcal{P}_1(LX) \to \mathrm{Gr}_{G,G\ltimes H}\text{.}
\end{equation*}
According to the direct product structure of the target Lie groupoid, this functor splits in
\begin{enumerate}
\item 
a diffeological function $h_{F}: LX \to G$ and
\item
a smooth functor $\mathcal{P}_1(LX) \to \mathcal{B}(G \ltimes H)$, which in turn corresponds by Theorem \ref{th4} and projection to the factors to
\begin{enumerate}
\item 
a 1-form $A_F\in \Omega^1(LX,\mathfrak{g})$ and
\item
a 1-form $\varphi_F \in \Omega^1(LX,\mathfrak{h})$.
\end{enumerate}
\end{enumerate}
Summarizing, we have, for any smooth manifold $X$ and any Lie 2-group $\mathfrak{G}$, a map
\begin{equation*}
\mathrm{H}_2:\alxydim{}{\left \lbrace \txt{ Smooth 2-functors\\$F:\mathcal{P}_2(X) \rightarrow \mathcal{B}\mathfrak{G}$} \right \rbrace}
\to D^{\infty}(LX,G) \times \Omega^1(LX,\mathfrak{g})\times \Omega^1(LX,\mathfrak{h})\text{.}
\end{equation*}
This map generalizes the map $\mathrm{H}_1$ from (\ref{67}) from smooth functors to smooth 2-functors. Let us describe the image $\mathrm{H_2}(F)=(h_F,A_F,\varphi_F)$. The diffeological function $h_{F}:LX\to G$ is clearly   $h_F = \mathrm{H}_{1}(F_{0,1})$ for the restriction $F_{0,1}$ of $F$ to objects and 1-morphisms. The differential forms $A_{F}$ and $\varphi_F$ can be characterized as described in the following proposition, also see Figure \ref{fig3}.

\begin{proposition}
\label{prop7}
Let $F:\mathcal{P}_2(X) \to \mathcal{B}\mathfrak{G}$ be a smooth 2-functor, let $A\in \Omega^1(X,\mathfrak{g})$ and $B\in \Omega^2(X,\mathfrak{h})$  the corresponding differential forms on $X$, and let $A_F$ and $B_F$ the differential forms on the loop space determined by $\mathrm{H}_2(F)$. Then
\begin{equation*}
A_F = b^{*}A
\quad\text{ and }\quad
B_F = \int_{S^1} (\alpha_{F \circ \gamma})_{*} (\mathrm{ev}^{*}B)\text{,}
\end{equation*}
where $b:LX \to X$ is the projection to the base point,  $\mathrm{ev}: LX \times S^1 \to X$ is the evaluation map, $\gamma: LX \times S^1 \to P^1X$ assigns to a loop $\tau\in LX$ and $z\in S^1$  the path obtained by parsing the loop $\tau$ from $z$ to $1$ counterclockwise, and $\alpha_{F \circ \gamma}: H \to H$ is the action of the crossed module $\mathfrak{G}$ along the map $F \circ \gamma: LX \times S^1 \to G$. 
\end{proposition}

\proof
Let $c: U \to LX$ be a plot of $LX$and let  $\Gamma:\R \to U$ be a smooth curve. Using all involved definitions we obtain
\begin{eqnarray*}
(A_F)_c|_{\Gamma(0)}\left ( \left .\frac{\partial \Gamma}{\partial t} \right|_0 \right) &=&- \left .\frac{\mathrm{d}}{\mathrm{d}t} \right|_0 (p_G \circ c^{*}(\mathrm{tr} \circ \Lambda F \circ \ell) \circ \Gamma_{*} \circ \gamma_{\R})(0,t)
\\
(b^{*}A)_c|_{\Gamma(0)}\left ( \left .\frac{\partial \Gamma}{\partial t} \right|_0 \right)&=&
-\left .\frac{\mathrm{d}}{\mathrm{d}t} \right|_0 ((b \circ c)^{*}F \circ \Gamma_{*} \circ \gamma_{\R})(0,t)\text{.}
\end{eqnarray*}
Then one observes that
\begin{equation*}
F \circ b_{*}=  p_G \circ \mathrm{tr} \circ \Lambda F \circ \ell
\end{equation*}
as maps from $P^1LX$ to $G$; this shows the first equality.

In order to proof the second equality we still use the plot $c$  and the smooth curve $\Gamma$, and consider the path $\gamma_{t} := c_{*}(\Gamma_{*}(\gamma_{\R}(0,t)))\in P^1LX$, where $\gamma_{\R}(s,t) \in P^1\R$ is the standard path from $s$ to $t$. We have
\begin{multline*}
(\varphi_F)_c|_{\Gamma(0)}\left ( \left . \frac{\partial\Gamma}{\partial t} \right|_0 \right ) = - \left . \frac{\mathrm{d}}{\mathrm{d}t} \right|_0 (p_H \circ c^{*}(\mathrm{tr} \circ \Lambda F \circ \ell) \circ  \Gamma_{*} \circ \gamma_{\R})(0,t)
\\= -\left . \frac{\mathrm{d}}{\mathrm{d}t} \right|_0 p_H(F(m(\gamma_t)))^{-1} = \left .  \frac{\mathrm{d}}{\mathrm{d}t} \right |_0 p_H(F(m^{\gamma_t}_{*}\Sigma_{\R}(1,1)))\text{,}
\label{79}
\end{multline*}
where we have used the definitions of the functor $\ell$  and the map $m$ from (\ref{78}). Let us remark that for the bigon $\Sigma^{\Gamma}(s,t) := m_{*}^{\gamma_1}\Sigma_{\R}(s,t)$ we have $\Sigma^{\Gamma}(1,t) = m_{*}^{\gamma_t}\Sigma_{\R}(1,1)$, so that we may write
\begin{equation}
\label{79}
 \left .  \frac{\mathrm{d}}{\mathrm{d}t} \right |_0 p_H(F(m^{\gamma_t}_{*}\Sigma_{\R}(1,1))) =  \int_0^1 \mathrm{d}\theta \;  \left.\frac{\partial^2}{\partial \rho \partial t} \right |_0 p_H(F(\Sigma^{\Gamma}_{1-\theta-\rho}(\theta+\rho,t)))\text{.}
\end{equation}
A by now standard calculation shows that
\begin{multline*}
p_H(F(\Sigma^{\Gamma}_{1-\theta-\rho}(\theta+\rho,t))) \\= p_H(F(\Sigma^{\Gamma}_{1-\theta}(\theta,t))) \cdot \alpha(F(\gamma(c(\Gamma(t)),\mathrm{e}^{-2\pi\mathrm{i}\theta})),F(\Sigma^{\Gamma}_{1-\theta-\rho}(\rho,t)))
\end{multline*}
where we have used the map $\gamma: LX \times S^1 \to P^1X$.
Now the derivative in (\ref{79}) becomes
\begin{equation*}
\left.\frac{\partial^2}{\partial \rho \partial t} \right |_0 p_H(F(\Sigma^{\Gamma}_{1-\theta-\rho}(\theta+\rho,t))) = 
-(\alpha_{F(\gamma(c(\Gamma(t)),\mathrm{e}^{-2\pi\mathrm{i}\theta})})_{*} \left ( \left.\frac{\partial^2}{\partial \rho \partial t} \right |_0 F(\Sigma^{\Gamma}_{1-\theta}(\rho,t))  \right )\text{.}
\end{equation*}
Let us now induce from the plot $c:U \to LX$ of $LX$ a plot of $LX \times S^1$, namely the map
\begin{equation*}
\tilde c: U \times (0,1) \to LX \times S^1: (u,\theta) \mapsto (c(u),\mathrm{e}^{2\pi\mathrm{i}\theta})\text{.}
\end{equation*}
We have $m^{\gamma_1}(1-\theta,0) = (\mathrm{ev} \circ \tilde c)(\Gamma(0),-\theta) \in X$
and the tangent vectors 
\begin{eqnarray*}
\left . \frac{\mathrm{d}}{\mathrm{d}t} \right |_0 m^{\gamma_1}(1-\theta,t) &=& \mathrm{d}(\mathrm{ev} \circ \tilde c)|_{(\Gamma(0),-\theta)}\left (  \left . \frac{\partial\Gamma}{\partial t} \right |_0 \right )
\\
\left . \frac{\mathrm{d}}{\mathrm{d}\rho} \right |_{1-\theta} m^{\gamma_1}(\rho,0) &=& \mathrm{d}(\mathrm{ev} \circ \tilde c)|_{(\Gamma(0),-\theta)} \left(\frac{\partial}{\partial\theta} \right)\text{.}
\end{eqnarray*}
 By Proposition \ref{prop4} we hence have
\begin{equation*}
\left.\frac{\partial^2}{\partial \rho \partial t} \right |_0 F(\Sigma^{\Gamma}_{1-\theta}(\rho,t)) =- (\mathrm{ev}^{*}B)_{\tilde c}|_{(\Gamma(0),-\theta)} \left ( \frac{\partial}{\partial\theta} , \left .\frac{\partial\Gamma}{\partial t} \right|_0 \right )\text{.}
\end{equation*}
Putting all pieces together and transforming $\theta\mapsto -\theta$, we have shown
\begin{equation*}
(\varphi_F)_c|_{\Gamma(0)}\left ( \left . \frac{\partial\Gamma}{\partial t} \right|_0 \right ) = \int_{-1}^{0} \mathrm{d}\theta \;  (\alpha_{(F \circ \gamma)(\tilde c(\Gamma(0),\theta))})_{*}(\mathrm{ev}^{*}B)_{\tilde c}|_{(\Gamma(0),\theta)} \left ( \frac{\partial}{\partial\theta} , \left .\frac{\partial\Gamma}{\partial t} \right|_0 \right )
\end{equation*}
this is the announced fibre integral written in the plot $\tilde c$ of $LX \times S^1$. 
\endofproof

\begin{figure}[h]
\begin{equation*}
\alxydim{@R=1.4cm@C=3cm}{F\ar@{|->}[d]_{\text{Theorem \ref{th2}}} \ar@{|->}[r] & \mathrm{tr} \circ \Lambda F \circ \ell  \ar@{|->}[d]^{\text{Theorem \ref{th4}}} \\ (A,B) \ar@{|->}[r]_{b^{*} \times \int_{S^1} (\alpha_{F \circ \gamma})_{*} \circ \mathrm{ev}^{*}} & (A_F,\varphi_F)}
\end{equation*}
\caption{A diagram for manipulations on a smooth 2-functor $F: \mathcal{P}_2(X) \rightarrow \mathcal{B}\mathfrak{G}$, whose commutativity is Proposition \ref{prop7}. The first row consists of functors, and the second row of differential forms. The first column contains structure on $X$, and the second one structure on $LX$. }
\label{fig3}
\end{figure}

To conclude, let us discuss the case  $\mathfrak{G}=\mathcal{B}U(1)$. A smooth 2-functor $F: \mathcal{P}_2(X) \to \mathcal{B}\mathcal{B}U(1)$ induces a smooth functor \begin{equation}
\label{80}
\mathrm{tr} \circ \Lambda F \circ \ell : \mathcal{P}_1(LX) \to \mathcal{B}U(1)\text{,}
\end{equation}
since $\Lambda\mathcal{B}\mathcal{B}U(1)=\mathcal{B}U(1)$; the functor $\mathrm{tr}: \mathcal{B}U(1) \to \mathcal{B}U(1)$ just inverts group elements. The image of $F$ under $\mathrm{H}_2$ is hence just a 1-form $\varphi_F\in\Omega^1(X)$, and this 1-form is by Proposition \ref{prop7} just the ordinary fibre integral \begin{equation}
\label{81}
\varphi_F= \int_{S^1} \mathrm{ev}^{*}B\text{.}
\end{equation} 
Let us now interpret the 2-functor $F$ as a trivial abelian gerbe $\mathcal{G}$ with connection over $X$ (see Example \ref{ex5} in Section \ref{sec4_1}),  and the associated functor (\ref{80})  as a trivial principal $U(1)$-bundle $L$ with connection $\varphi_{F}$ over the loop space $LX$, see  Theorem \ref{th5} . Equation (\ref{81}) shows that the line bundle $L$ is the line bundle over the loop space obtained by \emph{transgression} from the gerbe $\mathcal{G}$. Transgression of abelian gerbes as so far been realized in many ways, for example in \cite{gawedzki3, brylinski1, gomi2, schweigert2}, and we have seen here that 
\begin{equation*}
F \mapsto \mathrm{tr} \circ \Lambda F \circ \ell  
\end{equation*}
is just another way to realize transgression. It has one important advantage compared to all  the above previous realizations: it works also for non-abelian gerbes. A further discussion is postponed to the upcoming article  \cite{schreiber2}. 

\begin{appendix}

\tocsection{Appendix}

\def\thesection{A}
\setcounter{theorem}{0}
\setcounter{equation}{0}

\subsection{Basic 2-Category Theory}

\label{app1}

In this article we only consider \emph{strict} 2-categories, 2-functors, inverse 1-morphisms etc., in contrast to the general case. We  only use the qualifier \tql strict\tqr\ in this section and omit it elsewhere.
A general reference is \cite{power1}.
\begin{definition}
\label{def6}
A (small) \emph{2-category} consists of a set of objects, for each pair $(X,Y)$  of objects a set of 1-morphisms denoted $f:X \to Y$ and for each pair $(f,g)$ of 1-morphisms $f,g:X \to Y$ a set  of 2-morphisms denoted $\varphi:f \Rightarrow g$, together with the following structure:
\begin{enumerate}

\item
For every pair $(f,g)$ of 1-morphisms $f: X \to Y$ and $g:Y \to Z$, a 
1-morphism $g \circ f:X \to Y$, called the composition of $f$ and $g$.

\item
For every object $X$, a 1-morphism $\id_X:X \to X$, called the identity 1-morphism of $X$.

\item
For every pair $(\varphi,\psi)$ of 2-morphisms $\varphi: f \Rightarrow g$ and $\psi: g \Rightarrow h$, a 2-morphism $\psi \bullet \varphi: f \Rightarrow h$, called the vertical composition of $\varphi$ and $\psi$.

\item
For every 1-morphism $f$, a 2-morphism $\id_{f}:f \Rightarrow f$, called the identity 2-morphism  of $f$.

\item
For every triple $(X,Y,Z)$ of objects, 1-morphisms $f,f':X \to Y$ and  $g,g':Y \to Z$, and every pair $(\varphi,\psi)$ of 2-morphisms $\varphi:f \Rightarrow f'$ and $\psi:g \Rightarrow g'$, a 2-morphism $\psi \circ \varphi: g \circ f \Rightarrow g' \circ f'$, called the horizontal composition of $\varphi$ and $\psi$.

\end{enumerate}
This structure has to satisfy the following list of axioms:
\begin{enumerate}
\item[(C1)] The composition of 1-morphisms and vertical and horizontal composition of 2-morphisms are associative. 

\item[(C2)] 
The identity 1-morphisms are units with respect to the composition of 1-morphisms,  and identity 2-morphisms are units with respect to vertical composition, i.e.
\begin{equation*}
\varphi \bullet \id_f = \id_g \bullet \varphi
\end{equation*}
for every 2-morphism $\varphi:f \Rightarrow g$. Horizontal composition preserves the identity 2-morphisms in the sense that
\begin{equation*}
\id_g \circ \id_f = \id_{g \circ f}\text{.}
\end{equation*}

\item[(C3)]
Horizontal and vertical  compositions are compatible in the
sense that 
\begin{equation*}
(\psi_1 \bullet \psi_2) \circ (\varphi_1 \bullet \varphi_2) = (\psi_1 \circ\varphi_1)\bullet(\psi_2\circ\varphi_2)
\end{equation*}
whenever these compositions are well-defined. 

\end{enumerate}
\end{definition}

The axioms of a strict 2-category allow to use pasting diagrams for  2-morphisms: every pasting diagram corresponds to a unique 2-morphism. 
In a 2-category, a 2-morphism $\Sigma:\gamma_1 \Rightarrow \gamma_2$ is called
\emph{invertible} or \emph{2-isomorphism} if there exists another 2-morphism $\Sigma^{-1}:
\gamma_2 \Rightarrow \gamma_1$ such that $\Sigma^{-1} \bullet \Sigma = \id_{\gamma_1}$
and $\Sigma \bullet \Sigma^{-1} = \id_{\gamma_2}$. In this case, $\Sigma^{-1}$
is uniquely determined and called  the  \emph{inverse} of $\Sigma$. A 1-morphism
$\gamma: x \to y$ is called \emph{strictly invertible} or \emph{strict 1-isomorphism}, if there exists
another 1-morphism $\bar\gamma: y \to x$ such that $ \id_x = \bar\gamma \circ \gamma$ and $ \gamma \circ
\bar\gamma = \id_y$. A 2-category in which every 1-morphism is strictly invertible is called a \emph{strict 2-groupoid}.

To relate two 2-categories, we use the following definition of a 2-functor,
which is analogous to a functor between categories.

\begin{definition}
Let $S$ and $T$ be two strict 2-categories.
A \emph{strict 2-functor} $F:S \to T$ is an assignment
\begin{equation*}
F \quad:\quad 
\bigon{X}{Y}{f}{g}{\varphi}
\quad\longmapsto\quad
\bigon{F(X)}{F(Y)}{F(f)}{F(g)}{F(\varphi)}
\end{equation*}
such that 
\begin{enumerate}
\item[(F1)]
The vertical composition is respected in the sense that
\begin{equation*}
F(\psi \bullet \varphi) = F(\psi) \bullet F(\varphi)
\quad\text{ and }\quad
F(\id_f) = \id_{F(f)}
\end{equation*}
for all composable 2-morphisms $\varphi$ and $\psi$, and any 1-morphism $f$.
\item[(F2)]
The composition of 1-morphisms is respected in the sense that
\begin{equation*}
F(g) \circ F(f) = F(g \circ h)
\end{equation*}
for all composable 1-morphisms $f$ and $g$, and the horizontal composition of 2-morphisms is respected in the sense that
\begin{equation*}
F(\psi) \circ F(\varphi)=F(\psi \circ \varphi)
\end{equation*}
for all horizontally composable 2-morphisms $\varphi$ and $\psi$.

\end{enumerate}
\end{definition}

To compare 2-functors, we use the notion of a pseudonatural transformation, which
generalizes a natural transformation between  functors.

\begin{definition}
\label{app1_def1}
Let $F_1$ and $F_2$ be two strict   2-functors both  from $S$ to $T$.
A \emph{pseudonatural transformation}
$\rho: F_1 \to F_2$ is an assignment
\begin{equation*}
\rho \quad:\quad
\alxy{ X \ar[r]^{f} & Y}
\quad\longmapsto\quad
\alxydim{@C=1.2cm@R=1.2cm}{F_1(X) \ar[r]^{F_1(f)} \ar[d]_{\rho(X)} & F_1(Y) \ar[d]^{\rho(Y)}
\ar@{=>}[dl]|{\rho(f)} \\ F_2(X) \ar[r]_{F_2(f)} & F_2(Y)\text{,}}
\end{equation*}
of a 2-isomorphism $\rho(f)$ in $T$ to each 1-morphism $f:X \to Y$ in $S$ such that two axioms are satisfied:
\begin{enumerate}
\item[(T1)]
The composition of 1-morphisms in $S$ is respected:
\begin{equation*}
\alxydim{@C=1.2cm@R=1.2cm}{F_1(X) \ar[r]^{F_1(f)} \ar[d]_{\rho(X)} & F_1(Y) \ar[r]^{F_1(g)} \ar[d]|{\rho(Y)}
\ar@{=>}[dl]|{\rho(f)} &F_1(Z) \ar@{=>}[dl]|{\rho(g)} \ar[d]^{\rho(Z)} \\ F_2(X)  \ar[r]_{F_2(f)} & F_2(Y) \ar[r]_{F_2(g)} & F_2(Z)}
=
\alxydim{@C=1.2cm@R=1.2cm}{F_1(X) 
  \ar[r]^{F_1(g \circ f)}="2"  \ar[d]_{\rho(X)} & F_1(Z) \ar[d]^{\rho(Z)}
\ar@{=>}[dl]|{\rho(g \circ f)} \\ F_2(X) \ar[r]_{F_2(g \circ
 f)} & F_2(Z)\text{.}}
\end{equation*}

\item[(T2)]
It is compatible with 2-morphisms:
\begin{equation*}
\alxydim{@C=1.2cm@R=1.2cm}{F_1(X) \ar[r]^{F_1(f)} \ar[d]_{\rho(X)} & F_1(Y) \ar[d]^{\rho(Y)}
\ar@{=>}[dl]|{\rho(f)} \\ F_2(X) \ar@/_3pc/[r]_{F_2(g)}="2"   \ar[r]^{F_2(f)}="1" & F_2(Y) \ar@{=>}"1";"2"|{F_2(\varphi)}}
=
\alxydim{@C=1.2cm@R=1.2cm}{F_1(x) \ar@/^3pc/[r]^{F_1(f)}="1" \ar[r]_{F_1(g)}="2"
\ar@{=>}"1";"2"|{F_1(\varphi)} \ar[d]_{\rho(X)} & F_1(Y) \ar[d]^{\rho(Y)}
\ar@{=>}[dl]|{\rho(g)} \\ F_2(X) \ar[r]_{F_2(g)} & F_2(Y)\text{.}}
\end{equation*}
\end{enumerate}
\end{definition}
It follows that $\rho(\id_X)=\id_{\rho(X)}$ for every object $X$ in $S$.
Pseudonatural
transformations $\rho_1:F_1 \to F_2$ and $\rho_2:F_2 \to F_3$ can naturally
be composed to a pseudonatural transformation $\rho_2 \circ \rho_1: F_1 \to
F_3$:
\begin{equation}
\label{28}
\rho_2 \circ \rho_1 \quad:\quad
\alxy{X \ar[r]^{f} & Y}
\quad\longmapsto\quad
\alxydim{@C=1.2cm@R=1.2cm}{F_1(X) \ar[r]^{F_1(f)} \ar[d]_{\rho_1(X)} & F_1(Y) \ar[d]^{\rho_1(Y)}
\ar@{=>}[dl]|{\rho_1(f)} \\ F_2(X) \ar[d]_{\rho_2(X)} \ar[r]|{F_2(f)} & F_2(Y)  \ar[d]^{\rho_2(Y)}
\ar@{=>}[dl]|{\rho_2(f)} \\ F_3(X) \ar[r]_{F_3(f)} & F_3(Y)\text{.}}
\end{equation}
We need one more definition for situations where we have two pseudonatural
transformations.

\begin{definition}
Let $F_1,F_2:S \to T$ be two strict  2-functors and let $\rho_1,\rho_2: F_1\to F_2$
be pseudonatural transformations. A \emph{modification}
$\mathcal{A}:\rho_1 \Rightarrow \rho_2$ is an assignment
\begin{equation*}
\mathcal{A}\quad:\quad X 
\quad\longmapsto\quad
\bigon{F_1(X)}{F_2(Y)}{\rho_1(X)}{\rho_2(X)}{\mathcal{A}(X)}
\end{equation*}
of a 2-morphism $\mathcal{A}(X)$ in $T$ to any object $X$ in $S$
which satisfies
\begin{equation*}
\alxydim{@C=1.2cm@R=1.2cm}{F_1(X) \ar@/_3.5pc/[d]_{\rho_2(X)}="1" \ar[r]^{F_1(f)} \ar[d]|{\rho_1(X)}="2" \ar@{<=}"1";"2"|-{\mathcal{A}(X)} & F_1(Y) \ar[d]^{\rho_1(y)}
\ar@{=>}[dl]|{\rho_1(f)} \\ F_2(X) \ar[r]_{F_2(f)} & F_2(Y)}
=
\alxydim{@C=1.2cm@R=1.2cm}{F_1(X) \ar[r]^{F_1(f)} \ar[d]_{\rho_2(X)} & F_1(Y) \ar@/^3.5pc/[d]^{\rho_1(X)}="3" \ar[d]|{\rho_2(y)}="4" \ar@{=>}"3";"4"|-{\mathcal{A}(Y)}
\ar@{=>}[dl]|{\rho_2(f)} \\ F_2(X) \ar[r]_{F_2(f)} & F_2(Y)}
\end{equation*}
\end{definition}

Horizontal and vertical compositions of 2-morphisms in $T$ induce  accordant compositions on modifications.

\medskip

For two fixed strict 2-categories $S$ and $T$, we recognize the following structures:
\begin{enumerate}

\item 
For two  strict 2-functors $F_1,F_2:S \to T$, the pseudonatural transformations $\rho:F_1 \to F_2$ together with  modifications and their vertical composition, form a category $\mathrm{Hom}(F_1,F_2)$.

\item
Even more, strict 2-functors from $S$ to $T$, together
with pseudonatural transformations and their modifications, and the assignments
$\circ$ and $\bullet$ as defined above, form a strict 2-category $\mathrm{Funct}(S,T)$.

\end{enumerate}

\medskip

\begin{definition}
Let $S$ and $T$ be strict  2-categories. Strict  2-functors $F:S \to T$ and $G:T \to S$ are called  \emph{isomorphisms of 2-categories}, if $G \circ F=\id_{S}$ and $F \circ G=\id_T$. 
\end{definition}

\subsection{Lie 2-Groups and Smooth Crossed Modules}

\label{app2}

Any strict monoidal category $(\mathfrak{G},\boxtimes,\trivlin)$ defines a 2-category $\mathcal{B} \mathfrak{G}$: it has a single object, the 1-morphisms are the objects of $\mathfrak{G}$ and the 2-morphisms are the morphisms of $\mathfrak{G}$. The horizontal composition is given by the tensor functor $\boxtimes$, and the vertical composition is the composition in $\mathfrak{G}$. The identity 1-morphism of the single object is the tensor unit $1$, and the identity 2-morphism of a 1-morphism $X$ is just the identity morphism $\id_X$ of the object $X$ in $\mathfrak{G}$. The axioms for the 2-category $\mathcal{B}\mathfrak{G}$ follow from the properties of the tensor functor $\boxtimes$. 

In the following, we enhance this construction by two features. First, we assume that $\mathfrak{G}$ is a groupoid and that we have an additional functor $i:\mathfrak{G} \to \mathfrak{G}$ which is an inverse to the tensor functor $\boxtimes$ in the sense that 
\begin{equation}
\label{1}
X \boxtimes i(X) = \trivlin = i(X) \boxtimes X
\quad\text{ and }\quad
f \boxtimes i(f) = \id_{\trivlin} = i(f) \boxtimes f
\end{equation}
for all objects $X$ and all morphisms $f$ in $\mathfrak{G}$.
In this case the 2-category $\mathcal{B}\mathfrak{G}$ is even a  2-groupoid. Secondly, we assume that $\mathfrak{G}$ is a Lie category, and that the functors $\boxtimes$ and $i$ are smooth. Then, $\mathcal{B}\mathfrak{G}$ is a Lie 2-groupoid.

\begin{definition}
A \emph{Lie 2-group} is a strict monoidal Lie category $(\mathfrak{G},\boxtimes,\trivlin)$ together with a smooth functor $i:\mathfrak{G} \to \mathfrak{G}$ such that (\ref{1}) is satisfied. \end{definition}

We denote the  Lie 2-groupoid associated to a Lie 2-group $\mathfrak{G}$ by $\mathcal{B} \mathfrak{G}$.
An important source of Lie 2-groups are smooth crossed modules. 

\begin{definition}
A \emph{smooth crossed module} is a collection $(G,H,t,\alpha)$ of Lie groups $G$ and $H$,  and of a Lie group homomorphism  $t: H \to G$ and a smooth map $\alpha: G \times H \to H$, such that
\begin{enumerate}
\item
$\alpha$ is a left action of $G$ on $H$ by Lie group homomorphisms, i.e.
the smooth map $\alpha_g: H \to H$ defined by $\alpha_g(h):=\alpha(g,h)$ \begin{itemize}
\item[a)]
is a Lie group homomorphism for all $g\in G$.
\item[b)]
satisfies $\alpha_1 = \id_H$ and $\alpha_{gg'} = \alpha_{g} \circ \alpha_{g'}$ for all $g,g'\in G$.
\end{itemize}
\item
$\alpha$ and $t$ are compatible in the following two ways:
\begin{itemize}
\item[a)]
$t(\alpha(g,h))=gt(h)g^{-1}$ for all $g\in G$ and $h\in H$.
\item[b)] 
$\alpha(t(h),x) = hxh^{-1}$ for all $h,x\in H$.
\end{itemize}
\end{enumerate}
\end{definition}

Any smooth crossed module $(G,H,t,\alpha)$ defines a Lie 2-group $(\mathfrak{G},\boxtimes,1,i)$ in the following way. 
\begin{description}
\item[The category $\mathfrak{G}$:] 
We put $\mathrm{Obj}(\mathfrak{G}):= G$ and $\mathrm{Mor}(\mathfrak{G}):= G \ltimes H$, the semi-direct product of $G$ and $H$ defined    by $\alpha$, explicitly
\begin{equation}
\label{2}
(g_{2},h_{2}) \cdot (g_{1},h_{1}) := (g_{2} g_{1},h_2\alpha(g_2,h_1))\text{.}
\end{equation}
An element $(g,h)\in \mathrm{Mor}(\mathfrak{G})$ is considered   as a morphism from $g$ to $t(h)g$. The composition is given by
\begin{equation}
\label{4}
(g',h') \circ (g,h) := (g,h'h)\text{,}
\end{equation} 
where $g'=t(h)g$. It is obviously associative, and the identity morphisms are $\id_g=(g,1)$.
 All these definitions are smooth, so that $\mathfrak{G}$ is a Lie category.

\item[The monoidal structure $(\boxtimes,\trivlin)$:]
The functor $\boxtimes:\mathfrak{G} \times \mathfrak{G} \to \mathfrak{G}$ is defined on objects by $g_2 \boxtimes g_1:=g_2g_1$ and on morphisms by the product (\ref{2}). By axiom 2.a), the morphisms have the correct target.   It respects identity morphisms,
\begin{equation*}
(g_2,1) \boxtimes (g_1,1) = (g_{2}g_{1},1)
\end{equation*}
and by axiom 2.b) the composition
\begin{multline*}
((g_{2}',h_2') \boxtimes (g_1',h_1')) \circ ((g_{2},h_{2}) \boxtimes (g_{1},h_{1})) \\= ((g_{2}',h_2') \circ (g_{2},h_{2})) \boxtimes ((g_1',h_1') \circ (g_{1},h_{1}) )\text{.}
\end{multline*}
It is also strictly associative and the object $\trivlin:=1\in G$ is a left and right unit. 

\item[The functor $i$:] The functor $i:\mathfrak{G} \to \mathfrak{G}$ is defined on objects by $i(g):=g^{-1}$ and on morphisms by $i(g,h) := (g^{-1},\alpha(g^{-1},h^{-1}))$. It respects sources and targets by axiom 2.a), the identities and by axioms 1.a) and 2.b) the composition.  It is also smooth and satisfies the condition (\ref{1}).
\end{description}

Now we have completely defined the Lie 2-group associated to a smooth crossed module. Indeed, it is a well-known fact \cite{brown1}, also see \cite{baez5} for a review, that every Lie 2-group arises -- up to a certain notion of equivalence -- from a smooth crossed module in this way. 

Let us also write down the Lie 2-groupoid $\mathcal{B} \mathfrak{G}$ associated the the Lie 2-group $\mathfrak{G}$ coming from a smooth crossed module $(G,H,t,\alpha)$.
A 2-morphism is  a morphism $(g,h)\in \mathrm{Mor}(\mathfrak{G})$, denoted as
\begin{equation*}
\bigon{\ast}{\ast}{g}{g'}{h}
\end{equation*}
with
\begin{equation}
\label{10}
g'=t(h)g\text{.}
\end{equation}
The ladder equation is also called the \emph{target-matching-condition} for the 2-morphism $(g,h)$.
The vertical composition is
according to (\ref{4})
\begin{equation}
\label{12}
\alxydim{@C=2cm}{\ast \ar[r]|{g'}="2" \ar@/^2.5pc/[r]^{g}="1" \ar@/_2.5pc/[r]_{g''}="3" \ar@{=>}"1";"2"|{h}\ar@{=>}"2";"3"|{h'} & \ast} = \bigon{\ast}{\ast}{g}{g''}{h'h}
\end{equation}
with $g'=t(h)g$ and $g''=t(h')g'=t(h'h)g$, and the horizontal composition is according to (\ref{2})
\begin{equation}
\label{13}
\alxydim{@C=1.2cm}{\ast \ar@/^1.5pc/[r]^{g_{1}}="1" \ar@/_1.5pc/[r]_{g'_{1}}="2" \ar@{=>}"1";"2"|{h_{1}} & \ast \ar@/^1.5pc/[r]^{g_2}="3" \ar@/_1.5pc/[r]_{g'_2}="4" \ar@{=>}"3";"4"|{h_2} & \ast} = \alxydim{@C=3cm}{\ast \ar@/^2.5pc/[r]^{g_2g_1}="1" \ar@/_2.5pc/[r]_{g_2'g_1'}="2" \ar@{=>}"1";"2"|{h_2\alpha(g_2,h_1)} & \ast}
\end{equation}

The construction of Lie 2-groups from smooth crossed modules is convenient to discuss basic examples.

\begin{example}
\normalfont
\label{ex1}
Let $A$ be an abelian Lie group. We define a smooth crossed module by taking $G=\lbrace 1 \rbrace$  the trivial group and $H:= A$. This fixes the maps to $t(a):=1$ and $\alpha(1,a):=a$. All axioms are satisfied in a trivial manner except axiom 2.b), which is satisfied only because $A$ is abelian. The associated Lie 2-group $\mathfrak{G}$ is denoted by $\mathcal{B} A$, and the associated Lie 2-groupoid by $\mathcal{B}\mathcal{B} A$.
\end{example}

\begin{example}
\normalfont
\label{ex2}
Let $G$ be any Lie group. We obtain a smooth crossed module by taking $H := G$, $t=\id$ and $\alpha(g,h):=ghg^{-1}$. The associated Lie 2-group, which also underlies the construction of a geometric realization of $EG$ \cite{segal3} is here denoted by  $\mathcal{E}G$. It can be interpreted as the \emph{inner automorphism 2-group} of $G$ \cite{roberts1}, and its Lie algebra plays an important role in  \cite{sati1}.
\end{example}

Let us briefly exhibit the details of the associated Lie 2-groupoid $\mathcal{EB}G$. It has one objects, and the set of 1-morphisms is $G$ with the usual composition $g_2 \circ g_1 := g_2 g_1$. Between every pair $(g,g')$ of 1-morphisms there is a unique 2-morphism
\begin{equation*}
\bigon{\ast}{\ast}{g}{g'}{h}
\end{equation*}
determined by $h:=g'g^{-1}$.

\begin{example}
\normalfont
\label{ex4}
Let $H$ be a connected Lie group. The group of Lie group automorphisms of $H$ is again a Lie group $G:=\mathrm{Aut}(H)$ \cite{onishchik1}. Together with  $t(h)(x):=hxh^{-1}$ and $\alpha(\varphi,h):= \varphi(h)$, we have defined a smooth crossed module whose associated Lie 2-group $\mathfrak{G}$ is denoted by $\mathrm{AUT}(H)$. 
\end{example}

\subsection{Proof of Lemma \ref{lem4}}

\label{app3}

In this section we show that the map 
\begin{equation*}
k_{A,B}:BX \to H
\end{equation*}
defined for the construction of a smooth 2-functor from two differential forms $A\in \Omega^1(x,\mathfrak{g})$ and $B \in \Omega^2(X,\mathfrak{h})$, only depends on the thin homotopy class of a bigon $\Sigma\in BX$. 

\begin{figure}[h]
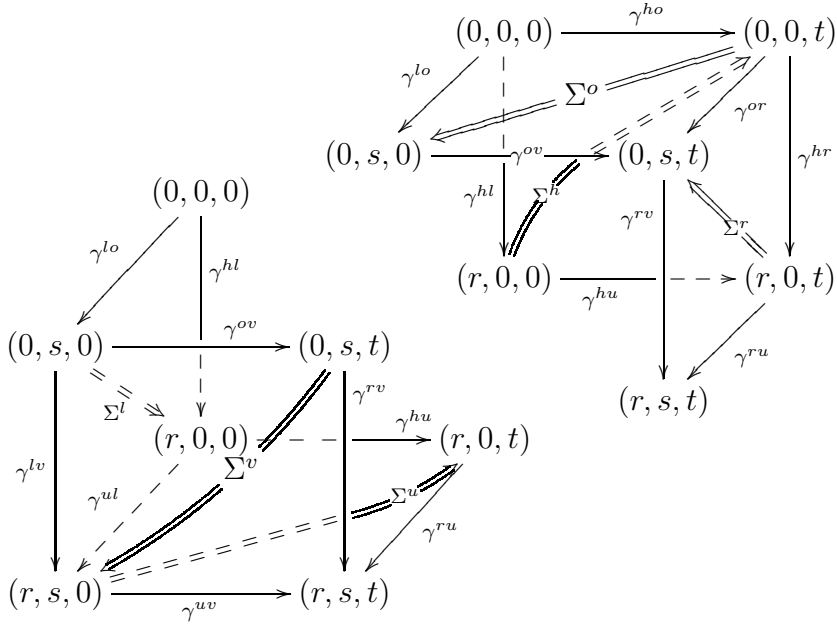

\begin{equation*}
\alxydim{@C=0.4cm@R=0.6cm}{&&&\\&&&\\&&&\\& (0,0,0) \ar@{-}[dd]^{\gamma^{hl}} \ar[ddl]_{\gamma^{lo}} && \\ &&&\\ (0,s,0) \ar@{==>}[dr]_{\Sigma^{l}} \ar[rr]^>>>>>>>{\gamma^{ov}} \ar[ddd]_{\gamma^{lv}} & \ar@{-->}[d] & (0,s,t) \ar@{=>}@/^0.8pc/[dddll]|>>>>>>>>>>>>>>>>*+{\Sigma^v} \ar[ddd]^<<<<{\gamma^{rv}} & \\ & (r,0,0) \ar@{-->}[ddl]_{\gamma^{ul}} \ar@{--}[r] & \ar[r]^{\gamma^{hu}} &(r,0,t) \ar[ddl]^{\gamma^{ru}} \\&& \ar@{=>}@/_0.35pc/[ur]|<<<<<{\Sigma^{u}}&\\ (r,s,0) \ar@{==}[urr] \ar[rr]_{\gamma^{uv}} && (r,s,t) &}
\hspace{-3cm}
\alxydim{@R=1cm@C=0.2cm}{& (0,0,0) \ar@{--}[d] \ar[rrr]^{\gamma^{ho}} \ar[ld]_{\gamma^{lo}} &&& (0,0,t) \ar@{=>}[dllll]|*+{\Sigma^{o}} \ar[dl]^{\gamma^{or}} \ar[dd]^{\gamma^{hr}} \\ (0,s,0) \ar[rrr]|<<<<<<<<<<<<<{\gamma^{ov}} & \ar[d]_<<<<<{\gamma^{hl}} & \ar@{==>}[urr] & (0,s,t)  \ar[dd]_<<<<<<{\gamma^{rv}} & \\ & (r,0,0) \ar@{=}@/^0.4pc/[ur]|>>>>>{\Sigma^{h}} \ar@{-}[rr]_<<<<<<{\gamma^{hu}} && \ar@{-->}[r] & (r,0,t) \ar@{=>}[ul]|<<<<<<{\Sigma^{r}} \ar[dl] ^{\gamma^{ru}} \\ &&&(r,s,t) \\ &&&&\\ &&&&}
\end{equation*}
\caption{The unit cube viewed as two bigons: $\Lambda_1: \mu => \nu$ on the right hand side, and $\Lambda_2:\nu \rightarrow \mu$ on the left hand side.}
\label{fig2}
\end{figure}

We first start with a general homotopy $h: [0,1]\times [0,1]^2 \to X$ between two bigons $\Sigma_1$ and $\Sigma_2$, i.e $h$ has the properties from Definition \ref{def1} except condition 2a) which constrains the rank of its differential. We shall represent the surface of the unit cube $[0,1]^3$ on which $h$ is defined as a bigon in $\R^3$. For this purpose, we define two paths $\mu(r,s,t)$ and $\nu(r,s,t)$ in $\R^3$ going from $0\in \R^3$ to $(r,s,t)$. With the notation introduced in Figure \ref{fig2}
these paths are
\begin{equation*}
\mu(r,s,t) := \gamma^{ru} \circ \gamma^{hu} \circ \gamma^{hl}
\quad \text{ and }\quad
\nu(r,s,t) := \gamma^{vr} \circ \gamma^{ov} \circ \gamma^{lo}\text{.}
\end{equation*}
Between these paths we have two bigons $\Lambda_1(r,s,t): \mu(r,s,t) \Rightarrow \nu (r,s,t)$ and $\Lambda_2(r,s,t):\nu(r,s,t) \Rightarrow \mu(r,s,t)$ defined by
\begin{equation}
\label{decomp1}
\Lambda_1 := (\id_{ \gamma^{rv}} \ast \Sigma^o) \bullet (\Sigma^r \ast \id_{\gamma^{ho}}) \bullet (\id_{\gamma^{ru}} \ast \Sigma^h)
\end{equation}
and 
\begin{equation}
\label{decomp2}
\Lambda_2 := (\Sigma^u \ast \id_{\gamma^{hl}}) \bullet (\id_{\gamma^{uv}} \ast \Sigma^l) \bullet (\Sigma^v \ast \id_{\gamma^{lo}}) 
\end{equation}
The vertical composition $\Lambda(r,s,t) := \Lambda_2(r,s,t) \bullet \Lambda_1(r,s,t)$ is then a bigon whose image is the surface of the cube. Notice that the two bigons $\Sigma_1$ and $\Sigma_2$ we started with can be found on the top and on the bottom of the unit cube, i.e.  $\Sigma_1 = h_{*}\Sigma^o$ and $\Sigma_2 = h_{*}(\Sigma^u)^{-1}$.

We evaluate the map $k_{A,B}$ on the bigon $h_{*}(\Lambda(r,s,t))$, defining a smooth function $u:[0,1]^3 \to H$. Since $k_{A,B}$ is by Lemma \ref{lem6} compatible with the vertical composition and the auxiliary horizontal composition $\ast$ we get from \erf{decomp1} and \erf{decomp2}
\begin{eqnarray}
u(r,s,t) &=& k_{A,B}(h_{*}\Sigma^u) \cdot \alpha(F_A(h_{*}\gamma^{uv}),k_{A,B}(h_{*}\Sigma^{l}))\nonumber\\&&\hspace{1cm} \cdot k_{A,B}(h_{*}\Sigma^v)  \cdot \alpha(F_A(h_{*}\gamma^{rv}), k_{A,B}(h_{*}\Sigma^o))\nonumber\\&&  \hspace{1cm}\cdot k_{A,B}(h_{*}\Sigma^r)  \cdot \alpha(F_A(h_{*}\gamma^{ru}), k_{A,B}(h_{*}\Sigma^h))\text{.}
\label{46}
\end{eqnarray}
On the right hand side we have omitted the arguments $(r,s,t)$ for simplicity. 
In the following we use a bigon $\Lambda_{r_0}(r,s,t)$ that is shifted by $r_0$ along the $r$-axis with respect to the bigon $\Lambda(r,s,t)$ from Figure \ref{fig2}, to which it reduces for $r_0=0$.
 Accordingly, we have a smooth function $u_{r_0}: [0,1]^3 \to H$ additionally depending on the shift $r_0$. In the same way, we have a smooth function $u_{r_0,s_0,t_0}: [0,1]^3 \to H$ associated to a bigon $\Lambda_{r_0,s_0,t_0}$ that is additionally shifted by $s_0$ and $t_0$ along the $s$-axis and the $t$-axis, respectively.
\begin{lemma}
The smooth function $u_{r_0}: [0,1]^3 \to H$ has the following properties:
\begin{enumerate}
\item[(a)]
$u_{0}(1,1,1) = k_{A,B}(\Sigma_2)^{-1}  \cdot  k_{A,B}(\Sigma_1)$.

\item[(b)]

$u_{r_0}(r,1,1) = u_{r_0+ r'}(r-r',1,1) \cdot u_{r_0}(r',1,1)$.

\item[(c)]
$\displaystyle
u_{r_0}(r,s+\sigma,1) = H_{1}(r_0,r,s) \cdot u_{r_0,s,0}(r,\sigma,1) \cdot H_2(r_{0},r,s)
$.

\item[(d)]

$\displaystyle\frac{1}{3}\left . \frac{\partial^2}{\partial r \partial \sigma} \right|_0 \frac{\partial}{\partial t} u_{r_0,s,0}(r,\sigma,t) =(h^{*}K)_{(r_0,s,t)} \left( \frac{\partial}{\partial r},\frac{\partial}{\partial s},\frac{\partial}{\partial t} \right)$.
\end{enumerate}
In (c), $H_1$ and $H_2$ are certain $H$-valued smooth functions that do not depend on $\sigma$.
In (d), the 3-form $K \in \Omega^3(X,\mathfrak{h})$ is given by $K:=\mathrm{d}B + \alpha_{*}(A \wedge B)$. 
\end{lemma}

\proof
Condition 1 of Definition \ref{def1} for the homotopy $h$ implies the vanishing of the group elements $k_{A,B}(h_{*}\Sigma^r),k_{A,B}(h_{*}\Sigma^{l})$ and $F_A(h_{*}\gamma^{rv}),F_A(h_{*}\gamma^{ru})$ in the product (\ref{46}), so that
\begin{equation}
\label{16}
k_{A,B}(h_{*}\Lambda(1,1,1)) = k_{A,B}(\Sigma_2)^{-1}  \cdot k_{A,B}(h_{*}\Sigma^v)  \cdot k_{A,B}(\Sigma_1)  \cdot k_{A,B}(h_{*}\Sigma^h)\text{.}
\end{equation}
By condition 2b), the bigons $h_{*}\Sigma^v: \gamma_1'\Rightarrow \gamma_2'$ and $h_{*}\Sigma^h:\gamma_1 \Rightarrow \gamma_2$ are thin homotopies between paths, i.e. the rank of their differentials is less or equal to 1. Accordingly, $\mathcal{A}_{h_{*}\Sigma^{v}}=\mathcal{A}_{h_{*}\Sigma^h}=0$ and hence $k_{A,B}(h_{*}\Sigma^v)=k_{A,B}(h_{*}\Sigma^h)=1$.
Now, assertion (a) follows from \erf{16}. The same vanishing arguments show that
\begin{equation}
\label{20}
u_{x}(y,1,1) = k_{A,B}(h_{*}\Sigma^u_{x+y}(1,1)) \cdot k_{A,B}(h_{*}\Sigma^o_{x}(1,1))\text{.}
\end{equation}
Using  formula \erf{20} by setting $(x,y)$  to $(r_0,r)$, $(r_0+r',r-r')$ and $(r_0,r')$, respectively, shows assertion (b). Still the same arguments show that
\begin{multline*}
u_{r_0,0,0}(r,s+\sigma,1) = k_{A,B}h_{*}\Sigma^u_{r_0+r,0,0}(s,1) \cdot u_{r_0,s,0}(r,\sigma,1) \\\cdot k_{A,B}h_{*}\Sigma^v_{r_0,s,0}(r,1) \cdot k_{A,B}h_{*}\Sigma^o_{r_0,0,0}(s,1)\text{,}
\end{multline*}
which is assertion (c) identifying 
\begin{eqnarray*}
H_{1}(r_0,r,s) &=& k_{A,B}h_{*}\Sigma^u_{r_0+r,0,0}(s,1)
\\
H_2(r_{0},r,s) &=&  k_{A,B}h_{*}\Sigma^v_{r_0,s,0}(r,1) \cdot k_{A,B}h_{*}\Sigma^o_{r_0,0,0}(s,1). \end{eqnarray*}
For (d), we write down  formula \erf{46} for $u_{r_0,s,0}(r,\sigma,t+\tau)$ and decompose bigons of length $t + \tau$ into two bigons of length $t$ and $\tau$, respectively. This gives
\begin{eqnarray*}
&&\hspace{-0.8cm}u_{r_0,s,0}(r,\sigma,t+\tau)\\&& = k_{A,B}h_{*} \Sigma^u_{r_0+r,s,t}(\sigma,\tau)\cdot \alpha(k_{A,B}h_{*}\gamma^{uv}_{r_0+r,s_0+\sigma,t}(\tau), k_{A,B}h_{*}\Sigma^u_{r_0+r,s,0}(\sigma,t))
\\&&\hspace{1cm}\cdot\alpha(k_{A,B}h_{*}\gamma^{uv}_{r_0+r,s+\sigma,t}(\tau),k_{A,B}h_{*}\Sigma^v_{r_0,s+\sigma,0}(r,t))\cdot k_{A,B}h_{*}\Sigma^v_{r_0,s+\sigma,t}(r,\tau)
\\&&\hspace{1cm}\cdot \alpha \left (k_{A,B}h_{*}\gamma^{rv}_{r_0,s+\sigma,t+\tau}(r),\alpha(k_{A,B}h_{*}\gamma^{or}_{r_0,s+\sigma,t}(\tau),k_{A,B}h_{*}\Sigma^o_{r_0,s,0}(\sigma,t))
\right . \\&& \hspace{7cm}\left .\cdot k_{A,B}h_{*}\Sigma^o_{r_0,s,t}(\sigma,\tau) \right)\\&&\hspace{1cm}\cdot k_{A,B}h_{*}\Sigma^l_{r_0,s,t+\tau}(r,\sigma) \cdot \alpha(k_{A,B}h_{*}\gamma^{ru}_{r_0+r,s,t+\tau}(\sigma),k_{A,B}h_{*}\Sigma^h_{r_0,s,t}(r,\tau)\\&&\hspace{1cm}\cdot \alpha(k_{A,B}h_{*}\gamma^{hu}_{r_0+r,s,t}(\tau),k_{A,B}h_{*}\Sigma^h_{r_0,s,0}(r,t)))\text{.}
\end{eqnarray*}
Now we take the derivative of this expression by the three variables $r$, $\sigma$ and $\tau$, evaluate at zero and use Proposition \ref{prop4} in order to identify second derivatives of $k_{A,B}$ with the 2-form $B$. This gives
\begin{eqnarray*}
&&\hspace{-1.5cm}\left . \frac{\partial^3}{\partial r \partial \sigma\partial\tau} \right|_0 u_{r_0,s,0}(r,\sigma,t+\tau) \\[\smallskipamount]&=&3 \mathrm{d}B_{h(r_0,s,t)}(v_r,v_s,v_t) + \alpha_{*} (A_{r_0,s,t}(v_r),B_{r_0,s,t}(v_s,v_t))
\\&&+ \alpha_{*} (A_{r_0,s,t}(v_s),B_{r_0,s,t}(v_t,v_r))+ \alpha_{*} (A_{r_0,s,t}(v_t),B_{r_0,s,t}(v_r,v_s))\text{.}
\end{eqnarray*}
Using the antisymmetry of $B$, this shows (d).
\endofproof

\begin{remark}
The 3-form $K=\mathrm{d}B + \alpha_{*}(B \wedge A)$ that drops out in (d) has to be interpreted as the curvature of the connection $(A,B)$ on a trivial, (non-abelian) gerbe, see Section \ref{sec4_1}.
\end{remark}

Now, if the homotopy $h$ is thin, i.e. satisfies condition 2a) of Definition \ref{def1}, we have by (d)
\begin{equation*}
\left . \frac{\partial^2}{\partial r\partial \sigma}\right |_0 u_{r_0,s,0}(r,\sigma,1) =  \int_{0}^1\mathrm{d}t 
\left . \frac{\partial^2}{\partial r \partial \sigma} \right|_0 \frac{\partial}{\partial t} u_{r_0,s,0}(r,\sigma,t) =0\text{.}
\end{equation*}
Performing this trick once more, we obtain
\begin{eqnarray*}
\left . \frac{\partial}{\partial r} \right|_0 u_{r_0}(r,1,1)&=& \int_0^1
\mathrm{d}s \left .  \frac{\partial}{\partial r} \right|_0\frac{\partial }{\partial s} u_{r_0}(r,s,1)
\\&\stackrel{\mathrm{(c)}}{=}& \int_0^1
\mathrm{d}s H_{1}(r_0,r,s) \cdot \left \lbrace  \left .  \frac{\partial^2}{\partial r\partial \sigma} \right|_0  u_{r_0,s,0}(r,\sigma,1) \right \rbrace \cdot H_2(r_{0},r,s)
\\&=& 0
\end{eqnarray*}
The multiplicative property (b) transfers this result to all values of $r$,
\begin{equation*}
\left . \frac{\partial}{\partial r}\right |_{r_0} u_{0}(r,1,1) = \left . \frac{\partial}{\partial r}\right |_0 u_{r_0}(r,1,1) \cdot u_{0}(r_0,1,1) =0\text{.}
\end{equation*}
This means that the function $u_0(r,1,1)$ is constant, and thus determined by its value at $r=0$,
\begin{equation*}
1 = u_0(0,1,1) = u_0(1,1,1)\stackrel{\mathrm{(a)}}{=}k_{A,B}(\Sigma_2)^{-1}  \cdot  k_{A,B}(\Sigma_1)\text{.}
\end{equation*}
This shows that $k_{A,B}$ takes the same values on thin homotopic bigons $\Sigma_1$ and $\Sigma_2$.

\end{appendix}

\tocsection{Table of Notations}

\newcommand{\notation}[3]{
\noindent
\begin{minipage}[t]{0.15\textwidth}#1\end{minipage}
\begin{minipage}[t]{0.68\textwidth}#2\vspace{0.3cm}\end{minipage}\hfill
\begin{minipage}[t]{0.105\textwidth}Page \pageref{#3}\end{minipage}
}

\notation{$\mathrm{AUT}(H)$}{the automorphism 2-group of a Lie group $H$.}{ex4}

\notation{$\mathcal{B}G$}{the Lie groupoid with one object associated to a Lie group $G$.}{not:oneobject}

\notation{$\mathcal{B}\mathfrak{G}$}{the Lie 2-groupoid with one object associated to a Lie 2-group $\mathfrak{G}$.}{sec3_1}

\notation{$BX$}{the diffeological space of bigons in $X$.}{not:bx}

\notation{$B^2X$}{the diffeological space of thin homotopy classes of bigons in $X$.}{def1}

\notation{$\mathcal{D}$}{the 2-functor that extracts differential forms from smooth 2-functors.}{not:d}

\notation{$D^{\infty}$}{the category of diffeological spaces.}{not:dinfty}

\notation{$\mathcal{E}G$}{the inner automorphism 2-group associated to a Lie group $G$.}{ex2}

\notation{$\mathrm{Funct}^{\infty}$}{the category of smooth functors between Lie categories.}{not:smoothfun}

\notation{$LX$}{the loop space $D^\infty(S^1,X)$ of a diffeological space $X$.}{sec4_2}

\notation{$\Lambda$}{the functor that makes a category out of a 2-category.}{not:lambda}

\notation{$\ell$}{the functor that regards a path in the loop space of $X$ as a bigon in $X$.}{not:ell}

\notation{$PX$}{the diffeological space of smooth paths (with sitting instants) in $X$.}{not:smoothpaths}

\notation{$P^1X$}{the diffeological space of thin homotopy classes of paths in $X$.}{not:thinhomotopyclasses}

\notation{$\mathcal{P}_1(X)$}{the path groupoid of $X$.}{not:pathgroupoid}

\notation{$\mathcal{P}_2(X)$}{the path 2-groupoid of $X$.}{not:p2g}

\notation{$\mathcal{P}$}{the 2-functor that integrates differential forms to smooth 2-functors.}{sec3_2}

\notation{$\diffco{G}{1}{X}$}{the category of $G$-connections on $X$.}{not:gconn}

\notation{$\diffco{\mathfrak{G}}{2}{X}$}{the 2-category of $\mathfrak{G}$-connections on $X$.}{def8}

\kobib

\end{document}